\documentclass[10pt]{amsart}
\usepackage{amsmath}
\usepackage{amsthm}
\usepackage{latexsym}
\usepackage{amssymb}
\usepackage{mathrsfs}
\usepackage{tikz}
\usepackage{wasysym}

\def\Dslash{\textrm{\DH}}
\allowdisplaybreaks
\DeclareMathOperator{\Char}{Char}

\DeclareMathOperator{\re}{Re}

\DeclareMathOperator{\ad}{ad}

\DeclareMathOperator{\sgn}{sign}

\DeclareMathOperator{\supp}{supp}

\DeclareMathOperator{\ord}{ord}

\newtheorem{theorem}{Theorem}
\newtheorem{proposition}{Proposition}
\newtheorem{lemma}{Lemma}
\newtheorem{corollary}{Corollary}
\newtheorem{definition}{Definition}

\newtheorem{remark}{Remark}

\def\jpopn#1#2{%
  \mathopen{%
    \setbox0=\hbox{$#1\langle$}%
    \setbox2=\hbox{%
            {\hbox{$#1\langle$}}%
            \kern -.6\wd0\box0%
    }%
            \box2%
  }%
}

\def\jpcls#1#2{%
  \mathclose{%
    \setbox0=\hbox{$#1\rangle$}%
    \setbox2=\hbox{%
            {\hbox{$#1\rangle$}}%
            \kern -.6\wd0\box0%
    }%
            \box2%
  }%
}

\def\tp{\ {}^{t}\kern -3pt}

\def\R {{\mathbb{R}}}
\def\N {{\mathbb{N}}}
\def\C {{\mathbb{C}}}

\def\Q {{\mathbb{Q}}}
\renewcommand{\thetheorem}{\thesection.\arabic{theorem}}
\renewcommand{\theproposition}{\thesection.\arabic{proposition}}
\renewcommand{\thelemma}{\thesection.\arabic{lemma}}
\renewcommand{\thedefinition}{\thesection.\arabic{definition}}
\renewcommand{\thecorollary}{\thesection.\arabic{corollary}}
\renewcommand{\theequation}{\thesection.\arabic{equation}}
\renewcommand{\theremark}{\thesection.\arabic{remark}}

\AtBeginDocument{
  \renewcommand{\setminus}{\mathbin{\backslash}}%
}
\def\uu#1{$\breve{\textrm{#1}}$}
\def\vv#1{$\check{\textrm{#1}}$}

\begin{document}
\renewcommand{\kappa}{\ensuremath{\varkappa}}
\renewcommand{\phi}{\ensuremath{\varphi}}
\renewcommand{\epsilon}{\ensuremath{\varepsilon}}

%
\title[Optimal Gevrey Regularity]{Optimal Gevrey
  Regularity for Certain Sums of Squares in Two Variables}
\author{Antonio Bove}
\address{Dipartimento di Matematica, Universit\`a
di Bologna, Piazza di Porta San Donato 5, Bologna,
Italy}
\email{bove@bo.infn.it}
\author{Marco Mughetti}
\address{Dipartimento di Matematica, Universit\`a
di Bologna, Piazza di Porta San Donato 5, Bologna,
Italy}
\email{marco.mughetti@unibo.it}
\date{\today}
\begin{abstract}
For $ q $, $ a $ integers such that $ a \geq 1 $, $ 1 <
q $, $ (x, y) \in U $, $ U $ a neighborhood of the origin in $ \R^{2}
$, we consider the operator
$$ 
D_{x}^{2} + x^{2(q-1)} D_{y}^{2} + y^{2a} D_{y}^{2} .
$$
Slightly modifying the method of proof of \cite{monom} we can see that
it is Gevrey $ s_{0} $  hypoelliptic,
where $ s_{0}^{-1} = 1 - a^{-1} (q - 1) q^{-1} $. Here we show that
this value is optimal, i.e. that there are solutions to $ P u = f $
with $ f $ more regular than $ G^{s_{0}} $ that are not better than
Gevrey $ s_{0} $.

The above operator reduces to the M\'etivier operator
(\cite{metivier81}) when $ a = 1 $,
 $ q = 2 $. We give a description of the
characteristic manifold of the operator and of its relation with the
Treves conjecture on the real analytic regularity for sums of
squares. 
\end{abstract}
\maketitle
\tableofcontents

\section{Introduction and Statement of the Result}

For $ q $, $ a $ integers such that $ a \geq 1 $, $ 1 <
q $, $ (x, y) \in \R^{2} $, we consider the operator
\begin{equation}
\label{eq:P}
P(x, y, D_{x}, D_{y}) =   
D_{x}^{2} + x^{2(q-1)} D_{y}^{2} + y^{2a} D_{y}^{2} .
\end{equation}
Here $ D_{x} = 1/\sqrt{-1} \  \partial_{x} $. 

The operator $ P $ generalizes the operator
$$ 
P_{M} (x, y, D_{x}, D_{y}) =   
D_{x}^{2} + x^{2} D_{y}^{2} + y^{2} D_{y}^{2},
$$
defined by M\'etivier and for which M\'etivier, in an article in the
Comptes Rendus de l'Acad\'emie des Sciences in 1981, stated that it is not
real analytic hypoelliptic, \cite{metivier81}. Actually in
\cite{metivier81} it was stated that $ P_{M} $ is Gevrey 2
hypoelliptic and not better, meaning that it is not Gevrey-$ s $
hypoelliptic for any $ s $, $ 1 \leq s < 2 $.

For the sake of completeness we give the definition of $
G^{s}(\Omega) $---the class of all Gevrey-$ s $ functions in the open set
$ \Omega $:  
\begin{definition}
\label{def:gevrey}
If $ \Omega \subset \R^{n} $ is an open set we say that the function $
u $ belongs to the Gevrey class of order $ s \geq 1 $, $ G^{s}(\Omega)
$, if $ u \in C^{\infty}(\Omega) $ and for every compact set $ K
\subset \Omega $ there exists a positive constant $ C_{K} $ such that
$$
\sup_{K} | \partial_{x}^{\alpha} u(x) | \leq C_{K}^{|\alpha| + 1}
\alpha!^{s} .
$$
Note that $ G^{1}(\Omega) $ is the class of the real analytic
functions in $ \Omega $.
\end{definition}
We point out explicitly that the characteristic set of $ P_{M} $ is
the non symplectic submanifold of $ T^{*}\R^{2}\setminus\{0\} $ given by
$$ 
\Char(P_{M}) = \{ (x, y; \xi, \eta) \ | \ x = y = 0, \xi = 0, \eta
\neq 0 \}.
$$
We have
$$ 
\Char(P) = \{ (x, y; \xi, \eta) \ | \ x = y = 0, \xi = 0, \eta
\neq 0 \} = \Char(P_{M}) .
$$
The problem of the analytic or Gevrey regularity of the operator
\eqref{eq:P} is part of the broader problem of the analytic regularity
for sums of squares of vector fields with analytic
coefficients. Technically speaking the operator \eqref{eq:P} is not a
sum of squares, but it is easy to see that it shares the same
properties with the sum of squares
$$ 
D_{x}^{2} + \left(x^{q-1} D_{y}\right)^{2} + \left(y^{a}
  D_{y}\right)^{2} .
$$
There has been a fair amount of literature on the problem in the
eighties and ninenties and we refer to the paper \cite{bm-surv} for a
survey of the main results.

We focus on the conjecture stated by Treves (see \cite{Treves}, \cite{btreves} and
\cite{trevespienza} for its statement) in 1999, characterizing the
analytic hypoellipticity of sums of squares.

The Treves conjecture defines a stratification of $
\Char(P) $ such that the strata satisfy the following properties
(everything is meant to be microlocal near a fixed point in $ \Char(P) $)
\begin{itemize}
\item[(i)]{}
  Real analytic submanifolds of $ \Char(P) $.
\item[(ii)]{}
  The symplectic form $ \sigma = d\xi \wedge dx $ restricted to each
  stratum has constant rank.
\item[(iii)]{}
  For each stratum there exists an integer $ \nu \in \N $ such that
  all the Poisson brackets of the symbols of the vector fields of
  length $ \leq \nu $ are identically zero, but there is a Poisson
  bracket of length $ \nu+1 $ which is not zero. By length of a
  Poisson bracket we mean the number of vector fields of which we take
  an iterated bracket.
\end{itemize}
According to the conjecture an operator is analytic hypoelliptic if
and only if every stratum of its characteristic variety is a
symplectic manifold.

For the operator $ P $ the stratification is made up of two half
lines, $ \Sigma_{+} \cup \Sigma_{-} $, where
$$ 
\Sigma_{\pm} = \{ (x, y; \xi, \eta) \ | \ x = y = 0, \xi = 0, \pm \eta
> 0 \} .
$$
$ \Sigma_{\pm} $ is an isotropic submanifold of $
T^{*}\R^{2}\setminus\{0\} $ and, as such, it coincides with its
Hamilton leaf along the half-fiber $ \pm \eta >0 $.

Thus the Treves stratification is not symplectic, which suggests that
$ P $ is not analytic hypoelliptic. We recall that if a manifold is
not symplectic there is always a foliation, called the Hamilton
foliation, whose leaves can or cannot be transverse to the fibers of
the cotangent bundle.

We point out that the operator in \eqref{eq:P} actually 
is one of the simplest model operator exhibiting a Hamilton leaf non
transverse to the fibers.

By adapting the technique of \cite{monom} one can prove that $ P $ is Gevrey $
s_{0} $ hypoelliptic, with
$$ 
s_{0}^{-1} = 1 - \frac{1}{a} \frac{q - 1}{q} .
$$
In the present paper we want to show that $ s_{0} $ is an optimal
index, i.e.
\begin{theorem}
\label{th:1}
The operator $ P $ is not Gevrey-$ s $ hypoelliptic for any $ s $ such
that $ 1 \leq s < s_{0} $. Hence $P$ is not analytic hypoelliptic.
\end{theorem}
We point out that the optimality proof in the non tranverse case is
technically much more involved than that in the transverse case and,
to our knowledge, in the non transverse case only M\'etivier
paper \cite{metivier81} is available. 

In the remaining part of this section we want to give a brief
exposition about our motivation to study the operator $ P $ and set it
in a more general context in the framework of the problem of the real
analytic regularity of the solutions to sums of squares
operators, of which $ P $ is a two variable example.

An essential tool in M\'etivier's paper, \cite{metivier81}, is the
expansion of a solution as a 
linear combination of the eigenfunctions of the corresponding harmonic
oscillator, $ D_{x}^{2} + x^{2} $. It is known that the eigenfunctions
of the harmonic oscillator
are a basis in $ L^{2}(\R) $, are rapidly decreasing at infinity and
satisfy finite recurrence relations allowing us to express e.g. the
derivative, or the multiplication by $ x $, of an eigenfunction as a
combination of two different eigenfunctions. These recurrence relations
are essential in the M\'etivier approach.

When the degree of the potential is larger than two, i.e. for the
anharmonic oscillator $ D_{x}^{2} + x^{2(q-1)} $, $ q > 2 $, Gundersen
in \cite{gundersen} proved the
\begin{theorem}[\cite{gundersen}]
\label{th:gundersen}
Consider the equation in $ \C $
$$ 
w''(z) + ( \lambda - p(z)) w(z) = 0, \qquad \lambda \in \R, \ z \in \C . 
$$
Here $ p $ denotes the polynomial $ p(z) = a_{2m}z^{2m} + a_{2m-2}
z^{2m-2} + \ldots + a_{2} z^{2} $. Assume that $ a_{i} \geq 0 $ for
every $ i $ and $ a_{2m} > 0 $.

Denote by $ (\psi_{n}(z))_{n \in \N} $ a set of solutions that is a
complete orthonormal basis in $ L^{2}(\R) $.

Then if $ m $ is even and $ g $ is a non zero polynomial we have that
$ g(z)$ $ \cdot \psi_{\ell}^{(k)}(z) $ is not a finite linear combination of
the $ \psi_{n} $, for $ k \in \N \cup \{0\} $, $ \ell \in \N $ fixed
and $ \deg(g) \geq 1 $ if $ k = 0 $.
\end{theorem}
Thus for $ q $ odd there is no finite recurrence relation for the
anharmonic oscillator.

For an even $ q $, Bender and Wang, \cite{bw}, showed that the
eigenfunctions for the operator
\begin{equation}
\label{eq:bender}
- u''(x) + x^{2N+2} u(x) = E x^{N}u(x) , \qquad N = -1, 0, 1, 2, \ldots,
\end{equation}
do have finite recurrence relations using confluent hypergeometric
functions. We refer to Chinni, \cite{chinni}, for a result in this sense.

Our proof of Theorem \ref{th:1}, as well as M\'etivier's approach, shares the
general pathway of constructing an asymptotic formal solution with all
the optimality proofs in the transverse case. This is done first at
the formal level and then substantiated by
introducing suitable cutoff functions to turn a formal solution into a
true solution of an equation with a better right hand side. 

We shall show that the Gevrey regularity
of $P$ is strictly related to the functions satisfying the equation 
\begin{equation}
\label{eq:oscill}
- u''(x) + x^{2(q-1)} u(x) = \lambda u(x) ,
\end{equation}
i.e. the eigenfunctions of the operator
\begin{equation}
  \label{eq:Q}
Q = D^{2} + x^{2(q-1)}.
\end{equation}
However we do not require any recurrence relations among
the eigenfunctions. 

\bigskip

From a more geometric point of view we point out that even though the
Treves conjecture has been shown not to hold for
$ n \geq 4 $ (we refer to \cite{abm}, \cite{bm} for a proof and to
\cite{bm-j} for a case that might suggest that strata are still the
object to be defined) there are no counterexamples when $ n = 2, 3 $.

We observe that in the proofs of \cite{abm}, \cite{bm} there is a
non symplectic submanifold of the characteristic manifold that is not
identified by the Treves conjecture stratification procedure and plays
a crucial role in carrying the non real analytic wave front set. This
submanifold has a foliation on the base of the cotangent bundle,
i.e. there is a Gevrey (or analytic) wave front set propagation in the
space variables.

When $ n=3 $ this may no longer occur due to a dimensional constraint,
but we think that e.g. the following example
\begin{equation}
\label{eq:nonT-dim-3}
D_{1}^{2} + x_{1}^{2(r-1)} D_{y}^{2} + D_{2}^{2} + x_{2}^{2(q-1)}
D_{y}^{2} + x_{2}^{2(p-1)} y^{2a} D_{y}^{2} ,
\end{equation}
where $ a $, $ r $, $ p $, $ q \in \N $, $ 1 < r < p < q $, should be
a candidate to violate the Treves conjecture in dimension $ 3 $, since
the characteristic manifold has a foliation whose leaves are the $
\eta $ fibers of the cotangent bundle that, in our opinion, should
carry analytic wave front set.


To our knowledge no proof is known either of its analytic regularity
or of its non-analytic hypoellipticity. We hope that the same
techniques we are using for the operator $ P $ above can be suitably
modified to prove that the above example is not analytic hypoelliptic.

While for $ n \geq 3 $ there is no conjecture about the stratification
to use in order to characterize analytic hypoellipticity, when $ n=2 $
the Treves conjecture seems to describe accurately the geometry of the
characteristic variety
and we think that it might actually be true (see \cite{monom} again).
Furthermore the operator studied in this paper is a microlocal model
for a class of operators in two variables.

In this perspective, proving Theorem \ref{th:1} seems a necessary step
to accomplish a proof of the conjecture when $ n=2 $.
However no proof is available up to now. 

Finally a few words about the proof of the theorem. We just sketch the
idea of the proof with a simplified notation.

First we construct a formal solution of $ P(A(u)) = 0 $, of the form
$$ 
A(u)(x, y) = \int_{0}^{+\infty} e^{i y \rho^{s_{0}}}  u(x
\rho^{\frac{s_{0}}{q}}, \rho) d\rho ,
$$
where
$$ 
u(x, \rho) = \sum_{j\geq0} u_{j}(x, \rho).
$$
This is accomplished by observing that if $ P(A(u)) = 0 $, then $ u $
must satisfy to
$$ 
\sum_{j=0}^{2a} \frac{1}{\rho^{j}} P_{j}(x, D_{x},
D_{\rho}) u(x, \rho) = 0 ,
$$
where the $ P_{j} $ are differential operators of order $ 2a $. This
is done in section 2. In section 3 we compute the $ u_{j} $ as
solutions to countably many PDEs of the form
$$ 
P_{0} u_{j} = - \sum_{k=1}^{\min\{j, 2a\}} \frac{1}{\rho^{k}} P_{k}
u_{j-k} ,
$$
where, essentially, $ P_{0} = D_{x}^{2} + x^{2(q-1)} + D_{\rho}^{2a}
$ on the half plane $ \R_{x} \times$ $ ]0, +\infty[ $. In order to
compute $ u_{j} $ we have to compute an inverse of $ P_{0} $ in the
half plane $ \rho > 0 $. This is accomplished by separating variables,
using the eigenfunctions of the anharmonic oscillator $ D_{x}^{2} +
x^{2(q-1)} $, and solving an ODE in $ \rho $ for $ \rho > 0 $.
It can be seen that when solving the above mentioned ODE corresponding
to higher eigenstates of the anharmonic oscillator we gain a
decreasing rate of the form $ \mathscr{O}(\rho^{-j} e^{-c\rho}) $. When
it comes to the fundamental eigenstate the situation is more involved
and, to obtain a similar decreasing rate, one has to decouple the
equation for $ u_{j} $, using the projector onto the fundamental
eigenstate, $ \pi $, and its complement $ 1 - \pi $ (see \eqref{eq:tr1},
\eqref{eq:tr2}).

In sections 4, 5, 6, 7 we derive weighted $ L^{2} $-type estimates for $
\pi u_{j} $, $ (1-\pi)u_{j} $ and hence for $ u_{j} $ (Theorem
\ref{th:est}). The reason why 
we resorted to $ L^{2} $, and not $ L^{\infty} $, estimates is that we exploit the
orthonormality of the eigenfunctions of the anharmonic oscillator,
since there are no finite recurrence relations.

In section 8 the pointwise estimates for the $ u_{j} $ are derived via
the Sobolev immersion theorem (Theorem \ref{th:pest-t}). 

In order to turn the formal solution $ A(u) $ of the equation $ P A(u)
= 0$ into a true solution, in section 9 we replace $ u $ with $ v $ where
$$ 
v(x, \rho) = \sum_{j\geq 0} \psi_{j}(\rho) u_{j}(x, \rho),
$$
where $ \psi_{j} $ denotes a suitable cutoff function whose precise
definition is given in Lemma \ref{lemma:psij}. Then $ A(v) $ solves an
equation of the form $ P A(v) = f $, for a function $ f $ in a
Beurling class of order $ s_{0} $ (see Definition \ref{def:Bs} for the
definition of such classes). From Theorem 3.1 of M\'etivier,
\cite{metivier80}, arguing by contradiction, if $ P $ is Gevrey
hypoelliptic of order less than $ s_{0} $, then $ A(v) $ belongs to
the Beurling class of order $ s_{0} $.

The purpose of section 10 is to obtain a contradiction to the fact
that $ A(v) $ is in a Beurling class of order $ s_{0} $ by using the
structure of $ A(v) $ given above. This is done by comparing the
growth rate of its $ y $-Fourier transform as a member of the Beurling class
of order $ s_{0} $ (see Lemma \ref{lemma:fourier}) as well as a
function given by the above expression (see Lemma \ref{lemma:below}),
when $ x $ is frozen at the origin. As a technical detail we mention
that the use of the Fourier transform forces us to replace the
Beurling class of order $ s_{0} $ with $ L^{2}(\R_{y}) $ intersected
with the global Beurling class of order $ s_{0} $ (see definition
\ref{def:Bs}). To show that actually $ A(v) $ belongs to the latter
class we use Lemma \ref{lemma:yd}. 

Finally we gathered in the Appendixes some results needed in the main
body of the paper. Appendixes A and B present some basically known
facts and fix the notation in our setting. Appendix C and D present
some estimates that are essential for sections 5 and 6.
We postponed them in order to allow the reader to better follow
the deployment of the proof in those sections.

\section{A Formal Solution}

\setcounter{equation}{0}
\setcounter{theorem}{0}
\setcounter{proposition}{0}
\setcounter{lemma}{0}
\setcounter{corollary}{0}
\setcounter{definition}{0}
\setcounter{remark}{0}

First we look for a formal solution to $ P v = 0 $ by taking $ v $ of
the form
\begin{equation}
\label{eq:Au}
A(u)(x, y) = \int_{0}^{+\infty} e^{i y \rho^{s_{0}}} \rho^{r} u(x
\rho^{\frac{s_{0}}{q}}, \rho) d\rho ,
\end{equation}
where
\begin{equation}
\label{eq:s0}
\frac{1}{s_{0}} = 1 - \frac{1}{a} \ \frac{q - 1}{q} ,
\end{equation}
$ r $ is a complex number to be chosen later and $ u $ denotes a
smooth function defined in $ \R \times \R $, with $ \supp u \subset \{
\rho > 0 \}$, and rapidly decreasing for $ \rho \to + \infty $.

We have
\begin{equation}
\label{eq:D2A}
D_{x}^{2} A(u)(x, y) = \int_{0}^{+\infty} e^{i y \rho^{s_{0}}} \rho^{r
+ 2 \frac{s_{0}}{q}} (D_{x}^{2} u)(\rho^{\frac{s_{0}}{q}} x, \rho) \ d\rho.
\end{equation}
\begin{equation}
\label{eq:x2(q-1)}
x^{2(q-1)} D_{y}^{2} A(u)(x, y) 
= \int_{0}^{+\infty} e^{i y \rho^{s_{0}}} \rho^{r
+ 2 \frac{s_{0}}{q}} \Big [x^{2(q-1)} u(x, \rho)\Big]_{x = x \rho^{s_{0}/q}} \ d\rho.
\end{equation}
Here the notation $ x = x \rho^{s_{0}/q} $ means that the variable $ x
$ inside the square brackets has to be replaced by the r.h.s. of the
equation. 

Finally
\begin{equation}
\label{eq:x2(p-1)}
y^{2a} D_{y}^{2} A(u)(x, y) 
= \int_{0}^{+\infty} e^{i y \rho^{s_{0}}} \rho^{r
+ 2 s_{0}} y^{2a} \Big [ u(x,
\rho)\Big]_{x = x \rho^{s_{0}/q}} \  d\rho.
\end{equation}
Since
$$ 
y e^{i y \rho^{s_{0}}} = \frac{1}{i s_{0} \rho^{s_{0} - 1}}
\partial_{\rho} e^{i y \rho^{s_{0}}} ,
$$
integrating by parts we may rewrite the above expression as
\begin{multline*}
y^{2a} D_{y}^{2} A(u)(x, y) 
\\
= \int_{0}^{+\infty} e^{i y \rho^{s_{0}}} \left( - \partial_{\rho}
  \frac{1}{i s_{0} \rho^{s_{0} - 1}}\right)^{2a}   \rho^{r
+ 2 s_{0} \left(\frac{q-1}{q} + \frac{1}{q}\right)} 
\Big [ u(x, \rho) \Big]_{x = x\rho^{s_{0}/q}} \ d\rho.
\end{multline*}
Remark that
\begin{equation}
  \label{eq:gamma0}
\left(\partial_{\rho}
  \frac{1}{i s_{0} \rho^{s_{0} - 1}}\right)^{2a} = \sum_{h=0}^{2a}
\gamma_{2a, h} \frac{1}{\rho^{2a s_{0} - h}} \partial_{\rho}^{h},
\end{equation}
where the $ \gamma_{2a, h} $ are constants satisfying
\begin{equation}
  \label{eq:gamma}
| \gamma_{2a, h} | \leq C_{\gamma}^{2a+h} (2a-h)! , \qquad \gamma_{2a,
2a} = \left(\frac{i}{s_{0}}\right)^{2a}, \ \gamma_{2a, 0} = 1. 
\end{equation}
Hence, writing for the sake of simplicity, $ \gamma_{h} $ instead of $
\gamma_{2a, h} $, we obtain
\begin{multline*}
y^{2a} D_{y}^{2} A(u)(x, y) 
\\
= \int_{0}^{+\infty} e^{i y \rho^{s_{0}}} 
\sum_{h=0}^{2a} \gamma_{h} \frac{1}{\rho^{2as_{0} - h}}
\partial_{\rho}^{h} \left(
 \rho^{r
+ 2 s_{0} \left(\frac{q-1}{q} + \frac{1}{q}\right)} 
\Big [u(x, \rho)\Big]_{x = x\rho^{s_{0}/q}} \right) \ d\rho
\\
=
\int_{0}^{+\infty} e^{i y \rho^{s_{0}}} 
\sum_{h=0}^{2a} \sum_{\alpha=0}^{h} \binom{h}{\alpha} \gamma_{h}
\left(r + 2 s_{0} \left( \frac{q-1}{q} + \frac{1}{q} \right)\right)_{h-\alpha}
\\
\cdot
\left.
\rho^{-2as_{0} + h + r + 2s_{0} \left( \frac{q-1}{q} +
    \frac{1}{q}\right) - h + \alpha}
  \partial_{\rho}^{\alpha} 
\Big [ u(x, \rho)\Big]_{x = x\rho^{s_{0}/q}} \right) \ d\rho ,
\end{multline*}
where $ (\lambda)_{\beta} $ denotes the Pochhammer symbol, defined by
\begin{equation}
\label{eq:pochhammer}
(\lambda)_{\beta} = \lambda (\lambda -1) \cdots (\lambda - \beta + 1),
\qquad (\lambda)_{0} = 1, \quad \lambda \in \C.
\end{equation}
Moreover, since
$$ 
a (1 - s_{0}) + s_{0} \frac{q-1}{q} = a\left( 1 - s_{0} \left( 1 -
    \frac{1}{a} \ \frac{q-1}{q}\right) \right) = 0,
$$
we obtain that
\begin{multline*}
y^{2a} D_{y}^{2} A(u)(x, y) 
= \int_{0}^{+\infty} e^{i y \rho^{s_{0}}}
\sum_{h=0}^{2a} \sum_{\alpha=0}^{h} \binom{h}{\alpha} \gamma_{h}
\left(r + 2 \frac{s_{0}}{q} + 2 s_{0} \frac{q-1}{q} \right)_{h-\alpha}
\\
\cdot
\rho^{r + 2 \frac{s_{0}}{q} + \alpha - 2a}
 \partial_{\rho}^{\alpha} 
\Big [u(x, \rho)\Big]_{x = x\rho^{s_{0}/q}}  \ d\rho .
\end{multline*}
Since
\begin{equation}
  \label{eq:drhobin}
\partial_{\rho}^{\alpha} v(\rho^{\theta}x, \rho) = \left [ \left(
  \frac{\theta}{\rho} x \partial_{x} + \partial_{\rho} \right)^{\alpha} v(x,
\rho) \right]_{x = x \rho^{\theta}} ,
\end{equation}
we have
\begin{multline*}
y^{2a} D_{y}^{2} A(u)(x, y) 
=
\sum_{h=0}^{2a} \sum_{\alpha=0}^{h} \binom{h}{\alpha} \gamma_{h}
\left(r + 2 \frac{s_{0}}{q} + 2 s_{0} \frac{q-1}{q} \right)_{h-\alpha}
\\
\cdot
\int_{0}^{+\infty} e^{i y \rho^{s_{0}}} \rho^{r + 2 \frac{s_{0}}{q} +
  \alpha - 2a}
\left[
  \left(\frac{s_{0}}{q \rho} 
    x \partial_{x} + \partial_{\rho}
\right)^{\alpha} u(x, \rho) \right]_{x = x\rho^{s_{0}/q}} \
d\rho .
\end{multline*}
Because of the identity
$$ 
\left(\frac{a}{\rho} + \partial_{\rho}\right)^{\alpha} = \rho^{-a}
\partial_{\rho}^{\alpha} \rho^{\alpha} =
\sum_{k=0}^{\alpha}\binom{\alpha}{k} (a)_{k} \rho^{-k}
\partial_{\rho}^{\alpha - k} ,
$$
we deduce that
\begin{multline}
  \label{eq:x2(p-1):2}
y^{2a} D_{y}^{2} A(u)(x, y) 
=
\sum_{h=0}^{2a} \sum_{\alpha=0}^{h} \binom{h}{\alpha} \gamma_{h}
\left(r + 2 \frac{s_{0}}{q} + 2 s_{0} \frac{q-1}{q} \right)_{h-\alpha}
\\
\cdot
\int_{0}^{+\infty} e^{i y \rho^{s_{0}}} \rho^{r + 2 \frac{s_{0}}{q} +
  \alpha - 2a}
\left[
\sum_{k=0}^{\alpha} \binom{\alpha}{k} \left( \frac{s_{0}}{q} x
  \partial_{x} \right)_{k} \rho^{-k} \partial_{\rho}^{\alpha-k}
 u(x, \rho) \right]_{x = x \rho^{s_{0}/q}} d\rho, 
\end{multline}
where, in analogy with \eqref{eq:pochhammer}, we used the notation
$$ 
\left( \frac{s_{0}}{q} x \partial_{x} \right)_{k} = \frac{s_{0}}{q} x
\partial_{x} \left( \frac{s_{0}}{q} x \partial_{x} - 1 \right) \cdots
\left( \frac{s_{0}}{q} x \partial_{x} - k + 1 \right) .
$$
An inspection of \eqref{eq:x2(p-1):2} readily gives that the
differential operator above is a polynomial of degree $ k $ in $ x
\partial_{x} $ with uniformly bounded coefficients. We therefore may
write it in the form
$$ 
\left( \frac{s_{0}}{q} x \partial_{x} \right)_{k} = p_{k}(x
\partial_{x}) = \sum_{j=1}^{k} b_{k,j} (x \partial_{x})^{j},
$$
where
$$ 
b_{k, k} = \left( \frac{s_{0}}{q} \right)^{k}, \qquad b_{k, 0} =
\frac{s_{0}}{q} (-1)^{k-1}(k-1)!
$$
Hence
\begin{multline}
  \label{eq:x2(p-1):3}
y^{2a} D_{y}^{2} A(u)(x, y) 
=
\sum_{h=0}^{2a} \sum_{\alpha=0}^{h} \binom{h}{\alpha} \gamma_{h}
\left(r + 2 \frac{s_{0}}{q} + 2 s_{0} \frac{q-1}{q} \right)_{h-\alpha}
\\
\cdot
\int_{0}^{+\infty} e^{i y \rho^{s_{0}}} \rho^{r + 2 \frac{s_{0}}{q} +
  \alpha - 2a}
\left[
\sum_{k=0}^{\alpha} \binom{\alpha}{k} p_{k}(x\partial_{x}) \rho^{-k}
\partial_{\rho}^{\alpha-k} u(x, \rho) \right]_{x = x \rho^{s_{0}/q}} d\rho, 
\end{multline}
where $ p_{0}(x\partial_{x}) = 1 $ by convention.

From \eqref{eq:D2A}, \eqref{eq:x2(q-1)} and the above we may then
write
\begin{multline}
\label{eq:PA}
P A(u) (x, y) = \int_{0}^{+\infty} e^{i y \rho^{s_{0}}} \rho^{r + 2
  \frac{s_{0}}{q}} \left[ D_{x}^{2} + x^{2(q-1)} \right.
\\
  + \sum_{h=0}^{2a} \sum_{\alpha=0}^{h} \sum_{k=0}^{\alpha}  \binom{h}{\alpha} \gamma_{h}
\left(r + 2 \frac{s_{0}}{q} + 2 s_{0} \frac{q-1}{q} \right)_{h-\alpha}
\binom{\alpha}{k}
\\
\left.
p_{k}(x\partial_{x})
\rho^{\alpha - k -2a} \partial_{\rho}^{\alpha-k}  u(x, \rho)
\right]_{x = x \rho^{s_{0}/q}} d\rho. 
\end{multline}
Let us consider the sums on the r.h.s. of the above expression. First
we call $ \ell = \alpha - k $ and then rewrite the sums setting $ \ell
= 2a - j$. We obtain
\begin{multline*}
P A(u) (x, y) = \int_{0}^{+\infty} e^{i y \rho^{s_{0}}} \rho^{r + 2
  \frac{s_{0}}{q}} \left[ D_{x}^{2} + x^{2(q-1)} \right.
\\
+ \sum_{h=0}^{2a} \sum_{\alpha=0}^{h} \sum_{j=2a-\alpha}^{2a}  \binom{h}{\alpha} \gamma_{h}
\left(r + 2 \frac{s_{0}}{q} + 2 s_{0} \frac{q-1}{q} \right)_{h-\alpha}
\binom{\alpha}{2a-j} 
\\
\left.
p_{\alpha+j-2a}(x\partial_{x})
\rho^{-j} \partial_{\rho}^{2a-j}  u(x, \rho)
\right]_{x = x \rho^{s_{0}/q}} d\rho. 
\end{multline*}
We may then interchange the sums in $ \alpha $ and $ j $ and after
that interchange the sums in $ h $ and $ j $ to get
\begin{multline}
\label{eq:PA:2}
P A(u) (x, y) = \int_{0}^{+\infty} e^{i y \rho^{s_{0}}} \rho^{r + 2
  \frac{s_{0}}{q}} \left[ D_{x}^{2} + x^{2(q-1)} \right.
\\
+ \sum_{j=0}^{2a} \rho^{-j} \left( \sum_{h=2a-j}^{2a} \sum_{\alpha=2a-j}^{h}  \binom{h}{\alpha} \gamma_{h}
\left(r + 2 \frac{s_{0}}{q} + 2 s_{0} \frac{q-1}{q} \right)_{h-\alpha}
\right.
\\
\left. \left.
\vphantom{\sum_{h=2a-j}^{2a}}
\binom{\alpha}{2a-j}
p_{\alpha+j-2a}(x\partial_{x})
\right) \partial_{\rho}^{2a-j}  u(x, \rho)
\right]_{x = x \rho^{s_{0}/q}} d\rho. 
\end{multline}
Denote by $ P_{0}(x, D_{x}, D_{\rho}) $ the differential operator
corresponding to $ j=0 $:
\begin{equation}
\label{eq:P0}
P_{0}(x, D_{x}, D_{\rho}) = D_{x}^{2} + x^{2(q-1)} + s_{0}^{-2a}
  D_{\rho}^{2a} ,
\end{equation}
the differential operator induced by $ P $---modulo some
normalization---on the $ \eta $ fibers. Furthermore set, for $ j \geq
1 $, 
\begin{align}
\label{eq:Pj}
P_{j}(x, \partial_{x}, \partial_{\rho})& = 
\sum_{h=2a-j}^{2a} \sum_{\alpha=2a-j}^{h}  \binom{h}{\alpha} \gamma_{h}
\left(r + 2 \frac{s_{0}}{q} + 2 s_{0} \frac{q-1}{q} \right)_{h-\alpha}
\notag
  \\
& \phantom{=}\cdot  
\binom{\alpha}{2a-j}
p_{\alpha+j-2a}(x\partial_{x})
\partial_{\rho}^{2a-j}
\\[5pt]
& =   \tilde{P}_{j}(x\partial_{x}) \partial_{\rho}^{2a-j} ,
\notag                                           
\end{align}
with $ \ord(\tilde{P}_{j}) = j $. 
We also have, see also \eqref{eq:gamma0},
\begin{multline}
\label{eq:P1}
P_{1}(x, \partial_{x}, \partial_{\rho}) =
\\
\left( 2a
  \gamma_{2a} \left(p_{1}(x\partial_{x}) + r +
    \frac{2s_{0}}{q} + 2 s_{0} \frac{q-1}{q}\right)
  +\gamma_{2a-1} \right) \partial_{\rho}^{2a-1}. 
\end{multline}
To keep the notation simple we write
\begin{equation}
\label{eq:Pjtilda}
\tilde{P}_{j}(x \partial_{x}) = \sum_{\ell=0}^{j} p_{j\ell} (x
\partial_{x})^{\ell} ,
\end{equation}
for $ j = 1, \ldots, 2a $ and suitable numbers $ p_{j\ell} $.

Equation \eqref{eq:PA:2} can thus be written as
\begin{multline}
\label{eq:PA:3}
P A(u) (x, y)
\\
= \int_{0}^{+\infty} e^{i y \rho^{s_{0}}} \rho^{r + 2
  \frac{s_{0}}{q}} \left[\sum_{j=0}^{2a} \frac{1}{\rho^{j}} P_{j}(x, D_{x},
D_{\rho}) u(x, \rho)\right]_{x = x \rho^{s_{0} /q}} d\rho .
\end{multline}
\begin{remark}
\label{rem:0}
The whole computation has been carried out for the mo\-del in
\eqref{eq:P}. We might have allowed variable coefficients in
\eqref{eq:P}, at least to a certain extent, which would have resulted
in an infinite sum in \eqref{eq:PA:3}. The rest of the proof should
proceed essentially with minor modifications. We however preferred to stick to the
proposed model for a greater clarity in the exposition, with a very
tiny degree of generality lost.
\end{remark}
Our next step is to formally solve the equation
\begin{equation}
\label{eq:transp}
\sum_{j=0}^{2a} \frac{1}{\rho^{j}} P_{j}(x, D_{x},
D_{\rho}) u(x, \rho) = 0.
\end{equation}

\section{Computing the Formal Solution}

\setcounter{equation}{0}
\setcounter{theorem}{0}
\setcounter{proposition}{0}
\setcounter{lemma}{0}
\setcounter{corollary}{0}
\setcounter{definition}{0}
\setcounter{remark}{0}

We start this section by discussing how to solve the prototype
equation
\begin{equation}
\label{eq:P0u=f}
P_{0}(x, D_{x}, D_{\rho}) u(x, \rho) = f(x, \rho),
\end{equation}
where $ P_{0} $ is given by \eqref{eq:P0}.

As a preliminary result we consider the kernel of the operator
\begin{equation}
\label{eq:Qlam}
Q_{\mu}(x, D) = D^{2} + x^{2(q-1)} - \mu ,
\end{equation}
$ \mu $ being a real parameter.

The following proposition is well known (see e.g. \cite{berezin}.)
\begin{proposition}
\label{prop:ker-Q}
There exist countably many positive numbers, $ \mu_{j} $,
$ \mu_{j+1} > \mu_{j} $, $ j \geq 0 $,  such that
$$ 
\ker Q_{\mu_{j}} \neq \{ 0 \},
$$
and actually $ \dim \ker Q_{\mu_{j}} = 1 $. Here when we write $
\ker Q_{\mu_{j}} $ we mean the kernel of the operator $
Q_{\mu_{j}} $ as an unbounded operator in $ L^{2}(\R) $.
\end{proposition}
Let us denote by $ \mu_{j} $, $ \phi_{j}(x) $, $ j \geq 0 $,
the eigenvalues and the eigenfunctions of the operator $ Q $, defined
in \eqref{eq:Q}, constructed
in the above proposition.

They are an orthonormal basis in $ L^{2}(\R) $. 

\bigskip

Consider now the equation
$$ 
\left( D^{2} + x^{2(q-1)} +  s_{0}^{-2a}
  D_{\rho}^{2a}\right) u = f,
$$
with $ f \in L^{2} $. We want to find an expression for $ u \in
L^{2} $ intersected with the natural domain of $ P_{0} $. We note
explicitly that the operator $ P_{0} $ has no tempered distributions
in its kernel if we consider it in the $ (x, \rho) $-plane. We are
considering it in a half plane, and hence we can find non trivial
distributions in its kernel.

Write
$$ 
u(x, \rho) = \sum_{k\geq 0} u_{(k)}(\rho) \phi_{k}(x),
$$
where $ u_{(k)}(\rho) = \langle u, \phi_{k} \rangle
$. Here $ \langle \cdot , \cdot \rangle $ denotes the (complex) scalar product
in $ L^{2}(\R) $, in the only variable $ x $.

Hence, with an analogous expansion of $ f $,
\begin{multline*}
P_{0}(x, D, D_{\rho}) u = \sum_{k \geq 0} \left( u_{(k)} Q \phi_{k} +
   \phi_{k} s_{0}^{-2a} D_{\rho}^{2a} u_{(k)} \right)
\\
= \sum_{k \geq 0} \phi_{k} \left( u_{(k)} \mu_{k} +
  s_{0}^{-2a} D_{\rho}^{2a} u_{(k)} \right) = \sum_{k\geq 0}
f_{k} \phi_{k} ,
\end{multline*}
where $ f_{k} = \langle f, \phi_{k} \rangle $.

Identifying the coefficients of $ \phi_{k} $ in the
last equality above we may find the $ u_{(k)} $ as the solution of the
differential equations
\begin{equation}
\label{eq:ode-k}
\left( \partial_{\rho}^{2a} + (-1)^{a} s_{0}^{2a} \mu_{k} \right)
u_{(k)} = (-1)^{a} s_{0}^{2a} f_{k}.
\end{equation}
Let us denote again by $ f_{k} $ the r.h.s. of \eqref{eq:ode-k}, the
difference being just a multiplicative constant.

For $ j = 1, \ldots, 2a $, denote by $ \mu_{k j} $ the $ 2a $-roots
of $ (-1)^{a+1} s_{0}^{2a} \mu_{k} $.

Assume now that $ k \geq 1 $.

Applying the result of Appendix A, we find that
\begin{equation}
\label{eq:u_{k}}
u_{(k)}(\rho) = \sum_{j=1}^{2a} A_{kj} I_{kj}(f_{k})(\rho) ,
\end{equation}
where the $ A_{kj} $ are defined as in \eqref{eq:Al} and
\begin{multline}
\label{eq:Ikj}
I_{kj}(f_{k})(\rho)
\\
= -  \sgn\left( \re \mu_{kj}\right) \int_{\R} e^{\mu_{kj}
  (\rho - \sigma)} H\left(-\sgn\left(\re \mu_{kj}\right) (\rho -
  \sigma)\right)
\\
\cdot H(\sigma -R) f_{k}(\sigma) d\sigma ,
\end{multline}
where $ H $ denotes the Heaviside function, $ R $ is a positive
number that can be chosen depending on some parameter to be
precised. We note that 
\eqref{eq:Ikj} defines a function solving \eqref{eq:ode-k} for $ \rho
> R $. 

Consider the case $ k = 0 $.

We define $ u_{(0)} $ in a slightly different way. The motivation for
such a distinction will be clear in the subsequent sections.

Consider  $ \mu_{0} $, the smallest of the values $ \mu_{0} < \mu_{1} <
\cdots $ defined in Proposition \ref{prop:ker-Q}, and
denote by $ \mu_{0i} $, $ i = 1, \ldots, 2a $, the $ 2a
$-roots of
$$
(-1)^{a+1} s_{0}^{2a}  \mu_{0}.
$$
Then define $\tilde{\mu}_{0} $ by
\begin{align}
\label{eq:mitlda0}
  \tilde{\mu}_{0} &= \re \mu_{0 i^{*}} , \text{ where } 
  \\
  &\phantom{ = }
\ \  \mu_{0 i^{*}}
\text{ is a $ 2a $-root of $ (-1)^{a+1} s_{0}^{2a} \mu_{0} $ with
    maximum negative }%
    \notag
  \\
 &\phantom{ = }
  \  \text{ real part.}
   \notag
\end{align}
We remark explicitly that, if $ a > 1 $, there are always two (complex
conjugate) roots, $ \mu_{0i_{1}} $, $ \mu_{0i_{2}} $ of maximal
negative real part, since we are taking even roots. Of course the
definition of $ \tilde{\mu}_{0} $ is independent of the choice of the
root. 

The reason of the above choice for $ \tilde{\mu}_{0} $, which will
have important implications in the sequel, is that we have better
decay rates for the components along higher eigenfunctions than that
of the fundamental eigenfunction.

Then we define
\begin{align}
\label{eq:u0}
u_{(0)}(\rho) &= - \sum_{\begin{subarray}{c} i \in \{1, \ldots, 2a\} \\
  \re \mu_{0i} > 0 \\ \text{ or } \re \mu_{0i} =
  \tilde{\mu}_{0} \end{subarray}} A_{0i} \int_{\rho}^{+\infty}
  e^{\mu_{0i} (\rho - \sigma)} H(\sigma - R) f_{0}(\sigma) d\sigma \\
  &\phantom{=}\
+ \sum_{\begin{subarray}{c} i \in \{1, \ldots, 2a\} \\
    \re \mu_{0i} < \tilde{\mu}_{0}  \end{subarray}} A_{0i}
  \int_{R}^{\rho}  e^{\mu_{0i} (\rho - \sigma)} H(\sigma - R) f_{0}(\sigma)
  d\sigma .
  \notag
\end{align}
We point out that if $ f_{0}(\sigma) = \mathscr{O}(\sigma^{-1-\delta}
e^{\tilde{\mu}_{0}\sigma}) $, with $ \delta > 0 $, the integrals where
$ \re \mu_{0i} = \tilde{\mu}_{0} $ are well defined.

It is convenient to harmonize the notation for \eqref{eq:u_{k}},
\eqref{eq:Ikj} and \eqref{eq:u0}. To this end we point out that, for $
k > 0 $, $ \re \mu_{ki} > 0 $ implies that $ \re \mu_{ki} -
\tilde{\mu}_{0} > 0 $, and that if $ \re \mu_{ki} < 0 $ then $ \re
\mu_{ki} - \tilde{\mu}_{0} < \re \mu_{0i} - \tilde{\mu}_{0} \leq 0 $,
due to the fact that the sequence $ (\mu_{k})_{k \in \N \cup \{0\}} $
is strictly increasing. Since these inequalities for $ k > 0 $ can be
perturbed it is clear that there exists a small positive number, say $
\epsilon_{\mu} $, such that for any $ k \geq 0 $ we may write that 
$$ 
u_{(k)}(\rho) = \sum_{j=1}^{2a} A_{kj} I_{kj}(f_{k})(\rho) ,
$$
with
\begin{multline}
\label{eq:Ikjgl}
I_{kj}(f_{k})(\rho)
= -  \sgn\left( \re \mu_{kj} - \tilde{\mu}_{0} + \epsilon_{\mu} \right) \int_{\R} e^{\mu_{kj}
  (\rho - \sigma)}
\\[5pt]
\cdot
H\left(-\sgn\left(\re \mu_{kj} - \tilde{\mu}_{0} +
    \epsilon_{\mu}\right) (\rho -  \sigma)\right) 
\cdot H(\sigma - R) f_{k}(\sigma) d\sigma ,
\end{multline}
where the $ A_{0i} $ are chosen according to the above prescription. 

\bigskip

Thus we obtain an expression for $ u $:
\begin{align}
\label{eq:P0inv}
  u(x, \rho) & = P_{0}^{-1} (f)
\\
           &  = \sum_{k\geq 0} \left(
               \sum_{j=1}^{2a} A_{kj} I_{kj}(f_{k})(\rho) \right) \phi_{k}(x)
\notag  \\
             & = \sum_{k\geq 0} E_{k}(f_{k})(\rho) \phi_{k}(x), \notag
\end{align}
where for the sake of brevity we used the notation
\begin{equation}
\label{eq:Ek}
E_{k}(f_{k})(\rho) =  \sum_{j=1}^{2a} A_{kj} I_{kj}(f_{k})(\rho) .
\end{equation}
In order to find a formal solution to
$$ 
P(x, y, D_{x}, D_{y}) A(u) = 0 ,
$$
we look for a function $ u $ of the form
\begin{equation}
\label{eq:uj}
u(x, \rho) = \sum_{j \geq 0} u_{j}(x, \rho).
\end{equation}
such that, when we plug it into \eqref{eq:transp}, it gives
\begin{equation}
\label{eq:transp1}
\sum_{k=0}^{2a} \frac{1}{\rho^{k}} P_{k}(x, D_{x},
D_{\rho}) u(x, \rho)
= \sum_{k=0}^{2a} \frac{1}{\rho^{k}} P_{k}(x, D_{x},
D_{\rho}) \sum_{j \geq 0} u_{j}(x, \rho) = 0.
\end{equation}
It seems natural to split \eqref{eq:transp1} according to
\begin{equation}
  \label{eq:nattransp}
P_{0} u_{j} = - \sum_{k=1}^{\min\{j, 2a\}} \frac{1}{\rho^{k}} P_{k}
u_{j-k}, \qquad j \geq 0 .
\end{equation}
(The sum is understood to be zero if its upper index is zero).

Consider the first equation to be solved:
\begin{equation}
\label{eq:P0u0}
P_{0}(x, D, D_{\rho}) u_{0}(x, \rho) = 0.
\end{equation}
Arguing as we just did, we may take, see \eqref{eq:mitlda0}, 
\begin{equation}
\label{eq:u0-}
u_{0}(x, \rho) = \phi_{0}(x) e^{\mu_{0i^{*}} \rho}.
\end{equation}
We would like to compute the $ u_{j} $, $ j \geq 1 $, in such a way that
\eqref{eq:nattransp} is satisfied and that $ u_{j} =
\mathscr{O}(\rho^{-j} e^{\tilde{\mu}_{0} \rho}) $ (actually, for
technical reasons, we prove in Theorem \ref{th:pest-t} a slightly
weaker estimate). 

Solving \eqref{eq:nattransp}, using the spectral decomposition $
\phi_{k}(x) $, by \eqref{eq:ode-k}, boils down to solving an ode of
the form
$$ 
(\partial_{\rho}^{2a} + (-1)^{a} s_{0}^{2a} \mu_{k})w = g_{k} ,
$$
where $ g_{k} = \mathscr{O}(\rho^{-j} e^{\tilde{\mu}_{0} \rho}) $ and
$ \tilde{\mu}_{0} $ is defined in \eqref{eq:mitlda0}. If $ k \geq 1 $
it is easy to see that $ w = \mathscr{O}(\rho^{-j} e^{\tilde{\mu}_{0}
  \rho})  $, since $ \tilde{\mu}_{0} $ differs from the real part of
the roots of the characteristic equation.

On the other hand if $ k = 0 $ this is not true anymore and we get
that $ w = \mathscr{O}(\rho^{-j+1} e^{\tilde{\mu}_{0} \rho})  $, since
$ \tilde{\mu}_{0} $ is the real part of two roots of the
characteristic equation. This prevents us from proving that
$ u_{j} = \mathscr{O}(\rho^{-j} e^{\tilde{\mu}_{0} \rho}) $.

A way around this is to reshuffle the transport equations using the
projection onto the ground state function defined as
\begin{definition}
\label{def:pi}
Denote by $ \pi $ the orthogonal projection onto
$$
L^{2}(\R_{\rho}) \otimes  [\phi_{0}] ,
$$
whose action is described by
$$ 
\pi(f)(x, \rho) = \langle f(\cdot, \rho), \phi_{0}
\rangle_{L^{2}(\R_{x})} \phi_{0}(x) .
$$
\end{definition}
\begin{proposition}
\label{prop:pi-and-P0}
We have that
\begin{equation}
\label{eq:pP}
(1 - \pi) P_{0} = P_{0} (1 - \pi), \qquad  \pi
P_{0} = P_{0} \pi ,
\end{equation}
where $ P_{0} $ is given by \eqref{eq:P0}.
\end{proposition}
\begin{proof}
 Let us prove the first relation. The second has the same proof.
 Compute first the l.h.s.
\begin{align*}
(1 - \pi) P_{0} v &=  (1 - \pi) P_{0} \sum_{k\geq 0} v_{k}
                          \phi_{k} \\
  &= (1 - \pi) \sum_{k \geq 0} \left( v_{k} Q \phi_{k} +
    s_{0}^{-2a} D_{\rho}^{2a} v_{k} \phi_{k}\right) \\
  &= (1 - \pi) \sum_{k \geq 0} \left( \mu_{k} v_{k} +
    s_{0}^{-2a} D_{\rho}^{2a} v_{k} \right)\phi_{k}  \\
  &=  \sum_{k \geq 1} \left( \mu_{k} v_{k} +
    s_{0}^{-2a} D_{\rho}^{2a} v_{k} \right)\phi_{k} .
\end{align*}
On the other hand
\begin{align*}
P_{0} (1 - \pi) v &= P_{0} (1 - \pi) \sum_{k \geq %
                            0} v_{k} \phi_{k} \\
  &= P_{0}  \sum_{k \geq 1} v_{k} \phi_{k}  \\
  &= \sum_{k \geq 1} \left( v_{k} Q \phi_{k} + \phi_{k}  s_{0}^{-2a}
    D_{\rho}^{2a} v_{k} \right) \\
  &= \sum_{k \geq 1} \left( \mu_{k} v_{k} + s_{0}^{-2a}
    D_{\rho}^{2a} v_{k} \right) \phi_{k} .
\end{align*}
This ends the proof.
\end{proof}
Next we are going to choose $ r $ in the definition \eqref{eq:Au} of $
A(u) $. This is necessary in order to bootstrap a formal recursive
calculation of all the $ u_{j} $.
\begin{proposition}
\label{prop:r}
It is possible to choose $ r $ in \eqref{eq:Au} so that
\begin{equation}
\label{eq:P1u0}
\pi P_{1} u_{0} = 0.
\end{equation}
\end{proposition}
\begin{proof}
It is just a computation. By \eqref{eq:P1} the above condition can be
written as
\begin{multline*}
\pi P_{1} u_{0} = \pi \left( \left[ \alpha x
  \partial_{x} + \beta + 2a \gamma_{2a} r\right]
\partial_{\rho}^{2a-1} e^{\mu_{0i^{*}}\rho}  \phi_{0}(x) \right)
\\
= (\partial_{\rho}^{2a-1} e^{\mu_{0i^{*}}\rho} ) \left[ \langle
  \alpha x \partial_{x} \phi_{0},
   \phi_{0} \rangle  + \beta + 2a \gamma_{2a} r \right]\phi_{0}.
\end{multline*}
Here we used the fact that a function of $ \rho $ only commutes with
the projectors and that $ \| \phi_{0} \| = 1 $.

Finally it is clear that the quantity in square brackets can be made
zero by suitably choosing $ r $, since $ \gamma_{2a} \neq 0 $ by
\eqref{eq:P1}, \eqref{eq:gamma}.
\end{proof}
In order to get the optimal decreasing rate $ u_{j} =
\mathscr{O}(\rho^{-j} e^{\tilde{\mu}_{0} \rho}) $, we split the
equation \eqref{eq:transp1} into two sets of equations using the
projection $ \pi $.

First we define $ u_{0} $ by \eqref{eq:u0} and then for $ j
\geq 1$ we solve recursively the equations
\begin{subequations}
\begin{gather}
\label{eq:tr1}    
(1 - \pi) P_{0} u_{j} = P_{0} (1 - \pi) u_{j} =
 -\sum_{k=1}^{\min\{j, 2a\}} \frac{1}{\rho^{k} } (1 - \pi) P_{k}
 u_{j-k} \\
 \label{eq:tr2} 
 \begin{align}  
\pi P_{0} u_{j} &=
P_{0} \pi u_{j} = - \frac{1}{\rho}
\pi P_{1} (1 - \pi) u_{j} - \frac{1}{\rho} \pi
                          P_{1} \pi u_{j-1}
\\
  &\phantom{= =}    - \sum_{k=1}^{\min\{j, 2a-1\}} \frac{1}{\rho^{k+1}}
    \pi P_{k+1} u_{j-k} . \notag
\end{align}
\end{gather}
\end{subequations}
We point out that the idea of the above splitting of the set of the
transport equations is due to M\'etivier, who first used it in
\cite{metivier81}. The motivation for the splitting is that one has to
match the decreasing rate of the $ u_{j} $ for $ \rho $ large.

Moreover the decay rate for large $ \rho $ of the $ (1-\pi)u_{j} $
will turn out to be slightly better than that of the $ \pi u_{j} $,
due to the choice of $ \tilde{\mu}_{0} $ in \eqref{eq:mitlda0}. 
\begin{proposition}
\label{prop:equiv}
The set of equations \eqref{eq:tr1}, \eqref{eq:tr2} for $ j \geq 1 $
as well as the equation $ P_{0} u_{0} = 0 $ 
are formally equivalent to \eqref{eq:transp1}.
\end{proposition}
\begin{proof}
We have, for $ j \geq 1 $,
\begin{align}
  \label{eq:P0uj}
P_{0} u_{j} &= P_{0}(1 - \pi) u_{j} + P_{0} \pi u_{j}
  \\
  &= - \sum_{\ell=1}^{\min\{j, 2a\}} \frac{1}{\rho^{\ell} } (1 - \pi) P_{\ell}
    u_{j-\ell} - \frac{1}{\rho} \pi P_{1} (1 - \pi) u_{j}
    \notag \\
  & \phantom{=} \ - \frac{1}{\rho} \pi P_{1} \pi u_{j-1}
  -  \sum_{\ell=1}^{\min\{j, 2a -1\}} \frac{1}{\rho^{\ell+1}} \pi
    P_{\ell+1} u_{j-\ell}
    \notag
\end{align}
Since
$$ 
\sum_{k=0}^{2a} \frac{1}{\rho^{k}} P_{k} \sum_{j \geq 0} u_{j} = 0 
$$
iff
$$ 
\sum_{j \geq 1} P_{0} u_{j} + \sum_{\substack{j \geq 0\\ k \geq 1}}
\frac{1}{\rho^{k}} P_{k} u_{j} = 0 ,
$$
from \eqref{eq:P0uj} we get
\begin{align*}
\sum_{j \geq 1} P_{0} u_{j} &= - \sum_{j \geq 1} \sum_{\ell=1}^{\min\{j, 2a\}} \frac{1}{\rho^{\ell} } (1 - \pi) P_{\ell}
    u_{j-\ell} - \frac{1}{\rho} \sum_{j \geq 1} \pi P_{1} (1 -
                              \pi) u_{j}
\\
&
\phantom{=} \ - \frac{1}{\rho} \sum_{j \geq 1} \pi P_{1} \pi u_{j-1}
  - \sum_{j \geq 1} \sum_{\ell=1}^{\min\{j, 2a-1\}} \frac{1}{\rho^{\ell+1}} \pi
    P_{\ell+1} u_{j-\ell} .
\end{align*}
The third term on the r.h.s. above is written as
\begin{align*}
\frac{1}{\rho} \sum_{j \geq 1} \pi P_{1} \pi u_{j-1} &=
        \frac{1}{\rho} \pi P_{1} \pi u_{0} + \frac{1}{\rho}
        \sum_{j \geq 2} \pi P_{1} \pi u_{j-1} \\
&=
\frac{1}{\rho} \pi P_{1} u_{0} + \frac{1}{\rho}
        \sum_{j \geq 1} \pi P_{1} \pi u_{j} .
\end{align*}
Hence
\begin{align*}
\sum_{j \geq 1} P_{0} u_{j} &= - \frac{1}{\rho} (1 - \pi) P_{1} u_{0}  
- \sum_{j \geq 2} \sum_{\ell=1}^{\min\{j, 2a\}} \frac{1}{\rho^{\ell} } (1 - \pi) P_{\ell}
    u_{j-\ell} \\
&\phantom{ = }\
- \frac{1}{\rho} \sum_{j \geq 1} \pi P_{1} (1 - \pi) u_{j} 
- \frac{1}{\rho} \pi P_{1} u_{0}
\\
&\phantom{ = }\
- \frac{1}{\rho} \sum_{j \geq 1} \pi P_{1} \pi u_{j}
- \sum_{j \geq 1} \sum_{\ell=1}^{\min\{j, 2a-1\}} \frac{1}{\rho^{\ell+1}} \pi
    P_{\ell+1} u_{j-\ell}
\\
&=
- \frac{1}{\rho} P_{1} u_{0} - \frac{1}{\rho} \sum_{j \geq 1}
     \pi P_{1} u_{j} - 
\sum_{j \geq 2} \sum_{\ell=1}^{\min\{j, 2a\}} \frac{1}{\rho^{\ell} } (1 - \pi) P_{\ell}
    u_{j-\ell}
  \\
  &\phantom{ = }\
- \sum_{j \geq 1} \sum_{\ell=1}^{\min\{j, 2a-1\}} \frac{1}{\rho^{\ell+1}} \pi
    P_{\ell+1} u_{j-\ell} .
\end{align*}
The third term in that last expression above when $ \ell = 1 $ is
$$ 
\sum_{j \geq 2} \frac{1}{\rho} (1 - \pi) P_{1} u_{j-1} =
\frac{1}{\rho} \sum_{j \geq 1} (1 - \pi) P_{1} u_{j} ,
$$
so that
\begin{align*}
\sum_{j \geq 1} P_{0} u_{j} &= - \frac{1}{\rho} P_{1} u_{0} - \frac{1}{\rho} \sum_{j \geq 1}
     \pi P_{1} u_{j} - \frac{1}{\rho} \sum_{j \geq 1} (1 -
                              \pi) P_{1} u_{j}
\\
&\phantom{=}\
- \sum_{j \geq 2} \sum_{\ell=2}^{\min\{j, 2a\}} \frac{1}{\rho^{\ell} } (1 - \pi) P_{\ell}
    u_{j-\ell}
\\
&\phantom{=}\
- \sum_{j \geq 1} \sum_{\ell=1}^{\min\{j, 2a-1\}} \frac{1}{\rho^{\ell+1}} \pi
    P_{\ell+1} u_{j-\ell} .
\end{align*}
Then
\begin{align*}
\sum_{j \geq 1} P_{0} u_{j} &= -  \frac{1}{\rho} P_{1} u_{0} - \frac{1}{\rho} \sum_{j \geq 1}
     P_{1} u_{j}  
- \sum_{j \geq 1} \sum_{\ell=1}^{\min\{j, 2a-1\}} \frac{1}{\rho^{\ell+1}} \pi
     P_{\ell+1} u_{j-\ell}
  \\
&\phantom{=}\ 
- \sum_{j \geq 1} \sum_{\ell=1}^{\min\{j+1, 2a\}-1}
     \frac{1}{\rho^{\ell+1} } (1 - \pi) P_{\ell+1} u_{j-\ell}.
\end{align*}
And finally, since $ \min\{j+1, 2a\}-1 = \min\{j, 2a-1\} $,
\begin{align*}
\sum_{j \geq 1} P_{0} u_{j} &= - \frac{1}{\rho} \sum_{j \geq 0} P_{1}
                              u_{j}
- \sum_{j \geq 1} \sum_{\ell=1}^{\min\{j, 2a-1\}}
     \frac{1}{\rho^{\ell+1} } P_{\ell+1} u_{j-\ell}
\\
&=
 - \frac{1}{\rho} P_{1} u_{0}  - \sum_{j \geq 1} \sum_{\ell=0}^{\min\{j, 2a-1\}}
     \frac{1}{\rho^{\ell+1} } P_{\ell+1} u_{j-\ell}
\\
&=
- \sum_{j \geq 0} \sum_{\ell=0}^{\min\{j, 2a-1\}}
     \frac{1}{\rho^{\ell+1} } P_{\ell+1} u_{j-\ell}                      
\end{align*}
Since $ P_{0} u_{0} = 0 $ we have
$$ 
\sum_{j \geq 0} P_{0} u_{j} + \sum_{j \geq 0} \sum_{\ell=0}^{\min\{j+1, 2a\}-1}
     \frac{1}{\rho^{\ell+1} } P_{\ell+1} u_{j-\ell} = 0.
$$
Changing the indices from $ (j, \ell) $ to $ (t, \ell) $, $ \ell \geq
0 $, $ j-\ell = t $, we conclude
$$ 
\sum_{t \geq 0} \sum_{\ell=0}^{2a} \frac{1}{\rho^{\ell}} P_{\ell}
u_{t} = 0.
$$
\end{proof}
The equations in \eqref{eq:tr1} and \eqref{eq:tr2} can be solved, in
principle. We must however make sure that the so obtained solutions
have estimates allowing us to turn our formal solution into a real
one.

Since the degeneration of the operator coefficients is higher than
quadratic, there is no possibility of using the three terms recurrence
relations valid for the quadratic case or in some special case, like
in \cite{metivier81} and \cite{chinni}. This leads us to getting
estimates directly by inspecting the form of the $ u_{j} $.

This technique, although more involved than the classical one, seems
however promising for more general (and generic) cases in two
variables. This will be the object of the next section.

\section{Sobolev Type Estimates of the $ u_{j} $}
\setcounter{equation}{0}
\setcounter{theorem}{0}
\setcounter{proposition}{0}
\setcounter{lemma}{0}
\setcounter{corollary}{0}
\setcounter{definition}{0}

In this section we deduce estimates of the functions $ u_{j} $ solving
\eqref{eq:tr1}, \eqref{eq:tr2} in the region $ \rho > \textrm{Const}\ j
$.

As it can be seen from \eqref{eq:P0inv} the functions $ u_{j} $, or
rather their projections, are computed as infinite expansions in the
eigenfunctions of a anharmonic oscillator. It is then natural to use
the $ L^{2} $ norms of the $ u_{j} $ and of their derivatives, to take
advantage of the orthonormal system of the $ \phi_{k} $.

We point out that in M\'etivier's case each $ u_{j} $ is given by a finite
expansion in the eigenfunctions and furthermore their derivatives are
a finite linear combination of the same eigenfunctions and allow to
get the $ L^{\infty} $ estimates of the $ u_{j} $ in a direct way.

In the present setting the pointwise estimates will be deduced
using the Sobolev embedding theorems in Section \ref{sec:pt}.

First we define the weight function (see \eqref{eq:mitlda0})
\begin{equation}
\label{eq:wj}
w_{j}(\rho) = e^{| \tilde{\mu}_{0} | \rho} \rho^{(j - \delta)\kappa},
\end{equation}
where
\begin{equation}
\label{eq:kappadelta}
0 < \delta < 1, \qquad  \frac{1}{s_{0}} < \kappa < 1, \qquad
\kappa \delta > \frac{1}{2}.
\end{equation}
We point out that the role of $ \delta $, $ \kappa $ is exquisitely
technical and due to the use of $ L^{2} $ norms, i.e.  making certain
integrals in Lemmas \ref{lemma:pi-mu0}, \ref{lemma:piP-1-pi}
absolutely convergent. 

We define the Sobolev spaces $ B^{k}(\R_{x}) $ as the space of
all $ L^{2} $ functions such that
\begin{equation}
\label{eq:spBk}
\| u\|_{k} = \max_{\frac{\beta}{q-1} + \alpha \leq k}  \| x^{\beta}
D_{x}^{\alpha} u \|_{0} .
\end{equation}
For the estimates we shall use the norms
\begin{equation}
\label{eq:norms}
\| w_{j}(\rho) f(x, \rho) \|_{k, A}^{2} =
\int_{A}^{+\infty} w_{j}^{2}(\rho) \| f( \cdot, \rho) \|_{k}^{2} d\rho ,
\end{equation}
where
$ A $ is a suitable positive constant to be chosen later.

We want to prove the
\begin{theorem}
\label{th:est}
There exist positive constants $ C_{u} $, $ R_{0} $ such that
\begin{multline}
\label{eq:est}
\| w_{j}(\rho) \partial_{\rho}^{\gamma} x^{\beta}
\partial_{x}^{\alpha} u_{j}(x, \rho)  \|_{0, A}
\\
\leq
C_{u}^{1 + j + \alpha + \beta + \gamma} 
\cdot
\left( \frac{j}{s_{0}} + \alpha \frac{q-1}{q} +
  \frac{\beta}{q} + \frac{\gamma}{a} \frac{q-1}{q} \right)! ,
\end{multline}
where
\begin{equation}
  \label{eq:AR}
\gamma \leq j + \gamma^{\#}, \qquad  A = A(j) = R_{0}(j + 1) ,
\end{equation}
and $ R_{0} > 0 $, $ \gamma^{\#} $ a positive constant.

Here we understand that for a positive number $ x $, $ x! = \Gamma(x+1) $.
\end{theorem}
Using the Sobolev immersion theorem, from Theorem \ref{th:est} we
deduce the desired pointwise estimates \eqref{eq:pest}. 

\vskip 2mm

We are going to prove a slightly different statement. Actually in the
theorem below both space derivatives and multiplication by $ x $ have
been replaced by the corresponding powers of the operator $ Q $. This
has the advantage that the action of $ Q $ onto an eigenfunction
expansion keeps the orthonormality of the basis, while this is not
true for both derivatives and multiplication by $ x $.

The transition from $ Q $ to the derivatives and multiplication by $ x
$ is done in Proposition \ref{prop:ineqQ}. 
\begin{theorem}
  \label{th:est'}
Let $ \nu $ denote a rational number $ \geq -1 $ and $ \gamma \in \N $
such that $ \nu + \gamma (2a)^{-1} \geq 0 $.
There exist positive constants $ C_{u} $, $ R_{0} $, $ \sigma $, $
\sigma' $ such that 
\begin{equation}
\label{eq:est'}
\| w_{j}(\rho) \partial_{\rho}^{\gamma} Q^{\nu}
 u_{j}(x, \rho)  \|_{0, A}
\leq
C_{0}^{1 + \sigma j + \sigma'( \nu + \gamma)} 
(\lambda(j, \nu, \gamma)+1)^{\lambda(j, \nu, \gamma)} ,
\end{equation}
where
$$ 
\lambda(j, \nu, \gamma) = \frac{j}{s_{0}} + \left(2 \nu +
  \frac{\gamma}{a}\right) \frac{q-1}{q},
$$
and
$ A = A(j) = R_{0}(j + 1)  $, $ R_{0}> 0 $, $ \gamma \leq j + \gamma^{\#} $.
\end{theorem}
\begin{remark}
\label{rem:c0}
We point out that the constants $ R_{0} $, $ C_{0} $ may be enlarged
independently from each other.
\end{remark}

We shall proof in the following that Theorem \ref{th:est'} actually
implies Theorem \ref{th:est}. First we prove Theorem \ref{th:est'}.


\bigskip
The proof proceeds by induction on $ j $. Consider first $ u_{0} $
defined in \eqref{eq:u0}. We have for $ A = A(0, \gamma) $, 
\begin{multline*}
\| w_{0}(\rho) \partial_{\rho}^{\gamma} Q^{\nu}
u_{0}(x, \rho)  \|_{0, A}^{2}
= \int_{A}^{+\infty} \rho^{-2\delta \kappa} e^{2 | \tilde{\mu}_{0} |
  \rho} \left|\partial_{\rho}^{\gamma} e^{\mu_{0i^{*}} \rho} \right|^{2}
\| Q^{\nu} \phi_{0}(x) \|_{0}^{2} d\rho
\\
\leq C^{\gamma} \| Q^{\nu}  \phi_{0}(x) \|_{0}^{2} 
\leq C^{\gamma} \mu_{0}^{\nu}
\leq
 C_{u}^{\gamma} \lambda(0, \nu, \gamma)^{\lambda(0, \nu, \gamma)} .
\end{multline*}
Assume now that \eqref{eq:est'} holds for any $ k $, $ 0 \leq k < j
$. We recall that
$$ 
u_{j} = (1 - \pi) u_{j} + \pi u_{j} ,
$$
where the summands in the r.h.s. above are defined by \eqref{eq:tr1}
and \eqref{eq:tr2} respectively.

We have
\begin{multline}
\label{eq:est''}
\| w_{j}(\rho) \partial_{\rho}^{\gamma} Q^{\nu} u_{j}(x, \rho)  \|_{0, A}
\\
\leq
\| w_{j}(\rho) \partial_{\rho}^{\gamma} Q^{\nu}
 (1 - \pi) u_{j}(x, \rho)  \|_{0, A} + \| w_{j}(\rho) \partial_{\rho}^{\gamma} Q^{\nu}
 \pi u_{j}(x, \rho)  \|_{0, A} .
\end{multline}
We are going to prove estimates for each of the two summands on the
right hand side of \eqref{eq:est''}.

\section{Estimate of  $(1 - \pi) u_{j} $}
\setcounter{equation}{0}
\setcounter{theorem}{0}
\setcounter{proposition}{0}
\setcounter{lemma}{0}
\setcounter{corollary}{0}
\setcounter{definition}{0}

This section is devoted to proving the estimates \eqref{eq:est'} for $
(1 - \pi) u_{j} $. Actually we shall prove a slightly better kind of
inequality only for $ (1 - \pi)u_{j} $. This improvement will be
crucial in inductively proving the estimates \eqref{eq:est'} for $ \pi
u_{j}$.  

\begin{theorem}
\label{th:est-1-p}
Let $ \nu $ denote a rational number $ \geq -1 $ and $ \gamma \in \N $
such that $ \nu + \gamma (2a)^{-1} \geq 0 $.
There exist positive constants $ C_{u} $, $ R_{0} $, $ \sigma $, $
\sigma' $ such that 
\begin{multline}
\label{eq:est-1-p}
\| w_{j}(\rho) \rho^{1-\kappa} \partial_{\rho}^{\gamma} Q^{\nu}
(1 - \pi) u_{j}(x, \rho)  \|_{0, A}
\\
\leq
C_{0}^{1 + \sigma j + \sigma'( \nu + \gamma)} 
(\lambda(j, \nu, \gamma)+1)^{\lambda(j, \nu, \gamma)} ,
\end{multline}
where
$$ 
\lambda(j, \nu, \gamma) = \frac{j}{s_{0}} + \left(2 \nu +
  \frac{\gamma}{a}\right) \frac{q-1}{q},
$$
and
$ A = A(j) = R_{0}(j + 1) $, $ R_{0} > 0 $, $ \gamma \leq j + \gamma^{\#} $.
\end{theorem}
\begin{proof}
We have, by \eqref{eq:tr1}, 
$$ 
(1 - \pi) u_{j} =
 -\sum_{k=1}^{\min\{j, 2a\}} P_{0}^{-1} \left( \frac{1}{\rho^{k} } (1 - \pi) P_{k}
 u_{j-k} \right) ,
$$
where $ P_{0}^{-1} $ has been defined in \eqref{eq:P0inv}. Hence
\begin{multline*}
\| w_{j}(\rho) \rho^{1-\kappa} \partial_{\rho}^{\gamma} Q^{\nu}
(1 - \pi) u_{j}(x, \rho)  \|_{0, A}
\\
\leq
\sum_{k=1}^{\min\{j, 2a\}} \| w_{j}(\rho) \rho^{1-\kappa}  \partial_{\rho}^{\gamma}
Q^{\nu} P_{0}^{-1} \left( \frac{1}{\rho^{k} } (1 - \pi) P_{k}
 u_{j-k} \right ) \|_{0, A} .
\end{multline*}
\begin{lemma}
\label{lemma:QnuP0}
We have
$$ 
Q^{\nu} P_{0}^{-1}(f) = P_{0}^{-1}( Q^{\nu} f), \qquad Q^{\nu} (1 -
\pi) = (1 - \pi) Q^{\nu}.
$$
\end{lemma}
\begin{proof}[Proof of Lemma \ref{lemma:QnuP0}]
From \eqref{eq:P0inv} we get (here $ f_{k} = \langle f, \phi_{k}
\rangle $) 
\begin{multline*}
Q^{\nu} P_{0}^{-1}(f) = Q^{\nu} \sum_{k\geq 0} \left(
  \sum_{j=1}^{2a} A_{kj} I_{kj}(f_{k})(\rho) \right) \phi_{k}(x)
\\
=
\sum_{k\geq 0} \left(
  \sum_{j=1}^{2a} A_{kj} I_{kj}(f_{k})(\rho) \right) \mu_{k}^{\nu}
\phi_{k}(x)
\\
=
\sum_{k\geq 0} \left(
  \sum_{j=1}^{2a} A_{kj} I_{kj}(\mu_{k}^{\nu} f_{k})(\rho) \right) 
\phi_{k}(x)
\\
=
\sum_{k\geq 0} \left(
  \sum_{j=1}^{2a} A_{kj} I_{kj}( ( Q^{\nu} f)_{k})(\rho) \right) 
\phi_{k}(x) .
\end{multline*}
The proof of the second relation is straightforward.
\end{proof}
Using the above lemma we get
\begin{multline}
  \label{eq:uj-2}
\| w_{j}(\rho) \rho^{1-\kappa} \partial_{\rho}^{\gamma} Q^{\nu}
(1 - \pi) u_{j}(x, \rho)  \|_{0, A}
\\
\leq
\sum_{k=1}^{\min\{j, 2a\}} \| w_{j}(\rho) \rho^{1-\kappa} \partial_{\rho}^{\gamma}
 P_{0}^{-1} \left( \frac{1}{\rho^{k} } (1 - \pi) Q^{\nu} P_{k}
 u_{j-k} \right ) \|_{0, A} .
\end{multline}
\begin{lemma}
\label{lemma:drhogamma}
Let $ f \in \mathscr{S}(\R) $ and $ \gamma \in \N $. Write $ \gamma =
2a r + \gamma_{1} $, with $ 0 \leq \gamma_{1} < 2a $. Then for $ k
\geq 0 $ we have
\begin{multline}
\label{eq:drhogamma}
\left(\frac{1}{s_{0}} D_{\rho}\right)^{\gamma} E_{k}(f_{k})(\rho)
=
\sum_{s=0}^{r-1} (- \mu_{k})^{s} \left(\frac{1}{s_{0}}
  D_{\rho}\right)^{\gamma - 2a (1+s)} f_{k}
\\
+ (- \mu_{k})^{r}
\left(\frac{1}{s_{0}} D_{\rho}\right)^{\gamma_{1}} E_{k}(f_{k})(\rho), 
\end{multline}
where $ E_{k}(f_{k}) $ is defined in \eqref{eq:Ek} and
$ f_{k} = \langle f, \phi_{k} \rangle $, $ \phi_{k} $ denoting the $ k
$-th eigenfunction of $ Q $.
\end{lemma}
The proof is just a computation using the fact that $ E_{k}(f_{k}) $
is a solution of \eqref{eq:ode-k} rapidly decreasing at infinity.

\begin{corollary}
\label{cor:drhogammaP0}
Let $ f \in \mathscr{S}(\R) $ and $ \gamma \in \N $. Write $ \gamma =
2a r + \gamma_{1} $, with $ 0 \leq \gamma_{1} < 2a $. Then
\begin{multline}
\label{eq:drhogammaP0}
\left(\frac{1}{s_{0}} D_{\rho}\right)^{\gamma} P_{0}^{-1}(f)
=
\sum_{s=0}^{r-1} (- 1)^{s} \left(\frac{1}{s_{0}}
  D_{\rho}\right)^{\gamma - 2a (1+s)} Q^{s} f
\\
+ (- 1)^{r}
\left(\frac{1}{s_{0}} D_{\rho}\right)^{\gamma_{1}} P_{0}^{-1}( Q^{r} f) . 
\end{multline}
\end{corollary}
Using Corollary \ref{cor:drhogammaP0} inequality \eqref{eq:uj-2}
becomes
\begin{multline}
\label{eq:uj-3}
\| w_{j}(\rho) \rho^{1-\kappa}  \Dslash_{\rho}^{\gamma} Q^{\nu}
(1 - \pi) u_{j}(x, \rho)  \|_{0, A}
\\
\leq
\sum_{k=1}^{\min\{j, 2a\}} \sum_{s=0}^{r-1} \| w_{j}(\rho) \rho^{1-\kappa}
\Dslash^{\gamma - 2a(1+s)}_{\rho} \left( \frac{1}{\rho^{k} } (1 - \pi)
  Q^{\nu +s} P_{k} u_{j-k} \right ) \|_{0, A}
\\
+
\sum_{k=1}^{\min\{j, 2a\}} \| w_{j}(\rho) \rho^{1-\kappa}  \Dslash^{\gamma_{1}}_{\rho}
P_{0}^{-1} \left( \frac{1}{\rho^{k} } (1 -  \pi)
  Q^{\nu +r} P_{k} u_{j-k} \right ) \|_{0, A} ,
\end{multline}
where we used the notation
\begin{equation}
  \label{eq:dslash}
\Dslash_{\rho} = \frac{1}{s_{0}} D_{\rho} = \frac{1}{i s_{0}}
\partial_{\rho} . 
\end{equation}
First consider the norm in the third line of \eqref{eq:uj-3}. We
denote by $ g_{jk} = Q^{\nu +r} P_{k} u_{j-k} $. Then for $ k \in \{1,
\ldots, \min\{j, 2a\}\}$,
\begin{multline*}
\| w_{j}(\rho)\rho^{1-\kappa}  \Dslash^{\gamma_{1}}_{\rho}
P_{0}^{-1} \left( \frac{1}{\rho^{k} } (1 -  \pi) g_{jk}  \right )
\|_{0, A}
\\
=
\| w_{j}(\rho) \rho^{1-\kappa}  \Dslash^{\gamma_{1}}_{\rho} \sum_{\ell \geq 1}
\sum_{i=1}^{2a} A_{\ell i} I_{\ell i}(\rho^{-k} \langle g_{j k} ,
\phi_{\ell} \rangle ) \phi_{\ell} \|_{0, A}
\\
=
\left( \sum_{\ell \geq 1} \| w_{j}(\rho) \rho^{1-\kappa} \Dslash^{\gamma_{1}}_{\rho}
  \sum_{i=1}^{2a} A_{\ell i} I_{\ell i}(\rho^{-k} \langle g_{j k} , 
  \phi_{\ell} \rangle ) \|_{A}^{2} \right)^{\frac{1}{2}} ,
\end{multline*}
where
$$ 
\| f(\rho) \|_{A}^{2} = \int_{A}^{+\infty} | f(\rho) |^{2} d \rho .
$$
Remark that, by \eqref{eq:Ikj}, \eqref{eq:Amuk}, since $ 0 \leq
\gamma_{1} < 2a $, 
\begin{equation}
\label{eq:Drhogamma1}
\Dslash^{\gamma_{1}}_{\rho} \sum_{i=1}^{2a} A_{\ell i} I_{\ell
  i}(f_{k}) = \frac{1}{(i s_{0})^{\gamma_{1}}} \sum_{i=1}^{2a} A_{\ell
  i} \mu_{\ell i}^{\gamma_{1}} I_{\ell i}(f_{k}) . 
\end{equation}
Hence
\begin{multline*}
\| w_{j}(\rho) \rho^{1-\kappa} \Dslash^{\gamma_{1}}_{\rho}
P_{0}^{-1} \left( \frac{1}{\rho^{k} } (1 -  \pi) g_{jk}  \right )
\|_{0, A}
\\
= s_{0}^{-\gamma_{1}} 
\left( \sum_{\ell \geq 1} \| w_{j}(\rho) \rho^{1-\kappa}
  \sum_{i=1}^{2a} A_{\ell i} \mu_{\ell i}^{\gamma_{1}}  I_{\ell i}(\rho^{-k} \langle g_{j k} , 
  \phi_{\ell} \rangle ) \|_{A}^{2} \right)^{\frac{1}{2}} 
\end{multline*}
Now
$$ 
A_{\ell i} = \prod_{\substack{r=1\\r\neq i}}^{2a} \frac{1}{\mu_{\ell
    i} - \mu_{\ell r}} .
$$
Since the roots $ \mu_{\ell r} $ are the vertices of a $ 2a $-regular
polygon whose circumscribed circle has radius $
\mu_{\ell}^{\frac{1}{2a}} $ (arithmetic root), we have that
$$ 
\frac{|A_{\ell i} \mu_{\ell i}^{\gamma_{1}
    } |}{\mu_{\ell}^{\frac{\gamma_{1}+ 1 - 2a}{2a}}} = c ,
$$
where $ c $ denotes a positive constant independent of $ \ell $, $ i
$, $ \gamma_{1} $. Thus
\begin{multline}
  \label{eq:uj-4}
\| w_{j}(\rho) \rho^{1-\kappa}  \Dslash^{\gamma_{1}}_{\rho}
P_{0}^{-1} \left( \frac{1}{\rho^{k} } (1 -  \pi) g_{jk}  \right )
\|_{0, A}
\\
\leq  C_{1} \sum_{i=1}^{2a} 
\left( \sum_{\ell \geq 1} \| w_{j}(\rho) \rho^{1-\kappa} 
 I_{\ell i}(\rho^{-k} \langle g_{j k} , 
 \mu_{\ell}^{\frac{\gamma_{1}+ 1 - 2a}{2a}}  \phi_{\ell} \rangle )
 \|_{A}^{2} \right)^{\frac{1}{2}}  
\end{multline}
\begin{lemma}
\label{lemma:muli>0}
Assume that $ \ell \geq 1 $, $ \re \mu_{\ell i} > 0 $. Then
\begin{equation}
\label{eq:muli>0}
\| w_{j}(\rho) \rho^{1-\kappa}  I_{\ell i}( \rho^{-k} f(\rho)) \|_{A}^{2} \leq C_{0} \|
w_{j-k}(\rho) \mu_{\ell}^{-\frac{1}{2a}} f(\rho) \|_{A}^{2} ,
\end{equation}
for a suitable constant $ C_{0} > 0 $ independent of $ \ell $, $ j $,
$ k $, $ i $.
\end{lemma}
\begin{remark}
In \eqref{eq:muli>0} the factor $ \rho^{1-\kappa} $ appears only in
the norms on the left hand side of the inequality because the right
hand side may contain both $ (1-\pi) u_{j-k} $ and $ \pi u_{j-k} $. 
\end{remark}
\begin{proof}
By \eqref{eq:Ikj}, \eqref{eq:mitlda0}, \eqref{eq:wj} we have
\begin{multline*}
\| w_{j}(\rho) \rho^{1-\kappa}  I_{\ell i}( \rho^{-k} f(\rho)) \|_{A}^{2}
=
\int_{A}^{+\infty} \rho^{2(j-\delta) \kappa + 2(1-\kappa)} e^{2|\tilde{\mu}_{0}|\rho} \Big|
\int_{\R} e^{\mu_{\ell i} (\rho - \sigma)}
\\
\cdot
H\left(- \sgn(\re \mu_{\ell
      i}) (\rho - \sigma)\right)
 \frac{1}{\sigma^{k}} \  H(\sigma - R) f(\sigma) d\sigma \Big|^{2}
 d\rho
 \\
 \leq
 \int_{A}^{+\infty} \Big(
  \int_{\R} e^{(\re \mu_{\ell i} + |\tilde{\mu}_{0}|) (\rho - \sigma)} H\left(- \sgn(\re \mu_{\ell
      i}) (\rho - \sigma)\right)
  \\
  \cdot
  \sigma^{ - (k - 1) (1 - \kappa) }
  e^{|\tilde{\mu}_{0}| \sigma}
 \sigma^{(j-k-\delta)\kappa} \  H(\sigma - R) | f(\sigma)| d\sigma \Big)^{2}
 d\rho
 \\
 \leq
 \int_{A}^{+\infty} \Big(
  \int_{\R} e^{(\re \mu_{\ell i} + |\tilde{\mu}_{0}|) (\rho - \sigma)} H\left(- \sgn(\re \mu_{\ell
      i}) (\rho - \sigma)\right)
  \\
  \cdot
  \sigma^{ - (k - 1) (1 - \kappa) }
  e^{|\tilde{\mu}_{0}| \sigma}
 \sigma^{(j-k-\delta)\kappa} \  H(\sigma - A) | f(\sigma)| d\sigma \Big)^{2}
 d\rho
\end{multline*}
where, since $ \sigma \geq \rho $, we bounded $ \rho $ by $ \sigma
$. Moreover for the same reason $ \sigma \geq A $ so that we may replace
$ H(\sigma-R)$ by $ H(\sigma-A) $.

We estimate the integral on the right hand side of the above
expression with an integral over the whole real line and apply Young
inequality to obtain
\begin{multline*}
\| w_{j}(\rho) \rho^{1-\kappa} I_{\ell i}( \rho^{-k} f(\rho)) \|_{A}^{2}
\\
\leq
 \int_{-\infty}^{+\infty} \Big(
  \int_{\R} e^{(\re \mu_{\ell i} + |\tilde{\mu}_{0}|) (\rho - \sigma)} H\left(- \sgn(\re \mu_{\ell
      i}) (\rho - \sigma)\right)
  \\
  \cdot
 w_{j-k}(\sigma)  H(\sigma - A) \sigma^{ - (k - 1) (1 - \kappa) } | f(\sigma)| d\sigma \Big)^{2}
 d\rho
 \\
 \leq
 C_{1}
 \| e^{(\re \mu_{\ell i} + |\tilde{\mu}_{0}|) \rho} H\left(- \sgn(\re \mu_{\ell
     i}) \rho \right) \|_{L^{1}(\R)}^{2}
\
 \|  w_{j-k}(\sigma) \sigma^{ - (k - 1) (1 - \kappa) } f(\sigma) \|_{A}^{2}
 \\
 \leq
 C_{2} \mu_{\ell}^{-\frac{1}{2a}}  \|  w_{j-k}(\sigma) f(\sigma)
 \|_{A}^{2} .
\end{multline*}
This concludes the proof of the lemma.
\end{proof}
\begin{lemma}
\label{lemma:muli<0}
Assume that $ \ell \geq 1 $, $ \re \mu_{\ell i} < 0 $. Then
\begin{equation}
\label{eq:muli<0}
\| w_{j}(\rho) \rho^{1-\kappa}  I_{\ell i}( \rho^{-k} f(\rho)) \|_{A}^{2} \leq C_{0} \|
w_{j-k}(\rho) \mu_{\ell}^{-\frac{1}{2a}} f(\rho) \|_{A}^{2} ,
\end{equation}
for a suitable constant $ C_{0} > 0 $ independent of $ \ell $, $ j $,
$ k $, $ i $, and $ R = A(j)= R_{0}(j+1)$.
\end{lemma}
\begin{proof}
As before by \eqref{eq:Ikj}, \eqref{eq:mitlda0}, \eqref{eq:wj} we have
\begin{multline*}
\| w_{j}(\rho) \rho^{1-\kappa}  I_{\ell i}( \rho^{-k} f(\rho)) \|_{A}^{2}
=
\int_{A}^{+\infty} \rho^{2(j-\delta)\kappa + 2(1-\kappa)} e^{2|\tilde{\mu}_{0}|\rho} \Big|
\int_{\R} e^{\mu_{\ell i} (\rho - \sigma)}
\\
\cdot
H\left(- \sgn(\re \mu_{\ell
      i}) (\rho - \sigma)\right)
 \frac{1}{\sigma^{k}} \  H(\sigma - A) f(\sigma) d\sigma \Big|^{2}
 d\rho
 \\
 \leq
 \int_{A}^{+\infty} \Big( \rho^{(j-\delta)\kappa + 1 -\kappa}
 \int_{\R} e^{(\re \mu_{\ell i} + |\tilde{\mu}_{0}|) (\rho - \sigma)}
 \\
 \cdot
 H\left(- \sgn(\re \mu_{\ell
      i}) (\rho - \sigma)\right)
  e^{|\tilde{\mu}_{0}| \sigma}
 \sigma^{-k} \  H(\sigma - A) | f(\sigma)| d\sigma \Big)^{2}
 d\rho
 \\
 =
 \int_{A}^{+\infty} \Big(
  \int_{\R} e^{(\re \mu_{\ell i} + |\tilde{\mu}_{0}|) (\rho - \sigma)} H\left(- \sgn(\re \mu_{\ell
      i}) (\rho - \sigma)\right)
  \\
  \cdot \left( 1 + \frac{\rho -
      \sigma}{\sigma}\right)^{(j-\delta)\kappa + 1 - \kappa}
  e^{|\tilde{\mu}_{0}| \sigma}
 \sigma^{(j - k - \delta) \kappa} \  H(\sigma - A)
 \sigma^{-(k-1)(1-\kappa)}| f(\sigma)| d\sigma \Big)^{2}  d\rho
\end{multline*}
Now, since $ (j - \delta)\kappa + 1 -\kappa \leq j $,
$$ 
\left( 1 + \frac{\rho - \sigma}{\sigma}\right)^{(j-\delta)\kappa + 1 -
\kappa} \leq \left(
  1 + \frac{\rho - \sigma}{\sigma}\right)^{j} ,
$$
keeping in mind that $ A \leq \sigma \leq \rho $, we have
\begin{multline*}
\| w_{j}(\rho) \rho^{1-\kappa} I_{\ell i}( \rho^{-k} f(\rho))
\|_{A}^{2}
\\
\leq
\int_{A}^{+\infty} \Big(
  \int_{\R} e^{(\re \mu_{\ell i} + |\tilde{\mu}_{0}|) (\rho - \sigma)} H\left(- \sgn(\re \mu_{\ell
      i}) (\rho - \sigma)\right)
  \\
  \cdot
  \left( \sum_{r=0}^{j} \binom{j}{r} \frac{(\rho - \sigma)^{r}}{\sigma^{r}}  \right)
  e^{|\tilde{\mu}_{0}| \sigma}
  \sigma^{(j - k - \delta)\kappa} \  H(\sigma - A)
  \sigma^{-(k-1)(1-\kappa)} | f(\sigma)| d\sigma \Big)^{2}
 d\rho
 \\
 =
\int_{A}^{+\infty} \Big( \sum_{r=0}^{j} \binom{j}{r} 
  \int_{\R} e^{(\re \mu_{\ell i} + |\tilde{\mu}_{0}|) (\rho - \sigma)} H\left(- \sgn(\re \mu_{\ell
      i}) (\rho - \sigma)\right)
  \\
  \cdot   \frac{(\rho - \sigma)^{r}}{\sigma^{r}} \ 
  e^{|\tilde{\mu}_{0}| \sigma}
 \sigma^{(j - k - \delta)\kappa} \  H(\sigma - A)  \sigma^{-(k-1)(1-\kappa)} | f(\sigma)| d\sigma \Big)^{2}
 d\rho .
\end{multline*}
Using the inequality
\begin{equation}
\label{eq:diseg}
\left( \sum_{i=1}^{n} a_{i} \right)^{2} \leq \sum_{i=1}^{n-1} 2^{i}
a_{i}^{2} + 2^{n-1} a_{n}^{2} \leq \sum_{i=1}^{n} 2^{i} a_{i}^{2} ,
\end{equation}
we get
\begin{multline*}
\| w_{j}(\rho) \rho^{1-\kappa} I_{\ell i}( \rho^{-k} f(\rho))
\|_{A}^{2}
\\
\leq
\sum_{r=0}^{j} 2^{r+1} \binom{j}{r}^{2} \int_{A}^{+\infty} \Big(
\int_{\R} e^{(\re \mu_{\ell i} + |\tilde{\mu}_{0}|) (\rho - \sigma)}
(\rho - \sigma)^{r}
 H\left(- \sgn(\re \mu_{\ell i}) (\rho - \sigma)\right)
\\
\cdot
\sigma^{-(k-1)(1-\kappa) - r}
  e^{|\tilde{\mu}_{0}| \sigma}
 \sigma^{(j - k - \delta)\kappa} \  H(\sigma - A) | f(\sigma)| d\sigma \Big)^{2}
 d\rho .
\end{multline*}
Now $ \sigma^{-r -(k-1)(1-\kappa) } \leq R_{0}^{-r }
j^{-r } $ and
$$ 
\| e^{(\re \mu_{\ell i} + |\tilde{\mu}_{0}|) \rho} \rho^{r}
\|_{L^{1}(\R^{+})} = \frac{r!}{\left |\re \mu_{\ell i} + |\tilde{\mu}_{0}|
  \, \right |^{r+1}  } .
$$
Hence, using Young inequality, we obtain
\begin{multline*}
\| w_{j}(\rho)  \rho^{1-\kappa} I_{\ell i}( \rho^{-k} f(\rho)) \|_{A}^{2}
\\
\leq C_{3}
\sum_{r=0}^{j} 2^{r+1}  \frac{1}{\left |\re \mu_{\ell i} + |\tilde{\mu}_{0}|
    \, \right |^{2(r+1)} }
\left(\binom{j}{r} \frac{r!}{ R_{0}^{r} j^{r}}
\right)^{2}
\| w_{j-k}(\rho) f(\rho) \|_{A}^{2} .
\end{multline*}
Remark that, for $ r > 0 $,
$$ 
\binom{j}{r} \frac{r!}{ R_{0}^{r} j^{r}} = \frac{1}{R_{0}^{r}} \left(1 -
\frac{1}{j}\right) \cdots \left( 1 - \frac{r-1}{j}\right) \leq
\frac{1}{R_{0}^{r}} 
$$
and that
$$ 
\frac{1}{\left |\re \mu_{\ell i} + |\tilde{\mu}_{0}|
  \, \right |^{r+1}} \leq \left(\frac{C_{4}}{\re \mu_{\ell
    i}}\right)^{r+1} \leq  C_{5}^{r+1} \mu_{\ell}^{- \frac{r+1}{2a}} ,
$$
we obtain the conclusion of the lemma recalling \eqref{eq:mitlda0},
provided
$$
R_{0} \geq \max \{4 C_{5} \mu_{\ell}^{-\frac{1}{2a}} , 2 \} .
$$
\end{proof}
Using the above lemmas \eqref{eq:uj-4} can be bound as
\begin{multline*}
\| w_{j}(\rho) \rho^{1-\kappa} \Dslash^{\gamma_{1}}_{\rho}
P_{0}^{-1} \left( \frac{1}{\rho^{k} } (1 -  \pi) g_{jk}  \right )
\|_{0, A}
\\
\leq  2 a C_{1} C_{0}  
\left( \sum_{\ell \geq 1} \| w_{j - k}(\rho)
 \mu_{\ell}^{\frac{\gamma_{1}}{2a} -1}  \langle g_{j k} ,
 \phi_{\ell}\rangle \|_{A}^{2} \right)^{\frac{1}{2}}  ,
\end{multline*}
where $ g_{jk} = Q^{\nu +r} P_{k} u_{j-k} $.

We have $ \mu_{\ell}^{\frac{\gamma_{1}}{2a} -1}  \langle g_{j k} ,
 \phi_{\ell}\rangle = \langle Q^{\nu + \frac{\gamma}{2a} - 1} P_{k}
 u_{j-k} , \phi_{\ell} \rangle $, see Lemma \ref{lemma:drhogamma}.

 Since the $ \phi_{\ell} $ are an orthonormal basis we obtain
\begin{multline*}
\| w_{j}(\rho) \rho^{1-\kappa} \Dslash^{\gamma_{1}}_{\rho} 
P_{0}^{-1} \left( \frac{1}{\rho^{k} } (1 -  \pi) Q^{\nu +r} P_{k} u_{j-k}   \right )
\|_{0, A}
\\
\leq C_{6}  
\| w_{j - k}(\rho)  Q^{\nu + \frac{\gamma}{2a} - 1} P_{k} u_{j-k} \|_{0,A} ,
\end{multline*}
As a consequence the last term on the r.h.s. of
\eqref{eq:uj-3} is estimated by
$$
\Lambda_{r} = 
C_{6} \sum_{k=1}^{\min\{j, 2a\}} \| w_{j - k}(\rho)  Q^{\nu +
  \frac{\gamma}{2a} - 1} P_{k} u_{j-k} \|_{0, A} .
$$
Recalling \eqref{eq:Pj}, \eqref{eq:Pjtilda} we write
$$
  \Lambda_{r}
  \leq C_{6} \sum_{k=1}^{\min\{j, 2a\}} \sum_{m=0}^{k} |p_{k
    m}|
  \| w_{j - k}(\rho)  \partial_{\rho}^{2a-k} Q^{\nu +
  \frac{\gamma}{2a} - 1} (x \partial_{x})^{m} u_{j-k} \|_{0,A} .
$$
Applying Proposition \ref{prop:AjD} and keeping in mind that the $
\rho $-derivative is conserved, we have to estimate
\begin{multline}
  \label{eq:Lambdar-2}
  \Lambda_{r} \leq C_{7} \sum_{k=1}^{\min\{j, 2a\}} \sum_{m=0}^{k}
  \Big(
  \| w_{j - k}(\rho)  \partial_{\rho}^{2a-k} Q^{\nu +
    \frac{\gamma}{2a} - 1 + m \frac{q}{2(q-1)} }  u_{j-k} \|_{0,A}
\\
  + \left( \nu +\frac{\gamma}{2a} - 1\right)^{2\frac{q-1}{q} (\nu +
  \frac{\gamma}{2a} - 1) + m } \| w_{j - k}(\rho)
\partial_{\rho}^{2a-k} u_{j-k} \|_{0,A} \Big)
\\
\leq
C_{8}
\sum_{k=1}^{\min\{j, 2a\}}
  \Big(
  \| w_{j - k}(\rho)  \partial_{\rho}^{2a-k} Q^{\nu +
    \frac{\gamma}{2a} - 1 + k \frac{q}{2(q-1)} }  u_{j-k} \|_{0,A}
\\
  + \left( \nu +\frac{\gamma}{2a} - 1\right)^{2\frac{q-1}{q} (\nu +
  \frac{\gamma}{2a} - 1) + k } \| w_{j - k}(\rho)
\partial_{\rho}^{2a-k} u_{j-k} \|_{0,A} \Big) , 
\end{multline}
where we estimated the coefficients of the polynomials $ \tilde{P}_{k}
$, $ p_{km} $, with a uniform constant and used the
estimate
$$
\|Q^{(m-k) \frac{q}{2(q-1)}} v \|_{0} \leq C \| v\|_{0} ,
$$
for $ 0 \leq m < k \leq 2a $. 

We also remark that, since $ A = A(j) = R_{0} (j+1) $, enlarging the
domain of integration we have norms over the half line $ ] A(j-k), +
\infty[ $, to which we may apply our inductive hypothesis. In fact the
exponent of $ Q $ is evidently $ \geq -1 $ and
$$ 
\nu + \frac{\gamma}{2a} - 1 + k \frac{q}{2(q-1)} + \frac{2a-k}{2a} \geq
\frac{k}{2} \frac{q}{q-1} \frac{1}{s_{0}} > 0.
$$
Thus, choosing $ \gamma^{\#} \geq 2a $,
\begin{multline*}
\| w_{j - k}(\rho)  \partial_{\rho}^{2a-k} Q^{\nu +
  \frac{\gamma}{2a} - 1 + k \frac{q}{2(q-1)} }  u_{j-k} \|_{0,A}
\\
\leq
C_{0}^{1+\sigma(j-k) +\sigma'(\nu + \frac{\gamma}{2a} - 1 + k
  \frac{q}{2(q-1)} + 2a - k)} (\lambda+1)^{\lambda} ,
\end{multline*}
where
\begin{multline*}
\lambda = \frac{j-k}{s_{0}} + \left( 2 \nu + \frac{\gamma}{a} - 2 + k
  \frac{q}{q-1} + 2 - \frac{k}{a}\right) \frac{q-1}{q}
\\
=
\frac{j}{s_{0}} + \left( 2\nu + \frac{\gamma}{a}\right) \frac{q-1}{q}
= \lambda(j, \nu, \gamma) .
\end{multline*}
As for the constant $ C_{0} $ we have
\begin{multline*}
\sigma(j-k) +\sigma'(\nu + \frac{\gamma}{2a} - 1 + k
\frac{q}{2(q-1)} + 2a - k)
\\
=
\sigma j + \sigma'(\nu + \frac{\gamma}{2a}) - \sigma k
+ \sigma' (2a-1) - \sigma' k \frac{q-2}{2(q-1)} 
\\
\leq
\sigma j + \sigma'(\nu + \gamma) - \sigma',
\end{multline*}
if we choose
\begin{equation}
  \label{eq:sigma}
\sigma = 2a \sigma'.
\end{equation}
An analogous computation can be made for the summands of the second
type in \eqref{eq:Lambdar-2}.

This completes the estimate of the second line on the right hand side
of \eqref{eq:uj-3}.

Next we are going to estimate the first line of the right hand side of
\eqref{eq:uj-3}:
\begin{multline*}
  \sum_{k=1}^{\min\{j, 2a\}}
  \Lambda_{k}
  = 
\sum_{k=1}^{\min\{j, 2a\}} \sum_{s=0}^{r-1} \| w_{j}(\rho) \rho^{1-\kappa}
\Dslash^{\gamma - 2a(1+s)}_{\rho}
\\
\left( \frac{1}{\rho^{k} } (1 - \pi)
  Q^{\nu +s} P_{k} u_{j-k} \right ) \|_{0, A} .
\end{multline*}
Since the sum over $ k $ is finite it is enough to consider just one
summand.
Observe that
\begin{multline*}
  \Lambda_{k}
  \leq C_{9}
\sum_{s=0}^{r-1}
  \sum_{m=0}^{\gamma - 2a(1+s)} \sum_{t=0}^{k}   \binom{\gamma - 2a(1+s)}{m}
  \frac{(k+m-1)!}{(k-1)!}
  \\
  \cdot
  \| w_{j}(\rho) \frac{1}{\rho^{k+m}} \rho^{1-\kappa} (1-\pi)
  \Dslash_{\rho}^{\gamma - 2a(1+s) - m + 2a - k} Q^{\nu+s} (x \partial_{x})^{t}
  u_{j-k} \|_{0,A}
  \\
  \leq
  C_{10}
  \sum_{s=0}^{r-1}
  \sum_{m=0}^{\gamma - 2a(1+s)} \sum_{t=0}^{k}   \binom{\gamma - 2a(1+s)}{m}
  \frac{(k+m-1)!}{(k-1)!} \frac{1}{A(j)^{m}}
  \\
  \cdot
  \| w_{j-k}(\rho) \rho^{-(k-1)(1-\kappa)}
  \partial_{\rho}^{\gamma - 2a s - m - k} Q^{\nu+s} (x \partial_{x})^{t}
  u_{j-k} \|_{0,A} 
   \\
  \leq
  C_{10}
  \sum_{s=0}^{r-1}
  \sum_{m=0}^{\gamma - 2a(1+s)} \sum_{t=0}^{k}   \binom{\gamma - 2a(1+s)}{m}
  \frac{(k+m-1)!}{(k-1)!} \frac{1}{A(j)^{m}}
  \\
  \cdot
  \| w_{j-k}(\rho)
  \partial_{\rho}^{\gamma - 2a s - m - k} Q^{\nu+s} (x \partial_{x})^{t}
  u_{j-k} \|_{0,A}
\end{multline*}
We apply Proposition \ref{prop:AjD} we obtain
\begin{multline}
  \label{eq:Lambdaks}
  \Lambda_{k}
  \leq
  C_{11}
  \sum_{s=0}^{r-1}
  \sum_{m=0}^{\gamma - 2a(1+s)} \sum_{t=0}^{k}   \binom{\gamma - 2a(1+s)}{m}
  \frac{(k+m-1)!}{(k-1)!}
  \\
  \cdot
  \frac{1}{A(j)^{m}}
  \Big( \|  w_{j-k}(\rho) \partial_{\rho}^{\gamma - 2a s - m -
  k} Q^{\nu+s+t \frac{q}{2(q-1)}} u_{j-k} \|_{0,A}
  \\
  +
  (\nu + s )^{2 \frac{q-1}{q} (\nu + s) + t} \| w_{j-k}(\rho)
  \partial_{\rho}^{\gamma - 2a s - m - k} u_{j-k} \|_{0,A} \Big)
  \\
  \leq
  C_{12}
  \sum_{s=0}^{r-1}
  \sum_{m=0}^{\gamma - 2a(1+s)}   \binom{\gamma - 2a(1+s)}{m}
  \frac{(k+m-1)!}{(k-1)!} \frac{1}{A(j)^{m}}
  \\
  \cdot
  \Big( \|  w_{j-k}(\rho) \partial_{\rho}^{\gamma - 2a s - m -
  k} Q^{\nu+s+k \frac{q}{2(q-1)}} u_{j-k} \|_{0,A}
  \\
  +
  (\nu + s )^{2 \frac{q-1}{q} (\nu + s) + k} \| w_{j-k}(\rho)
  \partial_{\rho}^{\gamma - 2a s - m - k} u_{j-k} \|_{0,A} \Big)
  .
\end{multline}
First we remark that, by \eqref{eq:AR},
\begin{multline*}
\binom{\gamma - 2a(1+s)}{m}
\frac{(k+m-1)!}{(k-1)!} \frac{1}{A(j)^{m}}
\\
=
\binom{k+m-1}{m} \frac{(\gamma - 2a(1+s))!}{(\gamma -2a(1+s) -m)!}
\frac{1}{(R_{0} (j+1))^{m}}
\\
\leq
\binom{k+m-1}{m} \frac{(\gamma - 2a(1+s)) \cdots (\gamma - 2a(1+s) - m
  + 1)}{R_{0}^{m} (j+1)^{m}}
\\
\leq
\binom{k+m-1}{m} \left(\frac{j+\gamma^{\#}}{R_{0}(j+1)}\right)^{m}
\leq 2^{2a-1} \left( \frac{2 \gamma^{\#}}{R_{0}}\right)^{m} ,
\end{multline*}
if $ R_{0} $ is suitably chosen.  
Moreover the quantities in the norms above verify the assumptions of
Theorem \ref{th:est'}, i.e. $ \nu \geq -1 $ and $ \nu + \gamma
(2a)^{-1} \geq 0 $. In fact the first is trivially true. As for the
second we have
\begin{multline*}
\nu + s + k \frac{q}{2(q-1)} + \frac{\gamma}{2a} - s - \frac{m}{2a} -
\frac{k}{2a}
\\
=
\nu + \frac{k}{2} \left( \frac{q}{q-1} - \frac{1}{a}\right) +
\frac{\gamma - m}{2a}
\\
\geq
\frac{k}{2} \frac{q}{q-1} \frac{1}{s_{0}} + \nu + 1 > 0 ,
\end{multline*}
since $ m \leq \gamma - 2a $ and $ \nu \geq -1 $. Moreover $ \gamma -
k \leq j - k + \gamma^{\#} $ if $ \gamma \leq j + \gamma^{\#} $. 

We may thus apply the inductive hypothesis to both norms in the right
hand side of \eqref{eq:Lambdaks}. Starting with the first on the next
to last line we have
\begin{multline*}
\|  w_{j-k}(\rho) \partial_{\rho}^{\gamma - 2a s - m -
  k} Q^{\nu+s+k \frac{q}{2(q-1)}} u_{j-k} \|_{0,A}
\\
\leq
C_{0}^{1+\sigma(j-k) + \sigma'(\nu+s+k \frac{q}{2(q-1)} +\gamma - 2a s - m -
  k)} (\lambda+1)^{\lambda} ,
\end{multline*}
where
\begin{multline*}
\lambda = \frac{j-k}{s_{0}} + \left( 2\nu + 2s + k \frac{q}{q-1} +
  \frac{\gamma}{a} - 2s - \frac{m}{a} - \frac{k}{a}\right)
\frac{q-1}{q} 
\\
=
\frac{j}{s_{0}} + \left( 2\nu + \frac{\gamma}{a}\right) \frac{q-1}{q}
- \frac{m}{a} \frac{q-1}{q} 
\\
\leq
\frac{j}{s_{0}} + \left( 2\nu + \frac{\gamma}{a}\right) \frac{q-1}{q}
= \lambda(j, \nu, \gamma) .
\end{multline*}
As for the exponent in the constant $ C_{0} $, we have
\begin{multline*}
\sigma(j-k) + \sigma'\left(\nu+s+k \frac{q}{2(q-1)} +\gamma - 2a s - m -
  k\right)
\\
=
\sigma j + \sigma' (\nu + \gamma) - \sigma k + \sigma' k
  \frac{q}{2(q-1)} - \sigma' s (2a-1) - \sigma' m
\\
\leq
\sigma j + \sigma' (\nu + \gamma) - \sigma' \left( 2a -
 1\right) - \sigma' s (2a-1)
,
\end{multline*}
since $ k \geq 1 $ and $ \sigma = 2a \sigma' $.

Consider now the second term in the right hand side of
\eqref{eq:Lambdaks}: by the inductive hypothesis
\begin{multline*}
(\nu + s )^{2 \frac{q-1}{q} (\nu + s) + k} \| w_{j-k}(\rho)
\partial_{\rho}^{\gamma - 2a s - m - k} u_{j-k} \|_{0,A}
\\
\leq
C_{0}^{1+\sigma(j-k) + \sigma'(\gamma-2as-m-k)} (\nu + s )^{2
  \frac{q-1}{q} (\nu + s) + k} (\lambda_{1}+1)^{\lambda_{1}},
\end{multline*}
where
$$ 
\lambda_{1} = \frac{j-k}{s_{0}} + \left(\frac{\gamma}{a} -2s
  -\frac{m}{a} - \frac{k}{a}\right) \frac{q-1}{q} .
$$
Both $ \nu + s $ and $ \lambda_{1} $ can be enlarged to be $
\lambda(j, \nu, \gamma) $, so that it is enough to check that the sum
of the exponents has the right value. As before this is
\begin{multline*}
\frac{j-k}{s_{0}} + \left(\frac{\gamma}{a} -2s
  -\frac{m}{a} - \frac{k}{a}\right) \frac{q-1}{q} + 2 \frac{q-1}{q}
(\nu + s) + k
\\
= \frac{j}{s_{0}} + \left( 2\nu + \frac{\gamma}{a}\right)
\frac{q-1}{q} - \frac{m}{a} \frac{q-1}{q} < \lambda(j, \nu, \gamma).
\end{multline*}
Hence
\begin{multline*}
  \Lambda_{k}
  \leq
   C_{12}
  \sum_{s=0}^{r-1}
  \sum_{m=0}^{\gamma - 2a(1+s)}   \binom{\gamma - 2a(1+s)}{m}
  \frac{(k+m-1)!}{(k-1)!} \frac{1}{A(j)^{m}}
  \\
  \cdot
  \Big( \|  w_{j-k}(\rho) \partial_{\rho}^{\gamma - 2a s - m -
  k} Q^{\nu+s+k \frac{q}{2(q-1)}} u_{j-k} \|_{0,A}
  \\
  +
  (\nu + s )^{2 \frac{q-1}{q} (\nu + s) + k} \| w_{j-k}(\rho)
  \partial_{\rho}^{\gamma - 2a s - m - k} u_{j-k} \|_{0,A} \Big)
  \\
  \leq
  C_{12} 2^{2a-1}
  \sum_{s=0}^{r-1}
  \sum_{m=0}^{\gamma - 2a(1+s)} \left(\frac{2 \gamma^{\#}}{R_{0}}\right)^{m}
  \\
  \cdot
  \Big( C_{0}^{1 + \sigma j + \sigma' (\nu + \gamma) - \sigma' \left( 2a -
 1\right) - \sigma' s (2a-1)} (\lambda(j, \nu, \gamma)+1)^{\lambda(j, \nu,
\gamma)}
\\
+ C_{0}^{1+\sigma j - 2a \sigma' + \sigma'(\gamma + \nu) - (2a-1)
  \sigma' s}
(\lambda(j, \nu, \gamma)+1)^{\lambda(j, \nu,
  \gamma)} \Big)
\\
\leq
  C_{0}^{1 + \sigma j + \sigma' (\nu + \gamma)} (\lambda(j, \nu, \gamma)+1)^{\lambda(j, \nu,
    \gamma)}
  \\
  \cdot
  C_{12} 2^{2a-1} C_{0}^{-\sigma'(2a-1)} 
  \sum_{s=0}^{\infty} C_{0}^{-\sigma'(2a-1)s}
  \sum_{m=0}^{\infty}  \left(
    \frac{2 \gamma^{\#}}{R_{0}}\right)^{m} .
\end{multline*}
Keeping into account that $ k $ ranges on a finite number of indices,
a suitable choice of both $ C_{0} $ and $ R_{0} $ completes the proof
of Theorem \ref{th:est-1-p}.
\end{proof}

\section{Estimate of $ \pi u_{j} $}
\setcounter{equation}{0}
\setcounter{theorem}{0}
\setcounter{proposition}{0}
\setcounter{lemma}{0}
\setcounter{corollary}{0}
\setcounter{definition}{0}

This section is devoted to proving estimates, analogous to those in
\eqref{eq:est-1-p}, for $ \pi u_{j} $. We use the same notation of the
preceding section. 

\begin{theorem}
\label{th:est-p}
Let $ \nu $ denote a rational number $ \geq -1 $ and $ \gamma \in \N $
such that $ \nu + \gamma (2a)^{-1} \geq 0 $.
There exist positive constants $ C_{0} $, $ R_{0} $, $ \sigma $, $
\sigma' $ such that 
\begin{equation}
\label{eq:est-p}
\| w_{j}(\rho) \partial_{\rho}^{\gamma} Q^{\nu}
\pi  u_{j}(x, \rho)  \|_{0, A}
\leq
C_{0}^{1 + \sigma j + \sigma'( \nu + \gamma)} 
(\lambda(j, \nu, \gamma)+1)^{\lambda(j, \nu, \gamma)} ,
\end{equation}
where
$$ 
\lambda(j, \nu, \gamma) = \frac{j}{s_{0}} + \left(2 \nu +
  \frac{\gamma}{a}\right) \frac{q-1}{q},
$$
and
$ A = A(j) = R_{0}(j + 1)  $, $ \gamma \leq j+\gamma^{\#} $, $ 2 \gamma^{\#} <
R_{0} $. 
\end{theorem}
\begin{proof}
Recall that from \eqref{eq:tr2} $ \pi u_{j} $ is obtained as a
solution of 
\begin{align}
  \label{eq:pi-u}
\pi P_{0} u_{j} &=
P_{0} \pi u_{j} = - \frac{1}{\rho}
\pi P_{1} (1 - \pi) u_{j} - \frac{1}{\rho} \pi
                          P_{1} \pi u_{j-1}
\\
  &\phantom{= =}    - \sum_{k=1}^{\min\{j, 2a-1\}} \frac{1}{\rho^{k+1}}
    \pi P_{k+1} u_{j-k} . \notag
\end{align}
First of all we observe that the decay rate in $ \rho $ of the first and third
terms on the right hand side of \eqref{eq:pi-u} is
$$
\rho^{-[(j-\delta)\kappa + (1-\kappa) +1]} , \qquad \rho^{-[(j-\delta)\kappa +
  ( 1- \kappa) k + 1]}
$$
respectively, while for the second term we have $ \rho^{-[(j -\delta)\kappa + 1 -
  \kappa]} $ which would not allow us to make an inductive argument.
\begin{lemma}
  \label{lemma:piPpi}
  We have
\begin{equation}
\label{eq:piPpi}
\pi P_{1} \pi v = 0.
\end{equation}
\end{lemma}
\begin{proof}
We just remark that $ \pi P_{1} \pi v = \pi P_{1} \left( \langle v,
  \phi_{0} \rangle \phi_{0}\right) $. Using the notation of
Proposition \ref{prop:r} we may write
\begin{eqnarray*}
  \pi P_{1} \left( \langle v, \phi_{0} \rangle \phi_{0}\right)
&  = &
  \pi \left( \left[ \alpha x \partial_{x} + \beta + 2a\gamma_{2a} r\right]
    \partial_{\rho}^{2a-1} \left( \langle v,
      \phi_{0} \rangle \phi_{0}\right) \right)
  \\
&  = & \pi \left( \left( \partial_{\rho}^{2a-1} \langle v, \phi_{0}
      \rangle \right) \left[ \alpha x \partial_{x} + \beta + 2a\gamma_{2a}
      r\right] \phi_{0} \right)
  \\
&  = & \left( \partial_{\rho}^{2a-1} \langle v, \phi_{0} \rangle \right)
  \left[ \langle \alpha x \partial_{x} \phi_{0} , \phi_{0}\rangle  +
    \beta + 2a\gamma_{2a} r \right] \phi_{0} = 0,
\end{eqnarray*}
due to the choice of $ r $ in Proposition \ref{prop:r}.
\end{proof}
As a consequence \eqref{eq:pi-u} becomes
\begin{equation}
  \label{eq:pi-u2}
P_{0} \pi u_{j} = - \frac{1}{\rho}
\pi P_{1} (1 - \pi) u_{j}
- \sum_{k=1}^{\min\{j, 2a-1\}} \frac{1}{\rho^{k+1}}
    \pi P_{k+1} u_{j-k} .
\end{equation}
Let us write
\begin{equation}
\label{eq:util}
\tilde{u}_{\ell}^{(j)} =
\begin{cases}
  (1 - \pi) u_{j} & \textrm{if $ \ell = j $} \\
  u_{\ell}        & \textrm{if $ \ell < j $}
\end{cases}
.
\end{equation}
Hence \eqref{eq:pi-u2} becomes
\begin{equation}
  \label{eq:pi-u3}
P_{0} \pi u_{j} =
- \sum_{k=0}^{\min\{j, 2a-1\}} \frac{1}{\rho^{k+1}}
    \pi P_{k+1} \tilde{u}_{j-k}^{(j)} .
\end{equation}
By the inductive hypothesis and by Theorem \ref{th:est-1-p} the
functions $ \tilde{u}_{j-k}^{(j)} $ verify the estimates
\eqref{eq:est'} and \eqref{eq:est-1-p} for $ k=0 $. 

To start with, by Lemma \ref{lemma:QnuP0},
\begin{multline}
\label{eq:puj-1}
\| w_{j}(\rho) \partial_{\rho}^{\gamma} Q^{\nu} \pi u_{j}(x, \rho)
\|_{0,A}
\\
\leq
\sum_{k=0}^{\min\{j, 2a-1\}} \| w_{j}(\rho) \partial_{\rho}^{\gamma}
P_{0}^{-1} \left( \frac{1}{\rho^{k+1}} \pi Q^{\nu} P_{k+1}
  \tilde{u}_{j-k}^{(j)} \right) \|_{0,A} .
\end{multline}
By Corollary \ref{cor:drhogammaP0} we have
\begin{multline}
\label{eq:puj-2}
\| w_{j}(\rho) \Dslash_{\rho}^{\gamma} Q^{\nu}
\pi u_{j}(x, \rho)  \|_{0, A}
\\
\leq
\sum_{k=0}^{\min\{j, 2a-1\}} \sum_{s=0}^{r-1} \| w_{j}(\rho)
\Dslash^{\gamma - 2a(1+s)}_{\rho} \left( \frac{1}{\rho^{k+1} }  \pi
  Q^{\nu +s} P_{k+1} \tilde{u}_{j-k}^{(j)} \right ) \|_{0, A}
\\
+
\sum_{k=0}^{\min\{j, 2a-1\}} \| w_{j}(\rho) \Dslash^{\gamma_{1}}_{\rho}
P_{0}^{-1} \left( \frac{1}{\rho^{k+1} }   \pi
  Q^{\nu +r} P_{k+1} \tilde{u}_{j-k}^{(j)} \right ) \|_{0, A} .
\end{multline}
Consider the terms on the third line of the above inequality.
Denote by $ \tilde{g}_{jk} = Q^{\nu +r} P_{k+1} \tilde{u}_{j-k}^{(j)}
$. Then for $ k \in \{0, \ldots, \min\{j, 2a-1\}\}$,
\begin{multline}
  \label{eq:dgamma1}
\| w_{j}(\rho) \Dslash^{\gamma_{1}}_{\rho}
P_{0}^{-1} \left( \frac{1}{\rho^{k+1} }  \pi \tilde{g}_{jk}  \right )
\|_{0, A}
\\
=
\| w_{j}(\rho) \Dslash^{\gamma_{1}}_{\rho} 
\sum_{i=1}^{2a} A_{0 i} I_{0 i}(\rho^{-(k+1)} \langle \tilde{g}_{j k} ,
\phi_{0} \rangle ) \phi_{0} \|_{0, A}
\\
=
\| w_{j}(\rho) \Dslash^{\gamma_{1}}_{\rho}
  \sum_{i=1}^{2a} A_{0 i} I_{0 i}(\rho^{-(k+1)} \langle \tilde{g}_{j k} , 
  \phi_{0} \rangle ) \|_{A} ,
\end{multline}
Preliminarily we need the analogous of lemmas \ref{lemma:muli>0},
\ref{lemma:muli<0} for the fundamental eigenfunction, corresponding to
the projection $ \pi $.

\begin{lemma}
\label{lemma:pi-muli>0}
Assume that $ \re \mu_{0 i} > 0 $. Then
\begin{equation}
\label{eq:pi-muli>0}
\| w_{j}(\rho)  I_{0 i}( \rho^{-k} f(\rho)) \|_{A}^{2} \leq \tilde{C}_{0} \mu_{0}^{-\frac{1}{a}} \|
w_{j-k}(\rho) f(\rho) \|_{A}^{2} ,
\end{equation}
for a suitable constant $ \tilde{C}_{0} > 0 $ independent of $ j $,
$ k $, $ i $, where we choose $ R = A(j) $, see \eqref{eq:Ikjgl}.
\end{lemma}
\begin{lemma}
\label{lemma:pi-muli<0}
Assume that $ \re \mu_{0 i} < \tilde{\mu}_{0}  < 0 $. Then
\begin{equation}
\label{eq:pi-muli<0}
\| w_{j}(\rho) I_{0 i}( \rho^{-k} f(\rho)) \|_{A}^{2} \leq \tilde{C}_{0} \mu_{0}^{-\frac{1}{a}} \|
w_{j-k}(\rho)  f(\rho) \|_{A}^{2} ,
\end{equation}
for a suitable constant $ \tilde{C}_{0} > 0 $ independent of $ j $,
$ k $, $ i $, and $ R = A(j)$.
\end{lemma}
Under the assumptions of lemmas \ref{lemma:pi-muli<0},
\ref{lemma:pi-muli>0}, $ \re \mu_{0i} - \tilde{\mu}_{0} \neq 0 $ and,
as a consequence, their proofs are completely analogous to those of
lemmas \ref{lemma:muli<0}, \ref{lemma:muli>0}, since the
factor $ \rho^{1-\kappa} $ plays no role. 
\begin{lemma}
\label{lemma:pi-mu0}
Assume that $ \re \mu_{0 i} = \tilde{\mu}_{0}  < 0 $ and $ 1 \leq k
\leq \min\{j, 2a - 1\}$. Then 
\begin{equation}
\label{eq:pi-mu0}
\| w_{j}(\rho) I_{0 i}( \rho^{-(k+1)} f(\rho)) \|_{A}^{2} \leq \tilde{C}_{0}^{2}
\frac{1}{j A(j)^{2(1-\kappa)}}  \|
w_{j-k}(\rho)  f(\rho) \|_{A}^{2} ,
\end{equation}
for a suitable constant $ \tilde{C}_{0} > 0 $ independent of $ j $, and $ R =
A(j)$.
\end{lemma}
\begin{proof}[Proof of Lemma \ref{lemma:pi-mu0}]
We have to estimate, using \eqref{eq:Ikjgl},
\begin{multline*}
  \| w_{j}(\rho) I_{0 i}( \rho^{-(k+1)} f(\rho)) \|_{A}^{2}
 =
\int_{A}^{+\infty} \rho^{2(j-\delta)\kappa} e^{-2 \tilde{\mu}_{0}
  \rho} \Big| \int_{\R} e^{\mu_{0i} (\rho - \sigma)}
\\
\cdot
H\left(-\sgn
  \left( \re \mu_{0i} - \tilde{\mu}_{0} + \epsilon_{\mu} \right) (\rho
  -\sigma) \right)
\frac{H(\sigma - A)}{\sigma^{k+1}} f(\sigma) d \sigma \Big|^{2} d\rho .
\end{multline*}
Since $  \re \mu_{0i} - \tilde{\mu}_{0} = 0 $ and  $ \epsilon_{\mu} >
0 $, we obtain
\begin{multline*}
  \| w_{j}(\rho) I_{0 i}( \rho^{-(k+1)} f(\rho)) \|_{A}^{2}
\\
 =
\int_{A}^{+\infty} \rho^{2(j-\delta)\kappa} e^{-2 \tilde{\mu}_{0}
  \rho} \Big| \int_{\R} e^{\mu_{0i} (\rho - \sigma)}
H\left(-  (\rho -\sigma) \right)
\frac{1}{\sigma^{k+1}} H(\sigma - A) f(\sigma) d \sigma \Big|^{2} d\rho
\\
\leq
\int_{A}^{+\infty} \rho^{2(j-\delta)\kappa} e^{-2 \tilde{\mu}_{0}
  \rho} \Big( \int_{\R} e^{\tilde{\mu}_{0} (\rho - \sigma)}
H\left(- (\rho -\sigma) \right)
\frac{1}{\sigma^{k+1}} H(\sigma - A) |f(\sigma)| d \sigma \Big)^{2}
d\rho
\\
=
\int_{A}^{+\infty} \rho^{2(j-\delta)\kappa} \Big( \int_{\rho}^{+\infty}
\frac{1}{\sigma^{k+1}} H(\sigma - A) e^{- \tilde{\mu}_{0} \sigma}
|f(\sigma)| d \sigma \Big)^{2}  d\rho
\\
=
\int_{A}^{+\infty} \rho^{2(j-\delta)\kappa}
\Big( \int_{\rho}^{+\infty}
\frac{H(\sigma - A)}{\sigma^{(j-\delta)\kappa + k(1-\kappa)+1}}   \sigma^{(j-k-\delta)\kappa}
e^{- \tilde{\mu}_{0} \sigma}
|f(\sigma)| d \sigma \Big)^{2}  d\rho . 
\end{multline*}
Using H\"older inequality on the inner integral above we get
\begin{multline*}
 \| w_{j}(\rho) I_{0 i}( \rho^{-(k+1)} f(\rho)) \|_{A}^{2}
\\
 \leq
  \| w_{j-k}(\rho) f(\rho) \|_{A}^{2}
 \int_{A}^{+\infty} \rho^{2(j-\delta)\kappa} \int_{\rho}^{+\infty}
 \frac{1}{\sigma^{2(j-\delta)\kappa + 2 k (1-\kappa)+2}} d\sigma
 \\
 =
 L_{j k} 
 \| w_{j-k}(\rho) f(\rho) \|_{A}^{2}  \int_{A}^{+\infty}
 \frac{1}{\rho^{2k(1-\kappa) + 1}} d\rho 
\end{multline*}
where we set
$$ 
\int_{\rho}^{+\infty} \frac{1}{\sigma^{2(j-\delta)\kappa + 2 k
    (1-\kappa) +2}} d\sigma = \frac{L_{j k}}{\rho^{2(j-\delta)\kappa + 2k
  (1-\kappa) + 1} } .
$$
The fact that $ k \geq 1 $ allows us to conclude the proof of the lemma.
\end{proof}
Finally we need one more result to estimate the first term in
\eqref{eq:pi-u3}, i.e. the term corresponding to $ k = 0 $ when $ \re
\mu_{0 i} = \tilde{\mu}_{0} $.
\begin{lemma}
\label{lemma:piP-1-pi}
Assume that $ \re \mu_{0 i} = \tilde{\mu}_{0}  < 0 $. 
\begin{equation}
\label{eq:piP-1-pi}
\| w_{j}(\rho) I_{0 i}( \rho^{-1} f(\rho)) \|_{A}^{2} \leq
\frac{C_{0}}{j A(j)^{2(1-\kappa)}} \|
w_{j}(\rho) \rho^{1-\kappa} f(\rho) \|_{A}^{2} ,
\end{equation}
for a suitable constant $ C_{0} > 0 $ independent of $ j $, and $ R =
A(j)$.
\end{lemma}
\begin{proof}[Proof of Lemma \ref{lemma:piP-1-pi}]
Arguing as in the beginning of the proof of Lemma \ref{lemma:pi-mu0}
we have
\begin{multline*}
  \| w_{j}(\rho) I_{0 i}( \rho^{-1} f(\rho)) \|_{A}^{2}
  =
  \int_{A}^{+\infty} \rho^{2(j-\delta)\kappa} e^{-2\tilde{\mu}_{0} \rho}
  \Big | \int_{\R} e^{\mu_{0i} (\rho - \sigma)}
\\
\cdot
H\left(-\sgn
  \left( \re \mu_{0i} - \tilde{\mu}_{0} + \epsilon_{\mu} \right) (\rho
  -\sigma) \right)
\frac{H(\sigma - A)}{\sigma} f(\sigma) d \sigma \Big|^{2} d\rho
\\
\leq
\int_{A}^{+\infty} \rho^{2(j-\delta)\kappa} \Big( \int_{\rho}^{+\infty}
\frac{1}{\sigma} H(\sigma - A) e^{- \tilde{\mu}_{0} \sigma}
|f(\sigma)| d \sigma \Big)^{2}  d\rho
\\
=
\int_{A}^{+\infty} \rho^{2(j-\delta)\kappa}
\Big( \int_{\rho}^{+\infty}
\frac{H(\sigma - A)}{\sigma^{(j-\delta) \kappa + 2 - \kappa}} \rho^{(j-\delta)\kappa + 1 - \kappa} e^{- \tilde{\mu}_{0} \sigma}
|f(\sigma)| d \sigma \Big)^{2}  d\rho
\end{multline*}
Using H\"older inequality on the inner integral above we get
\begin{multline*}
  \| w_{j}(\rho) I_{0 i}( \rho^{-1} f(\rho)) \|_{A}^{2}
  \\
 \leq
  \| w_{j}(\rho) \rho^{1-\kappa}  f(\rho) \|_{A}^{2}
 \int_{A}^{+\infty} \rho^{2(j-\delta)\kappa} \int_{\rho}^{+\infty}
 \frac{1}{\sigma^{2(j-\delta)\kappa + 2  (2-\kappa) }} d\sigma
 \\
 =
 L_{j} 
 \| w_{j}(\rho)\rho^{1-\kappa}   f(\rho) \|_{A}^{2}  \int_{A}^{+\infty}
 \frac{1}{\rho^{1 +2(1 - \kappa)}} d\rho 
\end{multline*}
where we set
$$ 
\int_{\rho}^{+\infty} \frac{1}{\sigma^{2(j-\delta)\kappa + 2
    (2-\kappa)}} d\sigma = \frac{L_{j}}{\rho^{2(j-\delta)\kappa + 
  2(1-\kappa)+1} } .
$$
This concludes the proof of the lemma.
\end{proof}
Going back to \eqref{eq:dgamma1} and using \eqref{eq:Drhogamma1} we
have
\begin{multline}
  \label{eq:gamma1}
\| w_{j}(\rho) \Dslash^{\gamma_{1}}_{\rho}
P_{0}^{-1}\left(\frac{1}{\rho^{k+1}} \pi \tilde{g}_{jk}\right) \|_{0,A} 
\\
=
\| w_{j}(\rho) \Dslash^{\gamma_{1}}_{\rho}
  \sum_{i=1}^{2a} A_{0 i} I_{0 i}(\rho^{-(k+1)} \langle \tilde{g}_{j k} , 
  \phi_{0} \rangle ) \|_{A} 
  \\
=
s_{0}^{-\gamma_{1}} 
\| w_{j}(\rho)
  \sum_{i=1}^{2a} A_{0 i} \mu_{0 i}^{\gamma_{1}}  I_{0 i}(\rho^{-(k+1)} \langle \tilde{g}_{j k} , 
  \phi_{0} \rangle ) \|_{A}
  \\
  \leq
C_{1}
\sum_{i=1}^{2a}  \| w_{j}(\rho)
    I_{0 i}(\rho^{-(k+1)} \langle \tilde{g}_{j k} ,
    \mu_{0}^{\frac{\gamma_{1} - 2a + 1}{2a}} \phi_{0} \rangle ) \|_{A} .
\end{multline}
To the last quantity we apply Lemmas \ref{lemma:pi-muli>0},
\ref{lemma:pi-muli<0} when $ \re \mu_{0i} \neq \tilde{\mu}_{0} $,
Lemma \ref{lemma:pi-mu0} when $ \re \mu_{0i} = \tilde{\mu}_{0} $ and $
1 \leq k \leq \min\{ j, 2a\} $ and Lemma \ref{lemma:piP-1-pi} when $
\re \mu_{0i} = \tilde{\mu}_{0} $ and $ k = 0 $. In the first case,
i.e. when $ \re \mu_{0i} \neq \tilde{\mu}_{0} $, we obtain a decay rate
better than that we get when $ \re \mu_{0i} = \tilde{\mu}_{0} $ so
that the gain must be neglected.

We also remark that there is no need to keep track of the precise
powers of $ \mu_{0} $ in the formula above because they may always be
absorbed by a constant. We did that only to have a more symmetric
formula. 

Consider first the terms where $ k = 0, 1, \ldots, \min\{j, 2a-1\} $ and
$ \re \mu_{0i}$  $ \neq \tilde{\mu}_{0}$. Applying Lemmas
\ref{lemma:pi-muli>0}, \ref{lemma:pi-muli<0} in \eqref{eq:gamma1}, we
get (see \eqref{eq:Pj})
\begin{multline*}
\sum_{\re \mu_{0i} \neq \tilde{\mu}_{0}} \| w_{j}(\rho)
    I_{0 i}(\rho^{-(k+1)} \langle \tilde{g}_{j k} ,
    \mu_{0}^{\frac{\gamma_{1} -2a + 1}{2a} } \phi_{0} \rangle
    ) \|_{A}
    \\
    \leq
 C_{2} \| w_{j-k-1}(\rho) \langle \tilde{g}_{j k} ,
    \mu_{0}^{\frac{\gamma_{1}}{2a} - 1} \phi_{0} \rangle \|_{A}
    \\
=
 C_{2} \| w_{j-k-1}(\rho) \langle Q^{\nu + r} P_{k+1} \tilde{u}_{j-k}^{(j)} ,
    \mu_{0}^{\frac{\gamma_{1}}{2a} - 1} \phi_{0} \rangle \|_{A}    
\\
=
C_{2} \| w_{j-k-1}(\rho) \partial_{\rho}^{2a-k-1} \langle  \tilde{u}_{j-k}^{(j)} ,
    \mu_{0}^{\nu + \frac{\gamma}{2a}-1} \ {}^{t}\tilde{P}_{k+1} \phi_{0} \rangle \|_{A} 
\\
=
C_{2} \mu_{0}^{\nu + \frac{\gamma}{2a}-1}
\| w_{j-k-1}(\rho)
\partial_{\rho}^{2a-k-1} \langle  Q^{- \frac{2a-k-1}{2a}} \tilde{u}_{j-k}^{(j)} ,  Q^{\frac{2a-k-1}{2a}}  \ {}^{t}\tilde{P}_{k+1} \phi_{0} \rangle \|_{A} 
\\
\leq
\frac{C_{2}\mu_{0}^{\nu + \frac{\gamma}{2a}-1}}{A(j)^{\kappa}}
\|Q^{\frac{2a-k-1}{2a}} \ {}^{t} \tilde{P}_{k+1} \phi_{0} \|_{0}
\| w_{j-k}(\rho)\partial_{\rho}^{2a-k-1} Q^{- \frac{2a-k-1}{2a}} \tilde{u}_{j-k}^{(j)} \|_{0,A}
\end{multline*}
We remark that the norm $ \|Q^{\frac{2a-k-1}{2a}}  \ {}^{t}\tilde{P}_{k+1}
\phi_{0} \|_{0} $ is an absolute constant since $ 0 \leq k \leq 2a-1
$. Furthermore we see that the conditions $ \nu \geq -1 $ and $ \nu +
\frac{\gamma}{2a} \geq 0 $, when $ \nu $, $ \gamma $ are the exponents
of $ Q $, $ \partial_{\rho} $ respectively, are satisfied; moreover $
2a - k -1 \leq j-k+\gamma^{\#} $ and hence we may
apply the inductive hypothesis when $ k \geq 1 $ and the result of the
preceding section for $ (1-\pi)u_{j} $ when $ k = 0 $, thus obtaining
the bound 
\begin{multline*}
\sum_{k=0}^{\min\{j, 2a-1\}} \sum_{\re \mu_{0i} \neq \tilde{\mu}_{0}} \| w_{j}(\rho)
    I_{0 i}(\rho^{-(k+1)} \langle \tilde{g}_{j k} ,
    \mu_{0}^{\frac{\gamma_{1} -2a + 1}{2a} } \phi_{0} \rangle
    ) \|_{A}
\\
\leq
\frac{C_{3}\mu_{0}^{\nu + \frac{\gamma}{2a}}}{A(j)^{\kappa}}
\sum_{k=0}^{\min\{j, 2a-1\}} \| w_{j-k}(\rho)\partial_{\rho}^{2a-k-1}
Q^{- \frac{2a-k-1}{2a}} \tilde{u}_{j-k}^{(j)} \|_{0,A}
\\
\leq
\frac{C_{3}\mu_{0}^{\nu + \frac{\gamma}{2a}}}{A(j)}
\| w_{j}(\rho) \rho^{1-\kappa} \partial_{\rho}^{2a-1}
Q^{- \frac{2a-1}{2a}} (1-\pi)u_{j} \|_{0,A}
\\
+
\frac{C_{3}\mu_{0}^{\nu + \frac{\gamma}{2a}}}{A(j)^{\kappa}}
\sum_{k=1}^{\min\{j, 2a-1\}} \| w_{j-k}(\rho)\partial_{\rho}^{2a-k-1}
Q^{- \frac{2a-k-1}{2a}} u_{j-k} \|_{0,A}
\\
\leq
\frac{C_{3}\mu_{0}^{\nu + \frac{\gamma}{2a}}}{A(j)^{\kappa}}
\sum_{k=0}^{\min\{j, 2a-1\}} C_{0}^{1 + \sigma(j-k)+ \sigma' (2a-k-1)
}
\\
\cdot
(\lambda(j-k, 0, 0) + 1)^{\lambda(j-k, 0, 0)} . 
\end{multline*}
Hence we get the estimate
\begin{multline}
\label{eq:k>0.}
\sum_{k=1}^{\min\{j, 2a-1\}} \| w_{j}(\rho) \Dslash^{\gamma_{1}}_{\rho}
  \sum_{\re \mu_{0i} \neq \tilde{\mu}_{0}} A_{0 i} I_{0 i}(\rho^{-(k+1)} \langle \tilde{g}_{j k} , 
  \phi_{0} \rangle ) \|_{A} 
  \\
  \leq
\frac{1}{2} C_{0}^{1 + \sigma j + \sigma'(\nu + \gamma) - \sigma'} (\lambda(j, \nu, \gamma)
+ 1)^{\lambda(j, \nu, \gamma)} ,
\end{multline}
when $ k \geq 1 $ and $ \re \mu_{0i} \neq \tilde{\mu}_{0} $, provided
we choose $ C_{0} $ large enough.

Let us now consider the case $ k=0 $ and $ \re \mu_{0i} \neq
\tilde{\mu}_{0} $. We have to estimate
\begin{multline*}
\| w_{j}(\rho) \Dslash^{\gamma_{1}}_{\rho}
  \sum_{\re \mu_{0i} \neq \tilde{\mu}_{0}} A_{0 i} I_{0 i}(\rho^{-1} \langle (1-\pi) u_{j} , 
  \phi_{0} \rangle ) \|_{A}
\\
\leq
\frac{C_{3}\mu_{0}^{\nu + \frac{\gamma}{2a}}}{A(j)^{\kappa}}
C_{0}^{1 + \sigma j+ \sigma' (2a-1)}
(\lambda(j, 0, 0) + 1)^{\lambda(j, 0, 0)}
\\
\leq
\frac{C_{3} C_{0}^{2a\sigma'} }{R_{0}^{\kappa}(j+1)^{\kappa}}
\frac{\mu_{0}^{\nu + \frac{\gamma}{2a}}}{ C_{0}^{\sigma'(\nu+\gamma)}
} \
C_{0}^{1 + \sigma j+ \sigma' (\nu+\gamma) - \sigma'}
(\lambda(j, \nu, \gamma) + 1)^{\lambda(j, \nu, \gamma)}
\end{multline*}
We make the choice $ C_{0}^{\sigma'} > \mu_{0} $ and $ R_{0}^{\kappa}
> 2 C_{3} C_{0}^{2a\sigma'} $, so that
\begin{multline}
\label{eq:kgeq0}
\sum_{k=0}^{\min\{j, 2a-1\}} \| w_{j}(\rho) \Dslash^{\gamma_{1}}_{\rho}
  \sum_{\re \mu_{0i} \neq \tilde{\mu}_{0}} A_{0 i} I_{0 i}(\rho^{-(k+1)} \langle \tilde{g}_{j k} , 
  \phi_{0} \rangle ) \|_{A} 
  \\
  \leq
C_{0}^{1 + \sigma j + \sigma'(\nu + \gamma) - \sigma'} (\lambda(j, \nu, \gamma)
+ 1)^{\lambda(j, \nu, \gamma)}
\end{multline}

\bigskip

So far we treated \eqref{eq:dgamma1} in the case when $ \re
\mu_{0i} \neq \tilde{\mu}_{0} $ for any value of $ k $. 

\bigskip

Consider now \eqref{eq:dgamma1}; we have to estimate
\eqref{eq:gamma1} when $ \re \mu_{0i} = \tilde{\mu}_{0} $.
\begin{equation}
  \label{eq:617}
 C_{1}  \sum_{k=0}^{\min\{j, 2a -1\}}
  \sum_{\re \mu_{0i} = \tilde{\mu}_{0}}  \| w_{j}(\rho)
    I_{0 i}(\rho^{-(k+1)} \langle \tilde{g}_{j k} ,
    \mu_{0}^{\frac{\gamma_{1} - 2a + 1}{2a}} \phi_{0} \rangle ) \|_{A} .
\end{equation}
Let us start considering the case when $ k=0 $:
$$ 
C_{1} \| w_{j}(\rho)
    I_{0 i}(\rho^{-1} \langle \tilde{g}_{j 0} ,
    \mu_{0}^{\frac{\gamma_{1} - 2a + 1}{2a}} \phi_{0} \rangle ) \|_{A} .
$$
Applying Lemma \ref{lemma:piP-1-pi} we get
\begin{multline*}
C_{1} \| w_{j}(\rho)
    I_{0 i}(\rho^{-1} \langle \tilde{g}_{j 0} ,
    \mu_{0}^{\frac{\gamma_{1} - 2a + 1}{2a}} \phi_{0} \rangle ) \|_{A}
    \\
\leq
C_{1} \frac{\tilde{C}_{0}}{j^{\frac{1}{2}} A(j)^{1-\kappa}} \| w_{j}(\rho) \rho^{1-\kappa}  \langle \tilde{g}_{j 0} ,
\mu_{0}^{\frac{\gamma_{1} - 2a + 1}{2a}} \phi_{0} \rangle \|_{A}
\\
=
C_{1} \frac{\tilde{C}_{0}}{j^{\frac{1}{2}} A(j)^{1-\kappa}}  \| w_{j}(\rho) \rho^{1-\kappa}  \langle Q^{\nu +r} P_{1} (1-\pi)u_{j} ,
\mu_{0}^{\frac{\gamma_{1} - 2a + 1}{2a}} \phi_{0} \rangle \|_{A}
\\
=
C_{1} \frac{\tilde{C}_{0}}{j^{\frac{1}{2}} A(j)^{1-\kappa}}  \| w_{j}(\rho) \rho^{1-\kappa}  \langle
\partial_{\rho}^{2a-1} Q^{- \frac{2a-1}{2a}} (1-\pi)u_{j},
Q^{\frac{2a-1}{2a}}
\mu_{0}^{\frac{\gamma - 2a + 1}{2a} + \nu} \ {}^{t}\tilde{P}_{1}
\phi_{0} \rangle \|_{A}
\\
\leq
C_{1} \frac{\tilde{C}_{0}}{j^{\frac{1}{2}} A(j)^{1-\kappa}}  \mu_{0}^{\frac{\gamma - 2a + 1}{2a} + \nu} \|
Q^{\frac{2a-1}{2a}} \
{}^{t}\tilde{P}_{1} \phi_{0} \|_{0}
\\
\cdot
\| w_{j}(\rho) \rho^{1-\kappa}
Q^{- \frac{2a-1}{2a}} \partial_{\rho}^{2a-1} (1-\pi) u_{j} \|_{0,A}
\end{multline*}
We may apply the estimate obtained in the preceding section for $
(1-\pi) u_{j} $, since $ - \frac{2a-1}{2a} \geq -1  $ and $ -
\frac{2a-1}{2a} + \frac{2a-1}{2a} = 0 $, $ 2a-1 \leq j + \gamma^{\#}
$. We thus get 
\begin{multline}
  \label{eq:pi-k=0}
C_{1} \| w_{j}(\rho)
    I_{0 i}(\rho^{-1} \langle \tilde{g}_{j 0} ,
    \mu_{0}^{\frac{\gamma_{1} - 2a + 1}{2a}} \phi_{0} \rangle ) \|_{A}
\\
\leq
 \frac{C_{3}}{j^{\frac{1}{2}} A(j)^{1-\kappa}} \mu_{0}^{\frac{\gamma}{2a} + \nu} 
\| w_{j}(\rho) \rho^{1-\kappa}
Q^{- \frac{2a-1}{2a}} \partial_{\rho}^{2a-1} (1-\pi) u_{j} \|_{0,A}
\\
\leq
\frac{C_{3}}{j^{\frac{1}{2}} A(j)^{1-\kappa}} \mu_{0}^{\frac{\gamma}{2a} + \nu}
C_{0}^{1 + \sigma j + \sigma' \frac{(2a-1)^{2}}{2a} } (\lambda(j, 0,
0)+1)^{\lambda(j, 0, 0)},
\end{multline}
where
$$ 
\lambda(j, 0, 0) = \frac{j}{s_{0}} .
$$
Hence from \eqref{eq:pi-k=0} we have
\begin{multline}
  \label{eq:pi-k=0--2}
C_{1} \| w_{j}(\rho)
    I_{0 i}(\rho^{-1} \langle \tilde{g}_{j 0} ,
    \mu_{0}^{\frac{\gamma_{1} - 2a + 1}{2a}} \phi_{0} \rangle ) \|_{A}
\\
\leq
\frac{C_{3} C_{0}^{2 a \sigma'} }{j^{\frac{1}{2}} A(j)^{1-\kappa}}
\mu_{0}^{\frac{\gamma}{2a} + \nu} C_{0}^{-\sigma'(\nu+\gamma)}
C_{0}^{1 + \sigma j + \sigma'(\nu + \gamma) - \sigma'} (\lambda(j, \nu, \gamma)
+ 1)^{\lambda(j, \nu, \gamma)}
\\
\leq
C_{0}^{1 + \sigma j + \sigma'(\nu + \gamma) - \sigma'} (\lambda(j, \nu, \gamma)
+ 1)^{\lambda(j, \nu, \gamma)} ,
\end{multline}
if we make the choice  $ C_{0}^{\sigma'} > \mu_{0} $, $
R_{0}^{1-\kappa} > C_{3} C_{0}^{2 a \sigma'} $.

To conclude the analysis of the last term on the right hand side of
\eqref{eq:puj-2} we must examine \eqref{eq:617} when $ \re
\mu_{0i} = \tilde{\mu}_{0} $ for $ k \in \{1, \ldots, \min\{j, 2a-1\}\}
$.

It is enough to estimate one out of the two terms occurring:
$$
C_{1}  \sum_{k=1}^{\min\{j, 2a -1\}}
  \| w_{j}(\rho)
    I_{0 i}(\rho^{-(k+1)} \langle \tilde{g}_{j k} ,
    \mu_{0}^{\frac{\gamma_{1} - 2a + 1}{2a}} \phi_{0} \rangle ) \|_{A} .
$$
Applying Lemma \ref{lemma:pi-mu0} we have
\begin{multline*}
C_{1}  \sum_{k=1}^{\min\{j, 2a -1\}}
  \| w_{j}(\rho)
    I_{0 i}(\rho^{-(k+1)} \langle \tilde{g}_{j k} ,
    \mu_{0}^{\frac{\gamma_{1} - 2a + 1}{2a}} \phi_{0} \rangle ) \|_{A}
\\
\leq
\tilde{C}_{0}C_{1} \sum_{k=1}^{\min\{j, 2a -1\}}  \frac{1}{j^{\frac{1}{2}} A(j)^{1-\kappa}}  \|
w_{j-k}(\rho) \langle \tilde{g}_{j k} ,
\mu_{0}^{\frac{\gamma_{1} - 2a + 1}{2a}} \phi_{0} \rangle   \|_{A}
\\
=
\tilde{C}_{0}C_{1} \sum_{k=1}^{\min\{j, 2a -1\}}
\frac{1}{j^{\frac{1}{2}} A(j)^{1-\kappa}}
\| w_{j-k}(\rho) \langle Q^{\nu + r} P_{k+1} u_{j-k} ,
\mu_{0}^{\frac{\gamma_{1} - 2a + 1}{2a}} \phi_{0} \rangle   \|_{A}
\\
=
\tilde{C}_{0}C_{1} \sum_{k=1}^{\min\{j, 2a -1\}}
\frac{1}{j^{\frac{1}{2}} A(j)^{1-\kappa}}
\\
\cdot
\| w_{j-k}(\rho) \partial_{\rho}^{2a-k-1} \langle
Q^{-\frac{2a-k-1}{2a}}  u_{j-k} ,
\mu_{0}^{\frac{\gamma - 2a + 1}{2a} + \nu} Q^{\frac{2a-k-1}{2a}}
\ {}^{t}\tilde{P}_{k+1} \phi_{0} \rangle   \|_{A}
\\
\leq
\tilde{C}_{0}C_{1} \sum_{k=1}^{\min\{j, 2a -1\}}
\frac{1}{j^{\frac{1}{2}} A(j)^{1-\kappa}} \mu_{0}^{\frac{\gamma - 2a +
    1}{2a} + \nu} 
\| Q^{\frac{2a-k-1}{2a}} \ {}^{t}\tilde{P}_{k+1} \phi_{0}  \|_{0}
\\
\cdot
\| w_{j-k}(\rho) \partial_{\rho}^{2a-k-1} Q^{-\frac{2a-k-1}{2a}}
u_{j-k}  \|_{0,A}
\\
\leq
\tilde{C}_{0}C_{2}
\frac{\mu_{0}^{\frac{\gamma - 2a +
    1}{2a} + \nu} }{j^{\frac{1}{2}} A(j)^{1-\kappa}} 
 \sum_{k=1}^{\min\{j, 2a -1\}} 
\| w_{j-k}(\rho) \partial_{\rho}^{2a-k-1} Q^{-\frac{2a-k-1}{2a}}
u_{j-k}  \|_{0,A} .
\end{multline*}
Here we absorbed the norm $ \| Q^{\frac{2a-k-1}{2a}} \ {}^{t}\tilde{P}_{k+1}
\phi_{0}  \|_{0} $ into an absolute constant $ C_{2} $ since $ 1\leq k
\leq 2a - 1$.

We point out that the conditions $ \nu \geq -1 $ and $ \nu +
\frac{\gamma}{2a} \geq 0 $, when $ \nu $, $ \gamma $ are the exponents
of $ Q $, $ \partial_{\rho} $ respectively, are satisfied and hence we
may apply the inductive hypothesis thus obtaining the bound
$$
\tilde{C}_{0}C_{2}
\frac{\mu_{0}^{\frac{\gamma - 2a +
    1}{2a} + \nu} }{j^{\frac{1}{2}} A(j)^{1-\kappa}} 
 \sum_{k=1}^{\min\{j, 2a -1\}} 
C_{0}^{1 + \sigma(j-k) + \sigma'(2a-k-1)} (\lambda+1)^{\lambda},
$$
where
$$ 
\lambda = \frac{j-k}{s_{0}} < \lambda(j, \nu, \gamma).
$$
Since $ \sigma = 2a \sigma' $ (see \eqref{eq:sigma}) we may bound the
above quantity by
\begin{multline}
  \label{eq:k>0--mu0i=mu0t}
C_{3} C_{0}^{-\sigma'} \frac{\mu_{0}^{\nu +
    \frac{\gamma}{2a}}}{C_{0}^{\sigma'(\nu + \gamma)}} C_{0}^{1 +
  \sigma j + \sigma'(\nu + \gamma) - \sigma'} (\lambda(j, \nu, \gamma) +
1)^{\lambda(j, \nu, \gamma)}
\\
\leq
C_{0}^{1 +
  \sigma j + \sigma'(\nu + \gamma) - \sigma'} (\lambda(j, \nu, \gamma) +
1)^{\lambda(j, \nu, \gamma)} .
\end{multline}
This concludes the estimation of the second term on the right hand
side of \eqref{eq:puj-2}. We are left with the estimate of the first
term in \eqref{eq:puj-2}.

We split the sum as
\begin{multline}
  \label{eq:kgrt1}
\sum_{s=0}^{r-1} \| w_{j}(\rho)
\Dslash^{\gamma - 2a(1+s)}_{\rho} \left( \frac{1}{\rho}  \pi
  Q^{\nu +s} P_{1} (1-\pi)u_{j} \right ) \|_{0, A}
\\
+
\sum_{k=1}^{\min\{j, 2a-1\}} \sum_{s=0}^{r-1} \| w_{j}(\rho)
\Dslash^{\gamma - 2a(1+s)}_{\rho} \left( \frac{1}{\rho^{k+1} }  \pi
  Q^{\nu +s} P_{k+1} u_{j-k} \right ) \|_{0, A}
\\
= B_{0} + \sum_{k=1}^{\min\{j, 2a-1\}} B_{k}.
\end{multline}
Let us examine $ B_{k} $, $ k > 0 $, first. 
\begin{multline*}
\sum_{s=0}^{r-1} \| w_{j}(\rho)
\Dslash^{\gamma - 2a(1+s)}_{\rho} \left( \frac{1}{\rho^{k+1} }  \pi
  Q^{\nu +s} P_{k+1} u_{j-k} \right ) \|_{0, A}
\\
\leq
C_{4}
\sum_{s=0}^{r-1} \sum_{m=0}^{\gamma - 2a(1+s)} \sum_{\ell=0}^{k+1} \binom{\gamma -
  2a(1+s)}{m} \frac{(k+m)!}{k!} \frac{1}{A(j)^{m}}
\\
\cdot
\| w_{j}(\rho)  \frac{1}{\rho^{k+1} } \partial_{\rho}^{\gamma -
  2as - m - k - 1} 
\pi  Q^{\nu +s} (x \partial_{x})^{\ell} u_{j-k}  \|_{0, A}
\\
\leq
C_{5}
\sum_{s=0}^{r-1} \sum_{m=0}^{\gamma - 2a(1+s)} \sum_{\ell=0}^{k+1}
\binom{k+m}{m}
\frac{(\gamma - 2a(1+s))!}{(\gamma - 2a(1+s) - m)!}
\\
\cdot
\frac{1}{A(j)^{m+1+k(1-\kappa) }} 
\| w_{j-k}(\rho)   \partial_{\rho}^{\gamma -
  2as - m - k - 1} 
\pi  Q^{\nu +s} (x \partial_{x})^{\ell} u_{j-k}  \|_{0, A} .
\end{multline*}
Now $ \pi Q^{\nu +s} (x \partial_{x})^{\ell} u_{j-k} = \langle Q^{\nu +s}
(x \partial_{x})^{\ell} u_{j-k} , \phi_{0} \rangle \phi_{0}   $, so
that we get
$$
\langle (x \partial_{x})^{\ell} u_{j-k} , \mu_{0}^{\nu+s}  \phi_{0} \rangle
\phi_{0} =  \langle u_{j-k} , \mu_{0}^{\nu+s}\ \  {}^{t}(x \partial_{x})^{\ell} \phi_{0} \rangle
\phi_{0} .
$$
Hence the above expression is bounded by
\begin{multline*}
C_{6}
\sum_{s=0}^{r-1}  \mu_{0}^{\nu+s} \sum_{m=0}^{\gamma - 2a(1+s)} \sum_{\ell=0}^{k+1}
\binom{k+m}{m}
\frac{(\gamma - 2a(1+s))!}{(\gamma - 2a(1+s) - m)!}
\\
\cdot
\frac{1}{A(j)^{m+1+k(1-\kappa) }} 
\| w_{j-k}(\rho)   \partial_{\rho}^{\gamma -
  2as - m - k - 1} u_{j-k}  \|_{0, A} ,
\end{multline*}
where we used the fact that since $ \ell $ takes up a finite number of
values, we can bound $ \| {}^{t}(x \partial_{x})^{\ell} \phi_{0}
\|_{0} $ by an absolute constant. Now
\begin{multline}
  \label{eq:factorials}
\binom{k+m}{m} \frac{(\gamma - 2a(1+s))!}{(\gamma -2a(1+s) -m)!}
  \frac{1}{A(j)^{m+1 +k(1-\kappa)}}
\\
=
\binom{k+m}{m} \frac{(\gamma - 2a(1+s))!}{(\gamma -2a(1+s) -m)!}
\frac{1}{(R_{0} (j+1))^{m+1+k(1-\kappa)}}
\\
\leq
\binom{k+m}{m} \frac{(\gamma - 2a(1+s)) \cdots (\gamma - 2a(1+s) - m
  + 1)}{R_{0}^{m+1+k(1-\kappa)} (j+1)^{m+1+k(1-\kappa)}}
\\
\leq
\binom{k+m}{m} \left(\frac{j+\gamma^{\#}}{R_{0}(j+1)}\right)^{m}
\leq 2^{2a-1} \left(\frac{2 \gamma^{\#}}{R_{0}}\right)^{m} .
\end{multline}
Thus the above quantity is estimated by
$$
C_{7} 
\sum_{s=0}^{r-1}  \mu_{0}^{\nu+s} \sum_{m=0}^{\gamma - 2a(1+s)} 
\left(\frac{2 \gamma^{\#}}{R_{0}}\right)^{m}
\| w_{j-k}(\rho)   \partial_{\rho}^{\gamma -
  2as - m - k - 1} u_{j-k}  \|_{0, A} .
$$
Using the inductive hypothesis we have that the above sum is bounded
by
$$
C_{7} 
\sum_{s=0}^{r-1}  \mu_{0}^{\nu+s} \sum_{m=0}^{\gamma - 2a(1+s)} 
\left(\frac{2 \gamma^{\#}}{R_{0}}\right)^{m}
C_{0}^{1+\sigma(j-k) + \sigma'(\gamma - 2as - m - k - 1)} (\lambda+1)^{\lambda},
$$
where
$$ 
\lambda = \frac{j-k}{s_{0}} + \left( \frac{\gamma - 2as - m - k -
    1}{a}\right) \frac{q-1}{q} .
$$
We see immediately that
$$ 
\lambda < \lambda(j, \nu, \gamma).
$$
As for the sums above we remark that, choosing $ R_{0} > 2 R_{1} $, we
may estimate them by
$$
C_{8} 
C_{0}^{-(\sigma + \sigma') k} \left(
  \frac{\mu_{0}}{C_{0}^{\sigma'}}\right)^{\nu} \ 
\sum_{s=0}^{r-1} \left(\frac{ \mu_{0}}{C_{0}^{\sigma' 2 a}}\right)^{s} 
C_{0}^{1+\sigma j + \sigma'(\gamma + \nu)} (\lambda(j,
\nu, \gamma) +1)^{\lambda(j, \nu, \gamma)}
$$
so that the choice $ C_{0}^{\sigma'} > \max\{\mu_{0},
\mu_{0}^{\frac{1}{2a}}\} $, $ C_{0}^{\sigma'} > 2a C_{8} $ allows us to obtain
from \eqref{eq:kgrt1}
\begin{equation}
\label{eq:Bk}
\sum_{k=1}^{\min\{j, 2a-1\}} B_{k} \leq C_{0}^{-\sigma'} 
C_{0}^{1+\sigma j + \sigma'(\gamma + \nu)} (\lambda(j,
\nu, \gamma) +1)^{\lambda(j, \nu, \gamma)} .
\end{equation}
Next we have to estimate in \eqref{eq:k>0--mu0i=mu0t} 
$$ 
B_{0} = \sum_{s=0}^{r-1} \| w_{j}(\rho)
\Dslash^{\gamma - 2a(1+s)}_{\rho} \left( \frac{1}{\rho}  \pi
  Q^{\nu +s} P_{1} (1-\pi)u_{j} \right ) \|_{0, A} .
$$
Taking the $ \rho $-derivative as above and recalling
\eqref{eq:dslash} we have
\begin{multline*}
B_{0} \leq C_{9} \sum_{s=0}^{r-1} \sum_{m=0}^{\gamma - 2a(1+s)}
\binom{\gamma -2a(1+s)}{m} m!
\\
\cdot
\| w_{j}(\rho)
 \frac{1}{\rho^{1+m} } \partial_{\rho}^{\gamma-2a(1+s) - m}   \pi
 Q^{\nu +s} P_{1} (1-\pi)u_{j}  \|_{0, A}
\\
=
C_{9} \sum_{s=0}^{r-1} \sum_{m=0}^{\gamma - 2a(1+s)}
\frac{(\gamma -2a(1+s))!}{(\gamma - 2a(1+s) - m)!} 
\\
\cdot
\| w_{j}(\rho)
 \frac{1}{\rho^{2 - \kappa +m} } \rho^{1-\kappa}  \partial_{\rho}^{\gamma-2a(1+s) - m}   \pi
 Q^{\nu +s} P_{1} (1-\pi)u_{j}  \|_{0, A}
 \\
\leq
C_{10} \sum_{s=0}^{r-1} \mu_{0}^{\nu+s} \sum_{m=0}^{\gamma - 2a(1+s)}
\frac{(\gamma -2a(1+s))!}{(\gamma - 2a(1+s) - m)!} \frac{1}{A(j)^{2-\kappa+m}}
\\
\cdot
\| w_{j}(\rho) \rho^{1-\kappa}  \partial_{\rho}^{\gamma-2a s - m -1} \langle
 (1-\pi)u_{j} , {}^{t}\tilde{P}_{1} \phi_{0} \rangle \phi_{0} \|_{0,
   A}
\\
\leq
C_{11} \sum_{s=0}^{r-1} \mu_{0}^{\nu+s} \sum_{m=0}^{\gamma - 2a(1+s)}
\frac{(\gamma -2a(1+s))!}{(\gamma - 2a(1+s) - m)!} \frac{1}{A(j)^{2-\kappa+m}}
\\
\cdot
\| w_{j}(\rho) \rho^{1-\kappa}  \partial_{\rho}^{\gamma-2a s - m -1} 
 (1-\pi)u_{j} \|_{0,A} ,
\end{multline*}
where we used \eqref{eq:Pj} and included in $ C_{11} $ the norm $ \|
{}^{t}\tilde{P}_{1} \phi_{0} \|_{0} $, which is an absolute
constant .

Arguing as we did in \eqref{eq:factorials} we get
\begin{multline*}
  B_{0} \leq
C_{12} \frac{1}{R_{0}} 
\sum_{s=0}^{r-1}  \mu_{0}^{\nu+s} \sum_{m=0}^{\gamma - 2a(1+s)} 
\left(\frac{2 \gamma^{\#}}{R_{0}}\right)^{m}
\\
\cdot
\| w_{j}(\rho)  \rho^{1-\kappa} \partial_{\rho}^{\gamma -
  2as - m - 1} (1-\pi) u_{j}  \|_{0, A} . 
\end{multline*}
Using the inductive hypothesis we may bound the above quantity by
$$
C_{12} \frac{1}{R_{0}} 
\sum_{s=0}^{r-1}  \mu_{0}^{\nu+s} \sum_{m=0}^{\gamma - 2a(1+s)} 
\left( \frac{2 \gamma^{\#}}{R_{0}}\right)^{m}
C_{0}^{1+\sigma j + \sigma'(\gamma - 2as - m - 1)} (\lambda+1)^{\lambda},
$$
where
$$ 
\lambda = \frac{j}{s_{0}} + \left( \frac{\gamma - 2as - m -
    1}{a}\right) \frac{q-1}{q} .
$$
We see immediately that
$$ 
\lambda < \lambda(j, \nu, \gamma).
$$
Taking $ R_{0} $ large enough and choosing $ C_{0} $ as above we
conclude that
\begin{equation}
\label{eq:B0}
B_{0} \leq C_{0}^{1+\sigma j + \sigma'(\gamma + \nu) - \sigma'}
(\lambda(j, \nu, \gamma) + 1)^{\lambda(j, \nu, \gamma)}. 
\end{equation}
Summing up, when we plug \eqref{eq:B0}, \eqref{eq:Bk},
\eqref{eq:k>0--mu0i=mu0t}, \eqref{eq:pi-k=0--2}, \eqref{eq:kgeq0} into
\eqref{eq:puj-2} and choose $ C_{0} $ sufficiently large we get the
desired estimate.

This ends the proof of Theorem \ref{th:est-p}.
\end{proof}

\section{End of the Proof of Theorem \ref{th:est}}
\setcounter{equation}{0}
\setcounter{theorem}{0}
\setcounter{proposition}{0}
\setcounter{lemma}{0}
\setcounter{corollary}{0}
\setcounter{definition}{0}

Using Theorem \ref{th:est-1-p}, \ref{th:est-p} and renaming $
C_{0} $ we obtain Theorem \ref{th:est'}.

Next we want to prove Theorem \ref{th:est}.

To this end we are going to use Proposition \ref{prop:ineqQ} with $
\theta = 0 $:
\begin{multline*}
\| w_{j}(\rho) \partial_{\rho}^{\gamma} x^{\beta}
\partial_{x}^{\alpha} u_{j}(x, \rho)  \|_{0, A}
\\
\leq
 \sum_{\ell = 0}^{m(0, I)}
C^{\ell + 1} \binom{m(0, I)}{\ell} \langle I \rangle^{\ell} \|
w_{j}(\rho) \partial_{\rho}^{\gamma} 
Q^{\frac{1}{2}\left(\langle I \rangle - \ell \frac{q}{q-1} \right)} u_{j}
\|_{0, A},
\end{multline*}
where
$$ 
\langle I \rangle = \alpha + \frac{\beta}{q-1},
$$
and
$$ 
m(0, I) = \left[ \langle I \rangle \frac{q-1}{q} \right] ,
$$
where the square brackets denote the integer part.

By Theorem \ref{th:est'} we have
\begin{multline*}
\| w_{j}(\rho) \partial_{\rho}^{\gamma} x^{\beta}
\partial_{x}^{\alpha} u_{j}(x, \rho)  \|_{0, A}
\\
\leq
 \sum_{\ell = 0}^{m(0, I)}
C^{\ell + 1} \binom{m(0, I)}{\ell} \langle I \rangle^{\ell}
C_{0}^{1+\sigma j + \sigma' (\gamma + \langle I \rangle - \ell
  \frac{q}{q-1} )}
\\
\cdot
\left(1 + \lambda(j, \frac{1}{2}\left(\langle I \rangle
  - \ell \frac{q}{q-1} \right), \gamma)\right)^{\lambda(j, \frac{1}{2}\left(\langle I \rangle
    - \ell \frac{q}{q-1} \right), \gamma)} .
\end{multline*}
Now
$$ 
\lambda(j, \frac{1}{2}\left(\langle I \rangle
  - \ell \frac{q}{q-1} \right), \gamma) = \frac{j}{s_{0}} + \langle I
\rangle \frac{q-1}{q} + \frac{\gamma}{a} \frac{q-1}{q} - \ell,
$$
and, since
$$
\langle I \rangle \leq \frac{q}{q-1} \lambda(j, \frac{1}{2}\left(\langle I \rangle
  - \ell \frac{q}{q-1} \right), \gamma) ,
$$
we obtain that
\begin{multline*}
\| w_{j}(\rho) \partial_{\rho}^{\gamma} x^{\beta}
\partial_{x}^{\alpha} u_{j}(x, \rho)  \|_{0, A}
\\
\leq
 \sum_{\ell = 0}^{m(0, I)}
C_{1}^{\ell + 1} \binom{m(0, I)}{\ell}
C_{0}^{1+\sigma j + \sigma' (\gamma + \langle I \rangle - \ell
  \frac{q}{q-1} )}
\left(1 + \lambda(j, \frac{1}{2} \langle I \rangle, \gamma)\right)^{\lambda(j, \frac{1}{2}\left(\langle I \rangle
  \right), \gamma)}
\\
\leq
C_{0}^{1+\sigma j + \sigma' (\gamma + \langle I \rangle )}
\left(1 + \lambda(j, \frac{1}{2} \langle I \rangle, \gamma)\right)^{\lambda(j, \frac{1}{2}\left(\langle I \rangle
  \right), \gamma)}
\\
\sum_{\ell = 0}^{m(0, I)}
C_{1}^{\ell + 1} \binom{m(0, I)}{\ell} C_{0}^{- \ell \frac{q}{q-1}\sigma'}
\\
\leq
C_{0}^{1+\sigma j + \sigma' (\gamma + \langle I \rangle )}
\left(1 + \lambda(j, \frac{1}{2} \langle I \rangle, \gamma)\right)^{\lambda(j, \frac{1}{2}\left(\langle I \rangle
  \right), \gamma)}
C_{1} \left( 1 + \frac{C_{1}}{C_{0}^{\frac{q}{q-1}\sigma'}} \right)^{m(0, I)} .
\end{multline*}
Choosing $ C_{0}^{\frac{q}{q-1}\sigma'} > C_{1} $ and recalling that $ m(0, I) <
\langle I \rangle \leq \alpha + \beta $ we have that
\begin{multline*}
\| w_{j}(\rho) \partial_{\rho}^{\gamma} x^{\beta}
\partial_{x}^{\alpha} u_{j}(x, \rho)  \|_{0, A}
\\
\leq
2^{\alpha+\beta} C_{1} C_{0}^{1+\sigma j + \sigma' (\gamma + \langle I \rangle )}
\left(1 + \lambda(j, \frac{1}{2} \langle I \rangle, \gamma)\right)^{\lambda(j, \frac{1}{2}\left(\langle I \rangle
  \right), \gamma)} .
\end{multline*}
Moreover
$$ 
\lambda(j, \frac{1}{2} \langle I \rangle, \gamma) = \frac{j}{s_{0}} +
\alpha \frac{q-1}{q} + \frac{\beta}{q} + \frac{\gamma}{a}
\frac{q-1}{q} ,
$$
so that
\begin{multline}
\label{eq:4.1}
\| w_{j}(\rho) \partial_{\rho}^{\gamma} x^{\beta}
\partial_{x}^{\alpha} u_{j}(x, \rho)  \|_{0, A}
\leq
C_{u}^{1 + j + \alpha + \beta + \gamma}
\\
\cdot
\left( 1 + \frac{j}{s_{0}} +
\alpha \frac{q-1}{q} + \frac{\beta}{q} + \frac{\gamma}{a}
\frac{q-1}{q} \right)^{\frac{j}{s_{0}} +
\alpha \frac{q-1}{q} + \frac{\beta}{q} + \frac{\gamma}{a}
\frac{q-1}{q}} .
\end{multline}
Renaming the constant $ C_{u} $, we finish the proof of Theorem \ref{th:est}.

\section{Pointwise Estimates of the $ u_{j} $}
\setcounter{equation}{0}
\setcounter{theorem}{0}
\setcounter{proposition}{0}
\setcounter{lemma}{0}
\setcounter{corollary}{0}
\setcounter{definition}{0}
\label{sec:pt}

First of all we point out that from the estimates of Theorem
\ref{th:est} it is straightforward to deduce the same type of
estimates for
\begin{multline}
\label{eq:norm}
\| \partial_{\rho}^{\gamma'}  \partial_{x}^{\alpha'} w_{j}(\rho) \partial_{\rho}^{\gamma} x^{\beta}
\partial_{x}^{\alpha} u_{j}(x, \rho)  \|_{0, A}
\leq
C_{u}^{1 + j + \alpha + \beta + \gamma}
\\
\cdot
\left(\frac{j}{s_{0}} +
\alpha \frac{q-1}{q} + \frac{\beta}{q} + \frac{\gamma}{a}
\frac{q-1}{q} \right)!
\end{multline}
with $ \alpha' + \gamma' \leq 2 $, with a different meaning of the
constant $ C_{u} $.

The purpose of adding two extra derivatives with
respect to $ (x, \rho) $ is to apply the Sobolev immersion theorem to
deduce pointwise estimates.

To this end we use cutoff functions in order to apply the Sobolev
immersion theorem on the whole plane.
\begin{lemma}
\label{lemma:chij}
There exist smooth functions $ \chi_{j}(\rho) $ such that
\begin{itemize}
\item[(i)  ]{}
$ \supp \chi_{j} \subset \{ \rho \geq R_{0}(j+1) \} $.
\item[(ii) ]{}
$ \chi_{j} \equiv 1 $ for $ \rho \geq 3 R_{0}(j+1) $.
\item[(iii)]{}
$ | D^{\gamma}_{\rho} \chi_{j}(\rho) | \leq C_{\chi}^{\gamma} $, for $
0 \leq \gamma \leq R_{0}(j+1) $.
\end{itemize}
\end{lemma}
\begin{proof}
Let $ \chi $ denote the characteristic function of the half line $ [2
R_{0}(j+1), + \infty [ $, and denote by $ \psi $ a function in $
C_{0}^{\infty}(\R) $, $ \supp \psi \subset \{ |x| \leq 1 \} $, and
such that $ \int_{\R}\psi(x) dx = 1 $. 

Define, assuming that $ R_{0} $ is an integer,
$$ 
\chi_{j} = \chi * \underbrace{\psi * \cdots * \psi}_{R_{0}(j+1) \text{
  times}}. 
$$
The support of $ \chi_{j} $ is evidently contained in $ \{ \rho \geq
R_{0}(j+1) \} $ and $ \chi_{j} \equiv 1 $ for $ \rho \geq 3 R_{0}(j+1)
$. 

Moreover
$$ 
D^{\gamma} \chi_{j}(\rho) = \chi * \underbrace{D\psi * \cdots * D\psi}_{\gamma \text{
  times}} * \psi * \cdots * \psi ,
$$
for $ \gamma \leq R_{0}(j+1) $. Hence, by Young inequality for the convolution,
$$ 
| D^{\gamma} \chi_{j}(\rho) | \leq \left( \| D\psi \|_{L^{1}(\R)}
  \right)^{\gamma} .
$$
This completes the proof of the lemma.
\end{proof}
Consider now
$$ 
\| \partial_{\rho}^{\gamma'}  \partial_{x}^{\alpha'} w_{j}(\rho) \partial_{\rho}^{\gamma} x^{\beta}
\partial_{x}^{\alpha} \chi_{j}(\rho)  u_{j}(x, \rho)  \|_{0} ,
$$
where $ \| \cdot \|_{0} $ denotes the norm in $ L^{2}(\R_{\rho} \times
\R_{x})$. We have
\begin{multline*}
\| \partial_{\rho}^{\gamma'}  \partial_{x}^{\alpha'} w_{j}(\rho) \partial_{\rho}^{\gamma} x^{\beta}
\partial_{x}^{\alpha} \chi_{j}(\rho)  u_{j}(x, \rho)  \|_{0}
\\
\leq
\| \chi_{j}(\rho) \partial_{\rho}^{\gamma'}  \partial_{x}^{\alpha'}
w_{j}(\rho) \partial_{\rho}^{\gamma} x^{\beta} \partial_{x}^{\alpha}
u_{j}(x, \rho)  \|_{0}
\\
+
2^{\gamma+\gamma'} \sum_{\substack{k=1\\ \gamma_{1}' + \gamma_{1} = \gamma' + \gamma -k}}^{\gamma+\gamma'} \binom{\gamma + \gamma'}{k}
\\
\cdot
\| D^{k}_{\rho}\chi_{j}(\rho) \cdot \partial_{\rho}^{\gamma'_{1}} \left( \partial_{x}^{\alpha'} w_{j}(\rho) \partial_{\rho}^{\gamma_{1}} x^{\beta}
\partial_{x}^{\alpha}   u_{j}(x, \rho) \right) \|_{0}
\\
\leq
\| \partial_{\rho}^{\gamma'}  \partial_{x}^{\alpha'}
w_{j}(\rho) \partial_{\rho}^{\gamma} x^{\beta} \partial_{x}^{\alpha}
u_{j}(x, \rho)  \|_{0, A}
\\
+
2^{\gamma+\gamma'} \sum_{\substack{k=1\\ \gamma_{1}' + \gamma_{1} = \gamma' + \gamma
    -k}}^{\gamma+\gamma'} \binom{\gamma + \gamma'}{k}
C_{\chi}^{k}
\| \partial_{\rho}^{\gamma'_{1}}  \partial_{x}^{\alpha'} w_{j}(\rho) \partial_{\rho}^{\gamma_{1}} x^{\beta}
\partial_{x}^{\alpha}   u_{j}(x, \rho)  \|_{0, A}
\\
\leq
C_{u}^{1 + j + \alpha + \beta + \gamma}
\left(\frac{j}{s_{0}} +
\alpha \frac{q-1}{q} + \frac{\beta}{q} + \frac{\gamma}{a}
\frac{q-1}{q} \right)!
\\
+ 2^{\gamma+\gamma'}
\sum_{\substack{k=1\\ \gamma_{1}' + \gamma_{1} = \gamma' + \gamma
    -k}}^{\gamma+\gamma'} \binom{\gamma + \gamma'}{k}
C_{\chi}^{k}
C_{u}^{1 + j + \alpha + \beta + \gamma_{1}}
\\
\cdot
\left(\frac{j}{s_{0}} +
\alpha \frac{q-1}{q} + \frac{\beta}{q} + \frac{\gamma_{1}}{a}
\frac{q-1}{q} \right)!
\\
\leq
C_{1}^{j+\alpha+\beta+\gamma' + \gamma} 
C_{u}^{1 + j + \alpha + \beta + \gamma + \gamma'}
\left(\frac{j}{s_{0}} +
\alpha \frac{q-1}{q} + \frac{\beta}{q} + \frac{\gamma}{a}
\frac{q-1}{q} \right)!
\\
\leq
((C_{1} C_{u})^{2})^{1 + j + \alpha + \beta + \gamma} \left(\frac{j}{s_{0}} +
\alpha \frac{q-1}{q} + \frac{\beta}{q} + \frac{\gamma}{a}
\frac{q-1}{q} \right)!
\end{multline*}
provided $ C_{u} > C_{\chi} $ and $ \gamma \leq j+\gamma^{\#}-2 $. Summing
up we obtain
\begin{multline}
\label{eq:estR}
\| \partial_{\rho}^{\gamma'}  \partial_{x}^{\alpha'} w_{j}(\rho) \partial_{\rho}^{\gamma} x^{\beta}
\partial_{x}^{\alpha} \chi_{j}(\rho)  u_{j}(x, \rho)  \|_{0}
\\
\leq
C_{1u}^{1 + j + \alpha + \beta + \gamma} \left(\frac{j}{s_{0}} +
\alpha \frac{q-1}{q} + \frac{\beta}{q} + \frac{\gamma}{a}
\frac{q-1}{q} \right)!
\end{multline}
for a suitable constant $ C_{1u} > 0 $.

The Sobolev embedding theorem yields
\begin{multline*}
 | w_{j}(\rho) \partial_{\rho}^{\gamma} x^{\beta}
\partial_{x}^{\alpha}  u_{j}(x, \rho) |
\leq
\| w_{j}(\rho) \partial_{\rho}^{\gamma} x^{\beta}
\partial_{x}^{\alpha} \chi_{j}(\rho) u_{j}(x, \rho)
\|_{L^{\infty}(\R^{2})} 
\\
\leq
C_{2u}^{1 + j + \alpha + \beta + \gamma} \left(\frac{j}{s_{0}} +
\alpha \frac{q-1}{q} + \frac{\beta}{q} + \frac{\gamma}{a}
\frac{q-1}{q} \right)!
\end{multline*}
for $ \rho \geq 3 R_{0}(j+1) $. As a consequence we get, recalling
\eqref{eq:wj}, \eqref{eq:mitlda0}, 
\begin{multline}
\label{eq:pest}
| \partial_{\rho}^{\gamma} x^{\beta}
\partial_{x}^{\alpha}  u_{j}(x, \rho) |
\\
\leq
K_{u}^{1 + j + \alpha + \beta + \gamma}  e^{\tilde{\mu}_{0} \rho}
\rho^{-(j-\delta) \kappa} \left(\frac{j}{s_{0}} +
\alpha \frac{q-1}{q} + \frac{\beta}{q} + \frac{\gamma}{a}
\frac{q-1}{q} \right)!
\end{multline}
So far we have proved the
\begin{theorem}
\label{th:pest-t}
Let the functions $ u_{j} $ denote the solutions of the equations \eqref{eq:tr1} and
\eqref{eq:tr2}, for $ j \geq 1 $ and let $ u_{0} $ denote the function
in \eqref{eq:u0-}, \eqref{eq:mitlda0}.  Then for every $ R_{1} > 0 $
there exist positive constants $ K_{u} $, $ R_{0} > R_{1} $, $ \delta
$, $ \kappa $, with $ 0 < \delta < 1 $, $ \frac{1}{s_{0}} < \kappa < 1
$, $ \kappa \delta > \frac{1}{2} $, such that for $ j \geq 1 $ the
estimates \eqref{eq:pest} are satisfied for $ \rho \geq 3 R_{0} (j+1)
$ and $ \gamma \leq j+ \gamma^{\#} $, where $ \gamma^{\#} $ denotes a
fixed positive constant $ \geq 2a $.  
\end{theorem}
Next we also need a lemma to estimate $ x $-derivatives of the
projections $ \pi $ and $ 1-\pi $ applied to $ \partial_{\rho}^{\gamma} x^{\beta}
\partial_{x}^{\alpha}  u_{j}(x, \rho) $. This will be used in the
proof of Theorem \ref{th:Beu}. 

\begin{lemma}
\label{lemma:derproj}
Under the same assumptions of Theorem \ref{th:pest-t} we have the
estimates
\begin{multline}
\label{eq:derproj1}
\left| \partial_{x}^{\alpha'} \pi \left(\partial_{\rho}^{\gamma} x^{\beta}
    \partial_{x}^{\alpha}  u_{j}(x, \rho) \right) \right|
\leq
M_{u}^{1 + j + \alpha + \alpha'+ \beta + \gamma}
\\
\cdot
e^{\tilde{\mu}_{0} \rho}
\rho^{-(j-\delta) \kappa}
\left(\frac{j}{s_{0}} +
(\alpha+\alpha') \frac{q-1}{q} + \frac{\beta}{q} + \frac{\gamma}{a}
\frac{q-1}{q} \right)!
\end{multline}
and
\begin{multline}
\label{eq:derproj2}
\left| \partial_{x}^{\alpha'} (1-\pi) \left(\partial_{\rho}^{\gamma} x^{\beta}
    \partial_{x}^{\alpha}  u_{j}(x, \rho) \right) \right|
\leq
M_{u}^{1 + j + \alpha + \alpha'+ \beta + \gamma}
\\
\cdot
e^{\tilde{\mu}_{0} \rho}
\rho^{-(j-\delta) \kappa}
\left(\frac{j}{s_{0}} +
(\alpha+\alpha') \frac{q-1}{q} + \frac{\beta}{q} + \frac{\gamma}{a}
\frac{q-1}{q} \right)!
\end{multline}
for $ \rho \geq 3 R_{0} (j+1) $ and $ \gamma \leq j+\gamma^{\#} $ and for
any $ \alpha $, $ \alpha' $.  
\end{lemma}
\begin{proof}
Let us start by proving \eqref{eq:derproj1}. Since
\begin{multline*}
\left| \partial_{x}^{\alpha'} \pi \left(\partial_{\rho}^{\gamma} x^{\beta}
    \partial_{x}^{\alpha}  u_{j}(x, \rho) \right) \right|
=
\left| \langle \partial_{\rho}^{\gamma} x^{\beta}
    \partial_{x}^{\alpha}  u_{j}(x, \rho) , \phi_{0} \rangle
    \partial_{x}^{\alpha'} \phi_{0}(x) \right|
\\
\leq
\| \partial_{\rho}^{\gamma} x^{\beta} \partial_{x}^{\alpha}  u_{j}(x,
\rho) \|_{L^{\infty}(\R_{x})}  \| \phi_{0} \|_{L^{1}(\R_{x})} \left|
  \partial_{x}^{\alpha'} \phi_{0}(x) \right| .
\end{multline*}
By Corollary \ref{cor:Dphi}, Appendix B, we have that
$$ 
| \partial_{x}^{\alpha'} \phi_{0}(x) | \leq C_{\phi}^{\alpha' + 1} \alpha'!^{\frac{q-1}{q}},
$$
and moreover $  \| \phi_{0} \|_{L^{1}(\R_{x})} \leq \tilde{C}_{\phi} $
for a suitable constant $ \tilde{C}_{\phi} > 0 $. Thus if $ \rho \geq
3 R_{0} (j+1) $ from Theorem \ref{th:pest-t} we conclude that
\begin{multline*}
\left| \partial_{x}^{\alpha'} \pi \left(\partial_{\rho}^{\gamma} x^{\beta}
    \partial_{x}^{\alpha}  u_{j}(x, \rho) \right) \right| \leq
\tilde{C}_{\phi} \cdot C_{\phi}^{\alpha' + 1} \alpha'!^{\frac{q-1}{q}}
\\
\cdot
K_{u}^{1 + j + \alpha + \beta + \gamma}  e^{\tilde{\mu}_{0} \rho}
\rho^{-(j-\delta) \kappa} \left(\frac{j}{s_{0}} +
\alpha \frac{q-1}{q} + \frac{\beta}{q} + \frac{\gamma}{a}
\frac{q-1}{q} \right)!
\end{multline*}
which gives \eqref{eq:derproj1} provided $ M_{u} $ is chosen
sufficiently large.

Finally to prove \eqref{eq:derproj2} it is enough to bound each term
of the difference, noting that the bound of the term with the identity
follows from \eqref{eq:pest}, while the bound of the term including the
projection has been already proved.

This ends the proof of the lemma.
\end{proof}

\section{Turning a Formal Solution into a True Solution}
\setcounter{equation}{0}
\setcounter{theorem}{0}
\setcounter{proposition}{0}
\setcounter{lemma}{0}
\setcounter{corollary}{0}
\setcounter{definition}{0}

In the preceding sections we constructed functions $ u_{j} $, $ j \geq
0$, such that the formal sum $ u = \sum_{j\geq 0} u_{j} $ formally
verifies the equation \eqref{eq:transp}
\begin{multline}
\label{eq:Pu-form}
P(x, y, D_{x}, D_{y}) A(u)(x, y)
\\
=
\int_{0}^{+\infty} e^{i y \rho^{s_{0}}} \rho^{r + 2 \frac{s_{0}}{q}}
\left[ \sum_{j=0}^{2a} \frac{1}{\rho^{j}} P_{j}(x, D_{x}, D_{\rho})
  u(x, \rho) \right]_{x = x\rho^{\frac{s_{0}}{q}}} d\rho = 0 .
\end{multline}
Our purpose is to turn $ u $ into a true solution of an equation of
the form
\begin{equation}
\label{eq:PAu=f}
P(x, y, D_{x}, D_{y}) A(u)(x, y) = f(x, y),
\end{equation}
for some smooth function $ f $ defined in $ \R^{2}_{(x, y)} $.

In order to do so we need a set of cutoff functions similar to those
discussed in Lemma \ref{lemma:chij}.
\begin{lemma}
\label{lemma:psij}
There exist smooth functions $ \psi_{j}(\rho) $, $ j \geq 0 $, such that
\begin{itemize}
\item[(i)  ]{}
$ \supp \psi_{j} \subset \{ \rho \geq 3 R_{0}(j+1) \} $.
\item[(ii) ]{}
$ \psi_{j} \equiv 1 $ for $ \rho \geq 6 R_{0}(j+1) $.
\item[(iii) ]{}
$ | D^{\gamma}_{\rho} \psi_{j}(\rho) | \leq C_{\psi}^{\gamma} $, for $
0 \leq \gamma \leq R_{0}(j+1) $.
\end{itemize}
\end{lemma}
The proof is the same as that of Lemma \ref{lemma:chij} and we omit
it.

Define
\begin{equation}
\label{eq:v}
v(x, \rho) = \sum_{j\geq 0} \psi_{j}(\rho) u_{j}(x, \rho) .
\end{equation}
The function $ v $ is well defined and smooth since the above sum is
locally finite in $ \rho $. 

The remaining part of the present section is devoted to computing $ P
A(v) $, where $ v $ is defined by \eqref{eq:v}. Obviously $ A(v) $ is
no more a null solution of $ P $, due to the errors introduced by the
cutoff functions $ \psi_{j} $.

First we recall the definition of the Beurling classes:
\begin{definition}
\label{def:Bs}
Let $ \Omega $ be an open subset of $ \R^{2} $ and $ s \geq 1 $. We
define the class $ 
\mathscr{B}^{s}(\Omega) $ (of Beurling type functions on $ \Omega $)
as the set of all smooth functions $ f(x, y) $ defined on $ \Omega $
such that for every $ \epsilon > 0 $ and for every $ K \Subset \Omega
$ compact, there exists a positive constant $ C = C(\epsilon, K) $
such that
\begin{equation}
\label{eq:Bs}
| \partial_{x}^{\alpha} \partial_{y}^{\beta} f(x, y)| \leq C
\epsilon^{\alpha+\beta} \alpha!^{s} \beta!^{s},
\end{equation}
for every $ (x, y) \in K $ and every $ \alpha $, $ \beta $.

If $ \Omega \subset \R^{2} $ we define the global class $
\mathscr{B}^{s}_{g}(\Omega) $ as the set of all smooth functions $ f(x, y) $
defined on $ \Omega $ such that for every $ \epsilon > 0 $ there
exists a positive constant $ C = C(\epsilon, \Omega) $ such that
\eqref{eq:Bs} holds for every $ (x, y) \in \Omega $ and every $ \alpha
$, $ \beta $. 
\end{definition}
\begin{remark}
\label{rem:1}
It is a consequence of the above definition that a function $ f \in
G^{s'}(\Omega)$ belongs to $ \mathscr{B}^{s}(\Omega) $, for $ s' < s $.
\end{remark}
We point out that $ v(x, \rho) = 0 $ if $ \rho \leq 3 R_{0} $ and
that, by \eqref{eq:pest}, the support properties of
the $ \psi_{j} $, using an analogous but simpler argument as in
the proof of Theorem \ref{th:Beu} below,
\begin{equation}
  \label{eq:convv}
| \rho^{k} \partial_{x}^{\alpha} v(x, \rho) | \leq C_{\alpha k}
e^{-\tilde{\lambda} \rho} ,
\end{equation}
for suitable positive constants $ C_{\alpha k} $, $ \tilde{\lambda} <
|\tilde{\mu}_{0}| $. As a consequence $ A(v) $ is a smooth function
defined on $ \R^{2} $.

The idea behind the next theorem is the following: applying $ P $ to
the function $ A(v) $ yields a function of the form $ A(\sum_{j} \tilde{v}_{j}) $,
where, because of the transport equations and the support properties
of the cutoff functions $ \psi_{j} $, $ \tilde{v}_{j} =
\mathscr{O}(\rho^{-j}) $, has compact $\rho $-support, where $  \rho \sim j $. This
implies that the sum is $ \mathscr{O}(e^{- c \rho \log \rho}) $. Via
Lemma \ref{lemma:intBs0} we see that $ P A(v) $ belongs to $
\mathscr{B}^{s_{0}}_{g}(\R^{2}) $. 
\begin{theorem}
\label{th:Beu}
We have that
\begin{equation}
\label{eq:Beu}
P(x, y, D_{x}, D_{y}) A(v)(x, y) \in \mathscr{B}^{s_{0}}_{g}(\R^{2}).
\end{equation}
\end{theorem}
\begin{proof}
By \eqref{eq:Pu-form} we have
\begin{multline*}
P(x, y, D_{x}, D_{y}) A(v)(x, y)
\\
=
\int_{0}^{+\infty} e^{i y \rho^{s_{0}}} \rho^{r + 2 \frac{s_{0}}{q}}
\left[ \sum_{k=0}^{2a} \frac{1}{\rho^{k}} P_{k}(x, D_{x}, D_{\rho})
  v(x, \rho) \right]_{x = x\rho^{\frac{s_{0}}{q}}} d\rho .
\end{multline*}
Since the sum $ v = \sum_{j\geq 0} \psi_{j}(\rho) u_{j}(x, \rho) $ is
a finite sum for any $ \rho $ in a compact set, we have
\begin{multline*}
\sum_{k=0}^{2a} \frac{1}{\rho^{k}} P_{k}(x, D_{x}, D_{\rho})
  v(x, \rho) = \sum_{k=0}^{2a} \frac{1}{\rho^{k}} P_{k}(x, D_{x},
  D_{\rho}) \sum_{j\geq 0} \psi_{j} u_{j}
\\
=  
\sum_{j\geq 0} \sum_{k=0}^{2a} \frac{1}{\rho^{k}} \sum_{\gamma=0}^{2a} D_{\rho}^{\gamma}
 \psi_{j} \frac{1}{\gamma!} P_{k}^{(\gamma)}(x, D_{x}, D_{\rho}) u_{j} ,
\end{multline*}
where $ P_{k}^{(\gamma)} $ denotes the symbol $
\partial_{\sigma}^{\gamma} P_{k}(x, \xi, \sigma) $. Hence
\begin{multline}
  \label{eq:Si}
\sum_{k=0}^{2a} \frac{1}{\rho^{k}} P_{k}(x, D_{x}, D_{\rho})
  v(x, \rho)
=
\sum_{j\geq 0} \sum_{k=0}^{2a} \frac{1}{\rho^{k}} \psi_{j}
  P_{k} u_{j}
\\
+
\sum_{j\geq 0} \sum_{k=0}^{2a} \frac{1}{\rho^{k}} \sum_{\gamma=1}^{2a} D_{\rho}^{\gamma}
\psi_{j} \frac{1}{\gamma!} P_{k}^{(\gamma)}(x, D_{x}, D_{\rho}) u_{j}
\\
= S_{1}(x, \rho) + S_{2}(x, \rho) .
\end{multline}
We remark that the $ \rho $-derivatives appearing in the above
expression have order at most $ 2a $.

For the first summand we are going to organize the terms using the
transport equations \eqref{eq:tr1} and \eqref{eq:tr2}, while for the
second summand we use the estimates \eqref{eq:pest} and the fact that
the derivatives of $ \psi_{j} $ have compact support.

Define
$$ 
f_{i}(x, y) = \int_{0}^{+\infty} e^{i y \rho^{s_{0}}} \rho^{r + 2 \frac{s_{0}}{q}}
\left[ S_{i}(x, \rho) \right]_{x = x\rho^{\frac{s_{0}}{q}}} d\rho .
$$
Let us start with $ f_{2} $. We have
$$
\partial_{x}^{\alpha} \partial_{y}^{\beta} f_{2}(x, y)
=
i^{\beta}
\int_{0}^{+\infty} e^{i y \rho^{s_{0}}} \rho^{r + (2+ \alpha)
  \frac{s_{0}}{q} +\beta s_{0} }
\left[ \partial_{x}^{\alpha} S_{2}(x, \rho) \right]_{x =
  x\rho^{\frac{s_{0}}{q}}} d\rho .
$$
In view of \eqref{eq:Pj}, the $ x $-derivatives of $ S_{2} $ are given by
$$
\partial_{x}^{\alpha} S_{2}(x, \rho)
=
\sum_{j\geq 0} \sum_{k=0}^{2a}  \sum_{\gamma=1}^{2a}
\sum_{\alpha'=0}^{k} \binom{\alpha}{\alpha'} 
\frac{1}{\rho^{k}}D_{\rho}^{\gamma}
\psi_{j} \frac{1}{\gamma!} P_{k, (\alpha') }^{(\gamma)}(x, D_{x},
D_{\rho}) \partial_{x}^{\alpha-\alpha'} u_{j} ,
$$
where $ P_{k, (\alpha') }^{(\gamma)} $ denotes the symbol $
\partial_{x}^{\alpha'} \partial_{\sigma}^{\gamma} P_{k}(x, \xi,
\sigma) $ and we can assume without loss of generality that $ \alpha >
2a$. 
Let us estimate $ | \partial_{x}^{\alpha} S_{2}(x, \rho) |  $.
We have
\begin{multline}
  \label{eq:daS2}
| \partial_{x}^{\alpha} S_{2}(x, \rho) |
\\  
\leq
\sum_{j\geq 0} \sum_{k=0}^{2a}  \sum_{\gamma=1}^{2a}
\sum_{\alpha'=0}^{k} \binom{\alpha}{\alpha'} 
| D_{\rho}^{\gamma} \psi_{j}| \frac{1}{\gamma!} \left | P_{k, (\alpha') }^{(\gamma)}(x, D_{x},
D_{\rho}) \partial_{x}^{\alpha-\alpha'} u_{j} \right | .
\end{multline}
Recalling the definitions \eqref{eq:P0}, \eqref{eq:P1}, \eqref{eq:Pj}
and \eqref{eq:Pjtilda}, we see that the above sum has a typical
summand of the form
$$ 
\binom{\alpha}{\alpha'} \frac{1}{\gamma!} | D_{\rho}^{\gamma}
\psi_{j}| \left| \partial_{\rho}^{2a -k - \gamma} x^{m-\alpha'}
  \partial_{x}^{m+\alpha-\alpha'} u_{j} \right|, 
$$
if $ k > 0 $, for $ m = 0, 1, \ldots, k $, with the proviso that the
term is zero if $ 2a - k - \gamma < 0 $ or $ m-\alpha' < 0 $, and of the form
$$ 
\binom{\alpha}{\alpha'} \frac{1}{\gamma!} | D_{\rho}^{\gamma}
\psi_{j}| \left|\partial_{\rho}^{2a - \gamma} 
  \partial_{x}^{\alpha} u_{j} \right| ,
$$
if $ k=0 $. Let us start by examining the first type. First observe
that, since $ \gamma > 0 $, the $ \rho $-support of every term is
contained in the interval $ 3 R_{0}(j+1) \leq \rho \leq 6 R_{0}(j+1) $
and that $ |D_{\rho}^{\gamma} \psi_{j}| $ is uniformly bounded since $
\gamma $ runs over a finite set of indices. As for the function $
u_{j} $ we may apply \eqref{eq:pest}, thus obtaining
\begin{multline*}
\left| \partial_{\rho}^{2a -k - \gamma} x^{m-\alpha'}
  \partial_{x}^{m+\alpha-\alpha'} u_{j} \right| 
\\
\leq
K_{u}^{1+j+m+\alpha-\alpha'+m-\alpha'+2a-k-\gamma} e^{\tilde{\mu}_{0}
  \rho} \rho^{-(j-\delta)\kappa}
\\
\cdot \left(\frac{j}{s_{0}} +
  (m+\alpha-\alpha')\frac{q-1}{q} + \frac{m-\alpha'}{q} +
  \frac{2a-k-\gamma}{a} \frac{q-1}{q}\right)!
\\
\leq
K_{u}^{1+j+k+\alpha-2\alpha'+2a-\gamma} e^{\tilde{\mu}_{0}
  \rho} \rho^{-(j-\delta)\kappa}
\\
\cdot
\left(\frac{j}{s_{0}} + \frac{k}{s_{0}} + 
  (\alpha-\alpha')\frac{q-1}{q} - \frac{\alpha'}{q} +
  \frac{2a-\gamma}{a} \frac{q-1}{q}\right)!
\end{multline*}
where we bound $ m $ by $ k $.

We can then see that there is a suitable positive constant, $ M_{1} $,
independent of $ j $, $ \alpha $, $ \alpha' $, such that
\begin{equation}
\label{eq:k>0}
\left| \partial_{\rho}^{2a -k - \gamma} x^{m-\alpha'}
  \partial_{x}^{m+\alpha-\alpha'} u_{j} \right| \leq
M_{1}^{1+j+\alpha} e^{\tilde{\mu}_{0}
  \rho} \rho^{-(j-\delta)\kappa}  j^{\frac{j}{s_{0}}} \alpha^{\alpha \frac{q-1}{q}} .
\end{equation}
Furthermore we can make the same argument to estimate the term with $
k = 0 $, i.e.
\begin{equation}
\label{eq:k=0}
\left|\partial_{\rho}^{2a - \gamma} \partial_{x}^{\alpha} u_{j}
\right| \leq
M_{1}^{1+j+\alpha} e^{\tilde{\mu}_{0}
  \rho} \rho^{-(j-\delta)\kappa}  j^{\frac{j}{s_{0}}} \alpha^{\alpha \frac{q-1}{q}} .
\end{equation}
Hence, when $ k > 0 $ we get
\begin{multline}
\label{eq:k>0-1}
\binom{\alpha}{\alpha'} \frac{1}{\gamma!} | D_{\rho}^{\gamma}
\psi_{j}| \left| \partial_{\rho}^{2a -k - \gamma} x^{m-\alpha'}
  \partial_{x}^{m+\alpha-\alpha'} u_{j} \right|
\\
\leq
M_{2}^{1+j+\alpha} | D_{\rho}^{\gamma} \psi_{j}|  e^{\tilde{\mu}_{0}
  \rho} \rho^{-(j-\delta)\kappa}  j^{\frac{j}{s_{0}}} \alpha^{\alpha
  \frac{q-1}{q}}
\\
\leq
M_{2}^{1+j+\alpha}  \left(
  \frac{1}{3R_{0}}\right)^{\frac{j}{s_{0}}}  | D_{\rho}^{\gamma}
\psi_{j}|  e^{\tilde{\mu}_{0} \rho} \rho^{\delta \kappa}
\rho^{-j( \kappa - \frac{1}{s_{0}} )} \alpha^{\alpha \frac{q-1}{q}}
\\
\leq
M_{2}^{1+j+\alpha} C_{1}  \left(
  \frac{1}{3R_{0}}\right)^{\frac{j}{s_{0}}}  e^{\tilde{\mu}_{0} \rho} \rho^{\delta \kappa} 
\rho^{- (\frac{\rho}{6R_{0}} - 1) ( \kappa - \frac{1}{s_{0}} )} \alpha^{\alpha
  \frac{q-1}{q}}
\\
\leq
M_{2}^{1+j+\alpha} C_{2}  \left(
  \frac{1}{3R_{0}}\right)^{\frac{j}{s_{0}}}  
e^{-\left( \frac{ \kappa - \frac{1}{s_{0}}}{6R_{0}} 
  \right) \rho \log \rho + \frac{\tilde{\mu}_{0}}{2}\rho} \alpha^{\alpha
  \frac{q-1}{q}} ,
\end{multline}
and an analogous estimate for the term containing \eqref{eq:k=0}.

Plugging the above estimate into \eqref{eq:daS2} and keeping into
account that the summations in $ k $, $ \gamma $, $ \alpha' $ run over
a finite range of indices, we obtain with a new positive constant $
M_{3} $ 
\begin{equation}
\label{eq:daS2-1}
| \partial_{x}^{\alpha} S_{2}(x, \rho) | \leq
\sum_{j\geq 0} M_{3}^{1+j+\alpha}  \left(
  \frac{1}{3R_{0}}\right)^{\frac{j}{s_{0}}}  
e^{-\left( \frac{ \kappa - \frac{1}{s_{0}}}{6R_{0}} 
  \right) \rho \log \rho+ \frac{\tilde{\mu}_{0}}{2}\rho } \alpha^{\alpha
  \frac{q-1}{q}} .
\end{equation}
Thus choosing $ (3 R_{0})^{\frac{1}{s_{0}}} > M_{3} $ we get
\begin{equation}
\label{eq:daS2-2}
| \partial_{x}^{\alpha} S_{2}(x, \rho) | \leq
M_{4}^{1+\alpha} \
e^{-\left( \frac{ \kappa - \frac{1}{s_{0}}}{6R_{0}} 
  \right) \rho \log \rho + \frac{\tilde{\mu}_{0}}{2}\rho } \alpha^{\alpha
  \frac{q-1}{q}} .
\end{equation}
We need the following
\begin{lemma}
\label{lemma:intBs0}
Let $ \mu > 0 $, $ s > 1 $. For any $ \epsilon > 0 $ there is a constant $
C_{\epsilon} > 0 $ such that
\begin{equation}
\label{eq:intBs0}
\int_{0}^{+\infty} e^{- \mu \rho \log \rho} \rho^{s \alpha}  d\rho
\leq C_{\epsilon} \epsilon^{\alpha} \alpha!^{s}.
\end{equation}
\end{lemma}
\begin{proof}[Proof of Lemma \ref{lemma:intBs0}]
Pick a positive $ M $ to be chosen later and write
\begin{eqnarray*}
\int_{0}^{+\infty} e^{- \mu \rho \log \rho} \rho^{s\alpha}  d\rho
 &  =  &
         \int_{0}^{M} e^{- \mu \rho \log \rho} \rho^{s\alpha}  d\rho
+ \int_{M}^{+\infty} e^{- \mu \rho \log \rho} \rho^{s\alpha}
d\rho \\
& = & I_{1} + I_{2}.
\end{eqnarray*}
Consider $ I_{2} $. Because $ e^{- \mu \rho \log \rho} \leq e^{- \mu \rho
  \log M}  $, we get
$$
I_{2} \leq \int_{0}^{+\infty} e^{- \mu \rho \log M} \rho^{s\alpha}
d\rho = \left( \frac{1}{\mu \log M}\right)^{s\alpha + 1} \alpha!^{s}.
$$
Choosing $ \mu^{-s} (\log M)^{-s} \leq \epsilon $ we prove the
assertion for $ I_{2} $.

Consider $ I_{1} $.
$$
I_{1} \leq e^{\mu/e} \frac{M^{s\alpha + 1}}{s \alpha + 1} \leq
e^{\mu/e} M
\frac{\left(\frac{M^{s}}{\epsilon}\right)^{\alpha}}{\alpha!^{s}
} \epsilon^{\alpha}  \alpha!^{s}.
$$
and this implies the assertion also for $ I_{1} $.
\end{proof}
Going back to the derivative of $ f_{2} $ we have
\begin{multline}
\label{eq:df2}
|\partial_{x}^{\alpha} \partial_{y}^{\beta} f_{2}(x, y)|
\\
\leq
\int_{0}^{+\infty}  \rho^{\re r + (2+ \alpha)
  \frac{s_{0}}{q} +\beta s_{0} }
\left| \partial_{x}^{\alpha} S_{2}(x, \rho) \right|_{x =
  x\rho^{\frac{s_{0}}{q}}} d\rho
\\
\leq
M_{4}^{1+\alpha} \ \alpha^{\alpha \frac{q-1}{q}}
\int_{0}^{+\infty}  \rho^{\re r + (2+ \alpha)
  \frac{s_{0}}{q} +\beta s_{0} }
e^{-\left( \frac{ \kappa - \frac{1}{s_{0}}}{6R_{0}} 
  \right) \rho \log \rho + \frac{\tilde{\mu}_{0}}{2}\rho} 
d\rho
\\
\leq
M_{5}^{1+\alpha} \ \alpha^{\alpha \frac{q-1}{q}}
\int_{0}^{+\infty}  \rho^{\left(\frac{\alpha}{q} +\beta\right) s_{0} }
e^{-\left( \frac{ \kappa - \frac{1}{s_{0}}}{6R_{0}} 
  \right) \rho \log \rho }  d\rho .
\end{multline}
Applying Lemma \ref{lemma:intBs0} we obtain that
\begin{multline}
\label{eq:df2-1}
|\partial_{x}^{\alpha} \partial_{y}^{\beta} f_{2}(x, y)|
\leq
M_{5}^{1+\alpha} \ \alpha^{\alpha \frac{q-1}{q}} C_{\epsilon}
\epsilon^{\frac{\alpha}{q} +\beta} \left(\frac{\alpha}{q}
  +\beta\right)!^{s_{0}}
\\
\leq
\tilde{C}_{\epsilon_{1}} \epsilon_{1}^{\alpha + \beta} \alpha^{\alpha
  \left(\frac{q-1}{q} +\frac{s_{0}}{q}\right) } \beta^{\beta s_{0}}
\leq
\tilde{C}_{\epsilon_{1}} \epsilon_{1}^{\alpha + \beta} \alpha^{\alpha
  s_{0}} \beta^{\beta s_{0}} ,
\end{multline}
since the inequality
$$ 
\frac{q-1}{q} + \frac{s_{0}}{q} < s_{0} 
$$
is obviously true, being equivalent to $ s_{0} > 1 $. This proves that
the assertion is true for $ f_{2} $.

Consider now $ f_{1} $:
$$ 
f_{1}(x, y) = \int_{0}^{+\infty} e^{i y \rho^{s_{0}}} \rho^{r + 2 \frac{s_{0}}{q}}
\left[ \sum_{j\geq 0} \sum_{k=0}^{2a} \frac{1}{\rho^{k}} \psi_{j}
  P_{k} u_{j}(x, \rho) \right]_{x = x\rho^{\frac{s_{0}}{q}}} d\rho .
$$
Using Proposition \ref{prop:equiv}, or rather its proof, we may
regroup the terms in the above summation according to the scheme of
\eqref{eq:tr1}, \eqref{eq:tr2}.
\begin{multline}
\label{eq:f1-1}
f_{1}(x, y) = \int_{0}^{+\infty} e^{i y \rho^{s_{0}}} \rho^{r + 2
  \frac{s_{0}}{q}}  
\left[\psi_{0} P_{0} u_{0} + \sum_{j\geq 1} \left( \sum_{k=0}^{\min\{j, 2a\}}
  \frac{1}{\rho^{k}} \psi_{j-k} (1 - \pi) P_{k} u_{j-k} \right)
\right.
\\
+
\sum_{j\geq 1}
\left(\vphantom{\sum_{k=1}^{\min\{j, 2a-1\}}}
\psi_{j} \pi P_{0} u_{j} + \frac{1}{\rho} \psi_{j} \pi P_{1} (1-\pi)
u_{j} + \frac{1}{\rho} \psi_{j-1} \pi P_{1} \pi u_{j-1}
\right.
\\
\left .
\left .
+ \sum_{k=1}^{\min\{j, 2a-1\}} \frac{1}{\rho^{k+1}} \psi_{j-k} \pi
P_{k+1} u_{j-k} \right) \right]_{x = x\rho^{\frac{s_{0}}{q}}} d\rho
\\
=
f_{10}(x, y) + f_{11}(x, y) + f_{12}(x, y) .
\end{multline}
Observe that $ f_{10} = 0 $ since $ P_{0} u_{0} = 0 $ because of
\eqref{eq:u0-}.

Consider $ f_{11} $. For a fixed $ j $, the support of
$$ 
\sum_{k=0}^{\min\{j, 2a\}} \frac{1}{\rho^{k}} \psi_{j-k} (1 - \pi)
P_{k} u_{j-k} 
$$
is contained in the interval $ [3 R_{0}((j-2a)_{+}+1), 6R_{0}(j+1)] $.
In fact for $ \rho \geq 6 R_{0}(j+1) $ the functions $ \psi_{j-k} \equiv
1$ and \eqref{eq:tr1} is satisfied. On the other hand if $ \rho \leq 3
R_{0}((j-2a)_{+}+1) $, then $ \psi_{j-k} \equiv 0 $.

Moreover
\begin{multline}
\label{eq:f11d}
\partial_{x}^{\alpha} \partial_{y}^{\beta} f_{11}(x, y)
=
i^{\beta}
\int_{0}^{+\infty} e^{i y \rho^{s_{0}}} \rho^{r + (2+ \alpha)
  \frac{s_{0}}{q} +\beta s_{0} }
\\
\cdot
\left[ \sum_{j\geq 1}  \sum_{k=0}^{\min\{j, 2a\}}
  \frac{1}{\rho^{k}} \psi_{j-k}  \partial_{x}^{\alpha} (1 - \pi) P_{k} u_{j-k}  \right]_{x =
  x\rho^{\frac{s_{0}}{q}}} d\rho .
\end{multline}
In view of Lemma \ref{lemma:derproj}, arguing as we did before for the
derivatives of $ f_{2} $, we conclude that $ f_{11} \in
\mathscr{B}^{s_{0}}_{g}(\R^{2}) $.

Finally consider $ f_{12} $. We first remark that since $ \pi P_{1}\pi
=0$, due to Lemma \ref{lemma:piPpi},
\begin{multline*}
\psi_{j} \pi P_{0} u_{j} + \frac{1}{\rho} \psi_{j} \pi P_{1} (1-\pi)
u_{j} + \frac{1}{\rho} \psi_{j-1} \pi P_{1} \pi u_{j-1}
\\
+
\sum_{k=1}^{\min\{j, 2a-1\}} \frac{1}{\rho^{k+1}} \psi_{j-k} \pi
P_{k+1} u_{j-k}
\\
=
\psi_{j} \pi P_{0} u_{j} +
\frac{1}{\rho} \psi_{j} \pi P_{1} u_{j}
+
\sum_{k=1}^{\min\{j, 2a-1\}} \frac{1}{\rho^{k+1}} \psi_{j-k} \pi
P_{k+1} u_{j-k} 
\end{multline*}
Again as above we have that, for a fixed $ j $, the $ \rho $-support
of the above quantity is contained in the interval $ [3
R_{0}((j-2a+1)_{+}+1), 6R_{0}(j+1)] $. 

Arguing as above we achieve the proof of the theorem.
\end{proof}
We remark that, as a consequence of Theorem \ref{th:Beu}, we found a
function $ f(x, y) \in \mathscr{B}^{s_{0}}_{g}(\R^{2}) $ such that
\begin{equation}
  \label{eq:PAv=f}
P A(v) = f, \qquad \text{ in } \R^{2}.
\end{equation}

We also need a technical variant of Theorem \ref{th:Beu} adding a
finite order vanishing rate at infinity to the property of being in $
\mathscr{B}^{s_{0}}_{g}(\R^{2}) $. 
\begin{corollary}
\label{cor:Beu-somm}
We use the same notation of Theorem \ref{th:Beu}. Let $ b > 1 $ denote
a fixed positive integer and $ \epsilon > 0 $. Then for every $ k $, $ 1 \leq k \leq b $,
there is a constant $ C_{\epsilon} > 0 $ such that
\begin{equation}
\label{eq:Beu-somm}
\left| \langle y \rangle^{k} \partial_{x}^{\alpha}
  \partial_{y}^{\beta} PA(v)(x, y) \right| \leq C_{\epsilon}
\epsilon^{\alpha+\beta} \alpha!^{s_{0}} \beta!^{s_{0}},
\end{equation}
for every $ \alpha $, $ \beta \in \N \cup \{ 0 \} $ and for any $ (x,
y) \in \R^{2} $.
\end{corollary}
\begin{proof}
The proof is made along the same lines as the proof of Theorem
\ref{th:Beu}, so we are going to give a sketch of it detailing those
parts where it is different from that of the above theorem.

We have, using \eqref{eq:gamma0} with $ 2a $ replaced by $ k $, 
\begin{multline*}
y^{k} \partial_{x}^{\alpha} \partial_{y}^{\beta} \left( P(x, y, D_{x},
  D_{y}) A(v)(x, y) \right)
\\
=
\int_{0}^{+\infty} e^{i y \rho^{s_{0}}}   \sum_{h=0}^{k}
\gamma_{k, h} \frac{(-1)^{h}}{\rho^{k s_{0} - h}} \partial_{\rho}^{h}
\\
\cdot
\left[\rho^{r + (2+\alpha) \frac{s_{0}}{q} + \beta s_{0}}
  \partial_{x}^{\alpha} \left( \sum_{j=0}^{2a} \frac{1}{\rho^{j}}
  P_{j}(x, D_{x}, D_{\rho})
  v(x, \rho) \right) \right]_{x = x\rho^{\frac{s_{0}}{q}}} d\rho .
\end{multline*}
Proceeding as we did for \eqref{eq:drhobin}, \eqref{eq:PA:2},
\eqref{eq:Pj} we obtain 
\begin{multline*}
y^{k} \partial_{x}^{\alpha} \partial_{y}^{\beta} \left( P(x, y, D_{x},
  D_{y}) A(v)(x, y) \right)
\\
=
\int_{0}^{+\infty} e^{i y \rho^{s_{0}}} \rho^{k(1-s_{0}) + N_{\alpha,
    \beta}}
\\
\cdot
\left[ \sum_{j'=0}^{k} \frac{1}{\rho^{j'}} P^{\#}_{j'}(x
  \partial_{x}) \partial_{\rho}^{k-j'} \partial_{x}^{\alpha}
  \sum_{j=0}^{2a} \frac{1}{\rho^{j}}
  P_{j}(x, D_{x}, D_{\rho}) v(x, \rho) \right]_{x = x\rho^{\frac{s_{0}}{q}}} d\rho ,
\end{multline*}
where $ P^{\#}_{j'} $ are polynomials of degree $ j' $ in $ x
\partial_{x} $ whose coefficients are $
\mathscr{O}(\hat{C}^{\alpha+\beta}) $ with $ \hat{C} $ a positive
absolute constant.

Using \eqref{eq:Si} we are reduced to estimating, for $ i= 1, 2 $, the
expressions
$$
\int_{0}^{+\infty} e^{i y \rho^{s_{0}}} \rho^{k(1-s_{0}) + N_{\alpha,
    \beta}}
\left[ \sum_{j'=0}^{k} \frac{1}{\rho^{j'}} P^{\#}_{j'}(x
  \partial_{x}) \partial_{\rho}^{k-j'} \partial_{x}^{\alpha}
  S_{i} (x, \rho) \right]_{x = x\rho^{\frac{s_{0}}{q}}} d\rho .
$$
Consider first the quantity involving $ S_{2} $. The typical summand
in the quantity in square brackets has the form
$$ 
\mathscr{O}(\hat{C}^{\alpha+\beta}) \frac{1}{\rho^{j'}} x^{\gamma}
\partial_{x}^{\gamma+\alpha} \partial_{\rho}^{k-j'} S_{2} ,
$$
where $ \gamma $, $ j' \leq k $.

Arguing as for the derivation of \eqref{eq:daS2-2}, choosing $
\gamma^{\#} $ large enough depending on $ k $, we obtain
\begin{multline}
\label{eq:daS2pp}
\left| \mathscr{O}(\hat{C}^{\alpha+\beta}) x^{\gamma}
\partial_{x}^{\gamma+\alpha} \partial_{\rho}^{k-j'} S_{2}(x, \rho)
\right|
\\
\leq
M_{2}^{1+\alpha+\beta} \
e^{-\left( \frac{ \kappa - \frac{1}{s_{0}}}{6R_{0}} 
  \right) \rho \log \rho + \frac{\tilde{\mu}_{0}}{2}\rho } \alpha^{\alpha
  \frac{q-1}{q}} .
\end{multline}
Then the conclusion follows as in the proof of Theorem \ref{th:Beu}.

Let us examine the term containing $ S_{1} $. As in the proof of
Theorem \ref{th:Beu} we see that
$$ 
S_{1}(x, \rho) = \sum_{j\geq 0} \sum_{k=0}^{2a} \frac{1}{\rho^{k}} \psi_{j}
  P_{k} u_{j},
$$
so that its $ \rho $-support is contained in the interval $ [3
R_{0}((j-2a)_{+}+1), 6R_{0}(j+1)] $. Then the proof goes along the
same lines as that of Theorem \ref{th:Beu}.

This completes the proof of the corollary.
\end{proof}

\section{End of the Proof}
\setcounter{equation}{0}
\setcounter{theorem}{0}
\setcounter{proposition}{0}
\setcounter{lemma}{0}
\setcounter{corollary}{0}
\setcounter{definition}{0}

To finish the proof of Theorem \ref{th:1} we argue by
contradiction. Assume that $ P $ is Gevrey-$ s $ hypoelliptic for an $
s $, with $ 1 \leq s < s_{0} $.

By Theorem 3.1 in \cite{metivier80} it follows from \eqref{eq:PAv=f}
that
\begin{equation}
  \label{eq:AvinB}
  A(v) \in \mathscr{B}^{s_{0}}(\R^{2}) .
\end{equation}
Here is a brief description of what we are going to do in what follows. First we show
that any $ y $-derivative  of $
A(v) $ is integrable over $ \R $ for $ x = 0 $. Then we prove that $
A(v)(0, y) \in \mathscr{B}_{g}^{s_{0}}(\R) $ with some decreasing rate
at infinity, and this allows us to show that for any $ \delta > 0 $
\begin{equation}
  \label{eq:Fav1}
\left | \mathscr{F}(A(v))(0, \eta) \right| \leq C_{\delta} e^{-
  \delta^{-1} |\eta|^{\frac{1}{s_{0}}}} ,
\end{equation}
for a suitable $ C_{\delta} > 0 $.

On the other hand the construction of $ A(v) $ implies that its
Fourier transform satisfies a bound from below of the form
\begin{equation}
  \label{eq:bbelow}
\left | \mathscr{F}(A(v))(0, \eta) \right| \geq C_{0} \eta^{\lambda'}
e^{\tilde{\mu}_{0} \eta^{\frac{1}{s_{0}}}}, \qquad \eta \geq (6
R_{0})^{s_{0}}, 
\end{equation}
for a suitable $ C_{0} > 0 $ and $ \lambda' \in \R $, where $
\tilde{\mu}_{0} $ is defined in 
\eqref{eq:mitlda0}. This is the desired contradiction.

\bigskip

\begin{lemma}
\label{lemma:global}
For any $ \alpha $, $ \beta $ there exists a positive constant $
C_{\alpha, \beta} $ such that 
\begin{equation}
  \label{eq:L2}
(1 + x^{2k} + y^{2k}) \left| \partial_{x}^{\alpha}
  \partial_{y}^{\beta} A(v)(x, y) \right| \leq C_{\alpha, \beta},
\end{equation}
for $ k \leq b $, $ b \in \N $ a fixed integer. In particular, if we
choose $ b $ suitably, the $ (\alpha, \beta) $-derivatives of $ A(v) $
are in $ L^{2}(\R^{2}) $.

\end{lemma}
\begin{proof}
Consider
$$
x^{2k} \partial_{x}^{\alpha} D_{y}^{\beta} A(v)(x, y)
=
\int_{0}^{+\infty} e^{i y \rho^{s_{0}}} \rho^{r+(\alpha-2k)\frac{s_{0}}{q}
+ \beta s_{0} }
 \left(x^{2k} \partial_{x}^{\alpha} v \right)(x\rho^{\frac{s_{0}}{q}}
 , \rho) d\rho .
$$
Proceeding analogously to the proof of Theorem \ref{th:Beu}, from
\eqref{eq:pest} we get
\begin{multline}
  \label{eq:xdalfa}
\left| x^{2k} \partial_{x}^{\alpha} v(x, \rho) \right|
\leq
\sum_{j \geq 0} \psi_{j}(\rho)  \left | x^{2k} \partial_{x}^{\alpha}
  u_{j}(x, \rho) \right|
\\
\leq
\sum_{j \geq 0} \psi_{j}(\rho) e^{\tilde{\mu}_{0} \rho} K_{u}^{1+j+\alpha+2k}
\rho^{-(j-\delta)\kappa} \left( \frac{j}{s_{0}} + \alpha
  \frac{q-1}{q} + \frac{2k}{q}\right)!
\\
\leq
C_{\alpha}  e^{\tilde{\mu}_{0} \rho}
\sum_{j \geq 0}   \tilde{K}_{u}^{1+j}
\left(3 R_{0}(j+1)\right)^{-(j-\delta)\kappa} (j+1)^{\frac{j}{s_{0}}}
\\
\leq
C_{1 \alpha}  e^{\tilde{\mu}_{0} \rho} \sum_{j \geq 0} \tilde{K}_{u}^{1+j}
(j+1)^{j(\frac{1}{s_{0}} -\kappa)}
=
C_{2 \alpha}  e^{\tilde{\mu}_{0}\rho} ,
\end{multline}
whence we conclude.

Consider then $ y^{2k} \partial_{x}^{\alpha} D_{y}^{\beta} A(v)(x, y)
$.
Arguing in the same way as we did when deducing \eqref{eq:PA:3}, and
disregarding the behaviour of the coefficients, we may write
\begin{multline*}
y^{2k} \partial_{x}^{\alpha} D_{y}^{\beta} A(v)(x, y)
\\
=
\int_{0}^{+\infty} e^{i y \rho^{s_{0}}} 
\sum_{h=0}^{2k}
\gamma_{2k, h} \frac{1}{\rho^{2k s_{0} - h}} \partial_{\rho}^{h}
\rho^{r+\alpha\frac{s_{0}}{q}
+ \beta s_{0} }
 \left( \partial_{x}^{\alpha} v \right)(x\rho^{\frac{s_{0}}{q}}
 , \rho) d\rho
\\
=
\int_{0}^{+\infty} e^{i y \rho^{s_{0}}} \rho^{N_{\alpha, \beta}}
\sum_{\ell=0}^{2k} \frac{1}{\rho^{\ell}} \left[ L_{\ell}(x \partial_{x})
\partial_{\rho}^{2k-\ell} \partial_{x}^{\alpha} v(x, \rho) \right]_{x = x
\rho^{\frac{s_{0}}{q}}} d\rho,
\end{multline*}
where $ N_{\alpha, \beta} $ is a suitable complex number independent
of $ \ell $ and $ L_{\ell} $ is a polynomial of degree $ \ell $ with respect to
$ x \partial_{x} $.

In order to estimate the above integral we proceed as before using
\eqref{eq:pest}. Because of \eqref{eq:v} the sum inside the above
integral reads
$$ 
\sum_{j \geq 0} \sum_{\ell=0}^{2k}  \frac{1}{\rho^{\ell}} \partial_{\rho}^{2k-\ell}
L_{\ell}(x \partial_{x}) \partial_{x}^{\alpha} \psi_{j}(\rho) u_{j}(x, \rho).
$$
Writing $ L_{\ell}(x \partial_{x}) = \sum_{m=0}^{\ell} c_{\ell m} (x
\partial_{x})^{m} $ we have, for fixed $ j $, a finite sum ($ \ell
\leq 2k \leq b $ and $ b $ is a fixed positive integer) whose typical
summand has the form
$$ 
\frac{c_{\ell m}}{\rho^{\ell}} \partial_{\rho}^{2k-\ell} x^{m}
\partial_{x}^{m+\alpha} \left( \psi_{j}(\rho) u_{j}(x, \rho)\right).
$$
Choosing $ R_{0} $, $ \gamma^{\#} $ so that $ 4k \leq 2 \gamma^{\#} \leq R_{0} $, we
obtain by \eqref{eq:pest} and Lemma \ref{lemma:psij} that the above
quantity is bounded by
$$ 
C_{2}^{1+\alpha+j} e^{\tilde{\mu}_{0} \rho} \rho^{-(j-\delta)\kappa}
\left( \frac{j}{s_{0}} + \alpha \frac{q-1}{q} \right)! 
$$
We may then argue as in \eqref{eq:xdalfa} and conclude that
\eqref{eq:L2} holds.
\end{proof}

To prove \eqref{eq:Fav1} we prove the following lemma, from which
\eqref{eq:Fav1} will be deduced.
\begin{lemma}
\label{lemma:yd}
For any $ \epsilon > 0 $, there exists a $ C_{\epsilon} > 0 $, such
that for any index $ \alpha $ we have
\begin{equation}
\label{eq:yd}
\left| \langle y \rangle^{a} \partial_{y}^{\alpha}
  \partial_{x}^{\beta} A(v)(0, y) \right|
\leq C_{\epsilon} \epsilon^{\alpha} \alpha!^{s_{0}},
\end{equation}
where $ \beta = 0, 1 $ and $ \langle y \rangle = ( 1 + y^{2})^{\frac{1}{2}} $.
\end{lemma}
\begin{proof}
Let $ \chi $ denote a smooth function in $ G^{s'}(\R) $, $ s' < s_{0} $, such that $
\chi(t) = 0 $ for $ |t| \leq 1 $ and $ \chi(t) = 1 $ for $ |t| \geq 2
$.

Observe that, by \eqref{eq:PAv=f}, 
$$ 
P \chi(y) A(v) = \chi(y) f - [P, \chi] A(v) .
$$
By Corollary \ref{cor:Beu-somm} and formula \eqref{eq:AvinB}, if we denote by $ g $ the right hand side of
the above equation and $ \psi(x) \in G_{0}^{s'}(\R) $ we have that
\begin{equation}
  \label{eq:Bs0PA}
\left| \langle y \rangle^{a} \partial_{y}^{\alpha}
  \partial_{x}^{\beta} (\psi g)(x, y) \right|  \leq C_{\delta}
    \delta^{\alpha+\beta} \alpha!^{s_{0}} \beta!^{s_{0}},
\end{equation}
for every $ (x, y) \in \R^{2} $ and $ \delta > 0 $. From the above
inequality it readily follows that, possibly renaming $ C_{\delta} $,
\begin{equation}
\label{eq:Bs0PA2}
\| \partial_{y}^{\alpha}
  \partial_{x}^{\beta} (\psi g)(x, y) \|_{0}  \leq C_{\delta}
    \delta^{\alpha+\beta} \alpha!^{s_{0}} \beta!^{s_{0}}
\end{equation}
We are going to use the maximal estimate for the above equation:
\begin{equation}
  \label{eq:maximal}
\sum_{i=1}^{3} \|X_{i} u \|_{0}^{2} \leq C \left( \langle P u, u
  \rangle + \| y^{a-1} u \|_{0}^{2} \right),
\end{equation}
where $ u \in L^{2}(\R^{2}) $, $ Pu \in L^{2}(\R^{2}) $ and $ y^{a-1}
u \in L^{2}(\R^{2}) $, and $ X_{1}$  $ = D_{x} $, $ X_{2} = x^{q-1}D_{y}
$, $ X_{3} = y^{a} D_{y} $.

Let $ \phi=\phi(x) \in G^{s'}_{0}(\R) $ denote a cutoff function near the
origin. 

We want to show that for any $ \epsilon >0 $, for any positive
integers $ \alpha $, $ \beta $, with $ 0 \leq \beta \leq N $, $ N $
arbitrarily fixed natural number, there exists a positive constant $
C_{\epsilon} $ such that
\begin{equation}
\label{eq:XA}
\| X_{i} \partial_{y}^{\alpha} \partial_{x}^{\beta} \phi  A(v)
\|_{0} \leq C_{\epsilon} \epsilon^{\alpha} \alpha!^{s_{0}},
\end{equation}
$ i = 1, 2, 3 $. Since $ A(v) \in \mathscr{B}^{s_{0}}(\R^{2}) $, it is
enough to prove \eqref{eq:XA} for $ A'(v) = \chi A(v) $.

Actually we are going to prove that for any $ \epsilon_{1} >0 $, $
\epsilon_{2} > 0 $,  for any positive
integers $ \alpha $ and $ \beta $, $ 0 \leq \beta \leq N $, there
exists a positive constant $ C_{\epsilon_{1} \epsilon_{2}} $ such that
\begin{equation}
\label{eq:XA'}
\sum_{i=1}^{3} \|X_{i} \partial_{y}^{\alpha} \partial_{x}^{\beta}
\phi^{(k)} A'(v)) \|_{0} \leq C_{\epsilon_{1} \epsilon_{2}}
\epsilon^{\alpha+\beta}_{1} \epsilon_{2}^{k}
(\alpha + \beta + k)!^{s_{0}}.
\end{equation}
The estimate \eqref{eq:XA} then follows keeping into account that $
\beta $ takes a finite number of values and choosing $ \epsilon_{1} $,
$ \epsilon_{2} $ suitably.

We proceed by induction on $ \alpha + \beta $. Actually
\eqref{eq:XA'} is true when $ \alpha + \beta = 0 $. In fact we use Lemma
\ref{lemma:global} as well as the fact that $ \phi \in G_{0}^{s'}(\R)
\subset \mathscr{B}^{s_{0}}_{g}(\R) $. 

Consider
\begin{equation}
  \label{eq:ori.norm}
\sum_{i=1}^{3} \|X_{i} \partial_{y}^{\alpha} \partial_{x}^{\beta}
\phi^{(k)} A'(v) \|_{0}^{2} .
\end{equation}
By \eqref{eq:maximal} we have
\begin{multline}
  \label{eq:max1}
\sum_{i=1}^{3} \|X_{i} \partial_{y}^{\alpha} \partial_{x}^{\beta}
\phi^{(k)}  A'(v)) \|_{0}^{2}
\\
\leq
C \left( \langle P \partial_{y}^{\alpha} \partial_{x}^{\beta}
\phi^{(k)}  A'(v) , \partial_{y}^{\alpha} \partial_{x}^{\beta}
\phi^{(k)} A'(v) \rangle
\right.
\\
\left.
+ \| y^{a-1} \partial_{y}^{\alpha} \partial_{x}^{\beta}
\phi^{(k)} A'(v) \|_{0}^{2} \right)
\\
\leq
C_{1} \left(
\underbrace{\langle \partial_{y}^{\alpha} \partial_{x}^{\beta} \phi^{(k)} P( A'(v)),
  \partial_{y}^{\alpha}  \partial_{x}^{\beta} \phi^{(k)} A'(v) \rangle}_{B_{1}}
\right .
  \\
  +
 2 \sum_{i=1}^{3} \Big( \underbrace{\langle [X_{i} , \partial_{y}^{\alpha}
  \partial_{x}^{\beta} \phi^{(k)} ] A'(v), X_{i}^{*} \partial_{y}^{\alpha} \partial_{x}^{\beta}
  \phi^{(k)}  A'(v) \rangle}_{B_{2i}} 
\\
+
\underbrace{%
\langle [X_{i} , [ X_{i} , \partial_{y}^{\alpha}
  \partial_{x}^{\beta} \phi^{(k)} ] ]  A'(v) ,  \partial_{y}^{\alpha} \partial_{x}^{\beta}
  \phi^{(k)} A'(v) \rangle}_{B_{3i}} \Big)
  \\
  +
\underbrace{%
    \langle [i a y^{2a-1} D_{y} , \partial_{y}^{\alpha}
  \partial_{x}^{\beta} \phi^{(k)} ] A'(v), \partial_{y}^{\alpha} \partial_{x}^{\beta}
  \phi^{(k)}  A'(v) \rangle}_{B_{2}'}
  \\
  \left .
  +
    \underbrace{%
 \| y^{a-1} \partial_{y}^{\alpha} \partial_{x}^{\beta}
 \phi^{(k)}  A'(v)) \|_{0}^{2}}_{B_{4}} 
 \right) .
\end{multline}
Let us start with $ B_{1} $. Let $ \psi \in G^{s'}_{0}(\R) $ be such
that $ \phi \psi = \phi $. Then,
\begin{multline*}
| B_{1}| \leq \left( \| \partial_{y}^{\alpha} \partial_{x}^{\beta} \phi^{(k)}
P( A'(v)) \|_{0} + \| \partial_{y}^{\alpha}  \partial_{x}^{\beta}
\phi^{(k)} A'(v) \|_{0}\right)^{2}
\\
\leq
C_{2}
\left(
\sum_{\ell=0}^{\beta} \binom{\beta}{\ell} \| 
\phi^{(k+\ell)} \partial_{y}^{\alpha} \partial_{x}^{\beta-\ell} \psi
P( A'(v)) \|_{0}
\right.
\\
\left.
+ \| X \partial_{y}^{\alpha'}  \partial_{x}^{\beta'}
\phi^{(k)} A'(v) \|_{0}
\vphantom{\sum_{\ell=0}^{\beta}}
\right)^{2} ,
\end{multline*}
where $ X $ denotes the vector field $ X_{1} $ or $ X_{3} $ according
to the fact that $ \beta > 0 $ or $ \alpha > 0 $ respectively,
since on the support of $ A'(v) $, $ |y| \geq c > 0 $ and $
\alpha'+\beta' = \alpha + \beta -1 $. 

Consider the first term in the sum above. By assumption
\eqref{eq:Bs0PA2} we have 
\begin{multline*}
\sum_{\ell=0}^{\beta} \binom{\beta}{\ell} \| 
\phi^{(k+\ell)} \partial_{y}^{\alpha} \partial_{x}^{\beta-\ell} \psi
P( A'(v)) \|_{0}
\\
\leq
\sum_{\ell=0}^{\beta} \binom{\beta}{\ell} C_{\delta_{2}} \delta_{2}^{k+\ell}
(k+\ell)!^{s_{0}} C_{\delta_{1}} \delta_{1}^{\alpha+\beta-\ell}
(\alpha+ \beta -\ell)!^{s_{0}}
\\
\leq
\frac{1}{\nu C_{2}}
C_{\epsilon_{1} \epsilon_{2}} \epsilon_{1}^{\alpha+\beta} \epsilon_{2}^{k} (\alpha+\beta+k)!^{s_{0}} ,
\end{multline*}
choosing $ \delta_{1} < \frac{\epsilon_{1}}{2} $, $ \delta_{2} <
\min\{\epsilon_{2}, \delta_{1}\} $ and $
C_{\epsilon_{1} \epsilon_{2}} >  \nu C_{\delta_{1}} C_{\delta_{2}} C_{2}$.

As for the second term we have, by applying the inductive assumption
\eqref{eq:XA'}, that 
\begin{multline*}
\| X \partial_{y}^{\alpha'}  \partial_{x}^{\beta'}
\phi^{(k)} A'(v) \|_{0} \leq 
C_{\epsilon_{1} \epsilon_{2}} \epsilon_{1}^{\alpha'+\beta'} \epsilon_{2}^{k}
(\alpha'+\beta'+k)!^{s_{0}}
\\
=
C_{\epsilon_{1} \epsilon_{2}} \frac{1}{(\alpha+\beta+k)^{s_{0}} \epsilon_{1}}
\epsilon_{1}^{\alpha+\beta} \epsilon_{2}^{k}  (\alpha+\beta+k)!^{s_{0}}
\\
\leq
\frac{1}{\nu} C_{\epsilon_{1} \epsilon_{2}}
\epsilon_{1}^{\alpha+\beta} \epsilon_{2}^{k} (\alpha+\beta+k)!^{s_{0}} ,
\end{multline*}
provided
\begin{equation}
\label{eq:abk1}
\alpha+\beta+k > (\nu \epsilon_{1}^{-1})^{\frac{1}{s_{0}}} .
\end{equation}
By Lemma \ref{lemma:global} we obtain the same estimate for any $
\alpha $, $ \beta $, $ 0 \leq \beta \leq N $, $ k $, if $
C_{\epsilon_{1} \epsilon_{2}} $ is chosen large enough.

Consider now $ B_{2i} $. Remark that if $ i = 1 $, $ [X_{1} , \partial_{y}^{\alpha}
  \partial_{x}^{\beta} \phi^{(k)} ] = \partial_{y}^{\alpha}
  \partial_{x}^{\beta}$ $ \frac{1}{i} \phi^{(k+1)}$. 
Hence for $ i = 1 $ we have
\begin{multline*}
\langle [X_{1} , \partial_{y}^{\alpha}
  \partial_{x}^{\beta} \phi^{(k)} ] A'(v), X_{1}^{*} \partial_{y}^{\alpha} \partial_{x}^{\beta}
  \phi^{(k)}  A'(v) \rangle
  \\
  =
 \langle \partial_{y}^{\alpha}
  \partial_{x}^{\beta} \frac{1}{i} \phi^{(k+1)} A'(v), X_{1}
  \partial_{y}^{\alpha} \partial_{x}^{\beta} \phi^{(k)}  A'(v) \rangle .
\end{multline*}
Thus
\begin{multline*}
\left |\langle [X_{1} , \partial_{y}^{\alpha}
  \partial_{x}^{\beta} \phi^{(k)} ] A'(v), X_{1}^{*} \partial_{y}^{\alpha} \partial_{x}^{\beta}
  \phi^{(k)}  A'(v) \rangle \right|
  \\
  \leq
 \| \partial_{y}^{\alpha}
  \partial_{x}^{\beta}  \phi^{(k+1)} A'(v)\|_{0}  \cdot \| X_{1}
  \partial_{y}^{\alpha} \partial_{x}^{\beta} \phi^{(k)}  A'(v) \|_{0}
\\
\leq
\mu \| X_{1}
  \partial_{y}^{\alpha} \partial_{x}^{\beta} \phi^{(k)}  A'(v)
  \|_{0}^{2} +
\frac{1}{\mu} \| \partial_{y}^{\alpha}
  \partial_{x}^{\beta}  \phi^{(k+1)} A'(v)\|_{0}^{2} .
\end{multline*}
The first term is absorbed on the left hand side of \eqref{eq:max1} if
$ \mu $ is small enough. We are left with the second term. As before
it can be bounded by
$$ 
C_{3} \frac{1}{\mu} \| X \partial_{y}^{\alpha'}
  \partial_{x}^{\beta'}  \phi^{(k+1)} A'(v)\|_{0}^{2},
$$
where $ \alpha' + \beta' = \alpha+\beta-1 $, and $ X $ denotes either
$ X_{1} $ or $ X_{3} $ according to the fact that $ \beta > 0 $ or $
\alpha >0 $ respectively, since, on the support of $ \phi^{(k+1)}
A'(v)$, $ y $ is bounded away from zero.

The above norm is bounded by
$$ 
C_{3} \frac{1}{\mu} C_{\epsilon_{1}\epsilon_{2}}
\epsilon_{1}^{\alpha'+\beta'} \epsilon_{2}^{k+1} (\alpha'+\beta'+k+1)!^{s_{0}},
$$
which becomes
$$ 
\left( C_{3} \frac{1}{\mu} \frac{\epsilon_{2}}{\epsilon_{1}} \right)  C_{\epsilon_{1}\epsilon_{2}}
\epsilon_{1}^{\alpha+\beta} \epsilon_{2}^{k} (\alpha+\beta+k)!^{s_{0}} ,
$$
so that, choosing
\begin{equation}
\label{eq:e1}
\frac{\epsilon_{2}}{\epsilon_{1}} < \frac{1}{\nu} \frac{\mu}{C_{3}} ,
\end{equation}
we obtain, modulo terms absorbed on the left,
\begin{multline}
\label{eq:B21}
\left |\langle [X_{1} , \partial_{y}^{\alpha}
  \partial_{x}^{\beta} \phi^{(k)} ] A'(v), X_{1}^{*} \partial_{y}^{\alpha} \partial_{x}^{\beta}
  \phi^{(k)}  A'(v) \rangle \right|
  \\
  \leq
  \frac{1}{\nu}
  C_{\epsilon_{1}\epsilon_{2}}
\epsilon_{1}^{\alpha+\beta} \epsilon_{2}^{k} (\alpha+\beta+k)!^{s_{0}} .
\end{multline}
Consider now $ B_{22} $. Since $ [X_{2} , \partial_{y}^{\alpha}
\partial_{x}^{\beta} \phi^{(k)} ] =  \partial_{y}^{\alpha} [ x^{q-1} , 
\partial_{x}^{\beta} ] \phi^{(k)}  D_{y}  $, we may write, for 
suitable positive constants $ C_{q} $, $ C_{1q} $, depending only on $ q $,
\begin{multline*}
\left | \langle [X_{2} , \partial_{y}^{\alpha}
  \partial_{x}^{\beta} \phi^{(k)} ] A'(v), X_{2} \partial_{y}^{\alpha} \partial_{x}^{\beta}
  \phi^{(k)}  A'(v) \rangle \right|
\\
\leq
C_{q} \sum_{\ell=1}^{\min\{\beta, q-1\}} \binom{\beta}{\ell}
\left| \langle
  x^{q-1-\ell} \partial_{x}^{\beta-\ell}
  \phi^{(k)} \partial_{y}^{\alpha+1} A'(v), X_{2} \partial_{y}^{\alpha} \partial_{x}^{\beta}
  \phi^{(k)}  A'(v) \rangle \right|
\\
\leq
C_{1q} \sum_{\ell=1}^{\min\{\beta, q-1\}} \frac{\beta!}{(\beta -
  \ell)!}
\| x^{q-1-\ell} \partial_{x}^{\beta-\ell}
 \partial_{y}^{\alpha+1}  \phi^{(k)}  A'(v) \|_{0} \cdot  \| X_{2} \partial_{y}^{\alpha} \partial_{x}^{\beta}
 \phi^{(k)}  A'(v)\|_{0}
 \\
 \leq
 C_{1q} \left(
   \vphantom{\left( \sum_{\ell=1}^{\min\{\beta, q-1\}}\right.}
   \mu \| X_{2} \partial_{y}^{\alpha} \partial_{x}^{\beta}
   \phi^{(k)}  A'(v)\|_{0}^{2}
\right .
\\
\left .
 +
 \frac{1}{\mu}  \left( \sum_{\ell=1}^{\min\{\beta, q-1\}} \frac{\beta!}{(\beta -
  \ell)!}
\| x^{q-1-\ell} \partial_{x}^{\beta-\ell}
 \partial_{y}^{\alpha+1}  \phi^{(k)}  A'(v) \|_{0} \right)^{2} \right)
\end{multline*}
The first term is absorbed on the left hand side of \eqref{eq:max1}
provided $ \mu $ is small but otherwise fixed,
so that we are left with the estimate of the square of the sum. The
latter can be bounded as
\begin{multline}
  \label{eq:commX2}
C_{4}  \sum_{\ell=1}^{\min\{\beta, q-1\}} \beta (\beta-1) \cdots (\beta-\ell+1)
\| X_{3} \partial_{x}^{\beta-\ell} \partial_{y}^{\alpha}  \phi^{(k)}
A'(v) \|_{0}
\\
\leq
C_{4} \sum_{\ell=1}^{\min\{\beta, q-1\}} \beta (\beta-1) \cdots (\beta-\ell+1)
C_{\epsilon_{1}\epsilon_{2}} \epsilon_{1}^{\beta-\ell+\alpha}
\epsilon_{2}^{k} (\alpha + \beta -\ell + k)!^{s_{0}}
\\
\leq
C_{\epsilon_{1}\epsilon_{2}} \epsilon_{1}^{\beta+\alpha}
\epsilon_{2}^{k} (\alpha + \beta+ k)!^{s_{0}}
\cdot
C_{4}
\sum_{\ell=1}^{\min\{\beta, q-1\}}
\left(\frac{1}{(\alpha+\beta+k-(q-1))^{s_{0}-1}
    \epsilon_{1}}\right)^{\ell}
\\
\leq
\frac{1}{\nu} \sqrt{\mu} C_{\epsilon_{1}\epsilon_{2}} \epsilon_{1}^{\beta+\alpha}
\epsilon_{2}^{k} (\alpha + \beta+ k)!^{s_{0}},
\end{multline}
provided
\begin{equation}
\label{eq:abk2}
\alpha+\beta+k -(q-1) \geq \left( \frac{2 \nu C_{4}}{\sqrt{\mu}
    \epsilon_{1}}\right)^{\frac{1}{s_{0}-1}} .
\end{equation}
By Lemma \ref{lemma:global} we obtain the same estimate of $ B_{22} $ for any $
\alpha $, $ \beta $, $ k $, if $ C_{\epsilon_{1} \epsilon_{2}} $ is chosen large
enough.

Consider now $ B_{23} $ in \eqref{eq:max1}.
Since $[X_{3} , \partial_{y}^{\alpha}
\partial_{x}^{\beta} \phi^{(k)} ] =  [ X_{3} , 
\partial_{y}^{\alpha} ] $  $\partial_{x}^{\beta} \phi^{(k)}
= [ y^{a} , \partial_{y}^{\alpha} ] D_{y} \partial_{x}^{\beta}
\phi^{(k)} $. Hence
\begin{multline*}
\left | \langle [X_{3} , \partial_{y}^{\alpha}
  \partial_{x}^{\beta} \phi^{(k)} ] A'(v), X_{3}^{*} \partial_{y}^{\alpha} \partial_{x}^{\beta}
  \phi^{(k)}  A'(v) \rangle \right|
\\
\leq
C_{1a} \sum_{\ell=1}^{\min\{\alpha, a\}} \binom{\alpha}{\ell}
\left | \langle y^{a-\ell} \partial_{y}^{\alpha-\ell+1} 
  \partial_{x}^{\beta} \phi^{(k)} A'(v), X_{3}^{*} \partial_{y}^{\alpha} \partial_{x}^{\beta}
  \phi^{(k)}  A'(v) \rangle \right|
\\
\leq
C_{1a} \sum_{\ell=1}^{\min\{\alpha, a\}} \binom{\alpha}{\ell} \|
y^{a-\ell} \partial_{y}^{\alpha-\ell+1} \partial_{x}^{\beta}
\phi^{(k)} A'(v) \|_{0}
\| X_{3}^{*} \partial_{y}^{\alpha} \partial_{x}^{\beta}
\phi^{(k)}  A'(v) \|_{0}
\\
\leq
C_{2a} \left(
\vphantom{\left(\sum_{\ell=1}^{\min\{\alpha, a\}} \right.}
\mu \| X_{3} \partial_{y}^{\alpha} \partial_{x}^{\beta}
\phi^{(k)}  A'(v) \|_{0}^{2} + \mu \|y^{a-1} \partial_{y}^{\alpha} \partial_{x}^{\beta}
\phi^{(k)}  A'(v) \|_{0}^{2}
\right .
\\
\left.
+  \frac{1}{\mu}
\left(
\sum_{\ell=1}^{\min\{\alpha, a\}} \binom{\alpha}{\ell} \|
y^{a-\ell} \partial_{y}^{\alpha-\ell+1} \partial_{x}^{\beta}
\phi^{(k)} A'(v) \|_{0}
\right)^{2} \right)
\end{multline*}
The first summand is absorbed on the left of \eqref{eq:max1}, so that
we have to treat the second as well as the sum. Consider the second
term above.

If $ \alpha = 0 $ we use Lemma \ref{lemma:global} and the fact that $
\beta $ takes a finite number of values to conclude that the second
term verifies the desired estimate. Assume $ \alpha > 0 $. Then
\begin{multline*}
C_{2a} \mu \|y^{a-1} \partial_{y}^{\alpha} \partial_{x}^{\beta}
\phi^{(k)}  A'(v) \|_{0}^{2}
%
\leq
C_{5} \|X_{3} \partial_{y}^{\alpha-1} \partial_{x}^{\beta}
\phi^{(k)}  A'(v) \|_{0}^{2}
\\
\leq
C_{5} \left(C_{\epsilon_{1} \epsilon_{2}} \epsilon_{1}^{\alpha+\beta-1}
  \epsilon_{2}^{k} (\alpha+\beta+k-1)!^{s_{0}} \right)^{2}
\\
=
C_{5}
\frac{1}{(\alpha+\beta+k)^{2s_{0}} \epsilon_{1}^{2} }
\left(C_{\epsilon_{1} \epsilon_{2}} \epsilon_{1}^{\alpha+\beta}
  \epsilon_{2}^{k} (\alpha+\beta+k)!^{s_{0}} \right)^{2}
\\
\leq
\frac{1}{\nu^{2}} \left(C_{\epsilon_{1} \epsilon_{2}} \epsilon_{1}^{\alpha+\beta}
  \epsilon_{2}^{k} (\alpha+\beta+k)!^{s_{0}} \right)^{2},
\end{multline*}
provided
\begin{equation}
\label{eq:abk3}
\alpha+\beta+k \geq \left( \frac{\nu \sqrt{C_{5}}}{\epsilon_{1}}
\right)^{\frac{1}{s_{0}}} .
\end{equation}
If \eqref{eq:abk3} is not satisfied we choose a larger $
C_{\epsilon_{1} \epsilon_{2}} $ as we did above, in view of Lemma
\ref{lemma:global}.

The same argument treats also $ B_{4} $ and $ B_{2}' $.

Consider now
$$ 
\sqrt{\frac{C_{2a}}{\mu}}
\sum_{\ell=1}^{\min\{\alpha, a\}} \binom{\alpha}{\ell}  \|
y^{a-\ell} \partial_{y}^{\alpha-\ell+1} \partial_{x}^{\beta}
\phi^{(k)} A'(v) \|_{0} . 
$$
The above sum is bounded as
\begin{multline}
  \label{eq:commX3}
\sqrt{\frac{C_{2a}}{\mu}} \sum_{\ell=1}^{\min\{\alpha, a\}}
\binom{\alpha}{\ell}
\| X_{3} \partial_{y}^{\alpha-\ell} \partial_{x}^{\beta}
\phi^{(k)} A'(v) \|_{0}
\\
\leq
\sqrt{\frac{C_{2a}}{\mu}} \sum_{\ell=1}^{\min\{\alpha, a\}}
\binom{\alpha}{\ell}
C_{\epsilon_{1} \epsilon_{2}} \epsilon_{1}^{\alpha+\beta-\ell}
\epsilon_{2}^{k} (\alpha+\beta+k-\ell)!^{s_{0}}
\\
\leq
C_{\epsilon_{1} \epsilon_{2}} \epsilon_{1}^{\alpha+\beta}
\epsilon_{2}^{k} (\alpha+\beta+k)!^{s_{0}}
\\
\cdot
\sqrt{\frac{C_{2a}}{\mu}}
\sum_{\ell=1}^{\min\{\alpha, a\}}
\left(\frac{1}{(\alpha+\beta+k-a)^{s_{0}-1} \epsilon_{1}}
\right)^{\ell} ,
\end{multline}
whence we conclude as for $ B_{22} $ (see \eqref{eq:abk2}).

\bigskip

We are left with the estimate of $ B_{3i} $ in \eqref{eq:max1}. Let us
start with $ B_{31} $. Since
$$
[X_{1} , [ X_{1} , \partial_{y}^{\alpha}
  \partial_{x}^{\beta} \phi^{(k)} ] ] = [X_{1} ,\partial_{y}^{\alpha}
  \partial_{x}^{\beta} \phi^{(k+1)} ] =   \partial_{y}^{\alpha}
  \partial_{x}^{\beta} \phi^{(k+2)} ,
$$
we have
\begin{multline*}
\left| \langle [X_{1} , [ X_{1} , \partial_{y}^{\alpha}
  \partial_{x}^{\beta} \phi^{(k)} ] ]  A'(v) ,  \partial_{y}^{\alpha} \partial_{x}^{\beta}
  \phi^{(k)} A'(v) \rangle \right|
\\
=
\left| \langle  \partial_{y}^{\alpha}
  \partial_{x}^{\beta} \phi^{(k+2)}   A'(v) ,  \partial_{y}^{\alpha} \partial_{x}^{\beta}
  \phi^{(k)} A'(v) \rangle \right| .
\end{multline*}
We may assume that $ \alpha+\beta \geq 2 $, since otherwise bounding
the original norm \eqref{eq:ori.norm} is easily done by choosing $
C_{\epsilon_{1} \epsilon_{2}} $ sufficiently large. Thus taking one
derivative to the right hand side the above scalar product becomes
$$
\left| \langle  \partial_{y}^{\alpha'}
  \partial_{x}^{\beta'} \phi^{(k+2)}   A'(v) ,  \partial_{y}^{\alpha''} \partial_{x}^{\beta''}
  \phi^{(k)} A'(v) \rangle \right| ,
$$
where $ \alpha' + \beta' = \alpha + \beta - 1 $ and $ \alpha'' +
\beta'' = \alpha + \beta + 1 $. Next we may write
\begin{multline*}
\left| \langle  \partial_{y}^{\alpha'}
  \partial_{x}^{\beta'} \phi^{(k+2)}   A'(v) ,  \partial_{y}^{\alpha''} \partial_{x}^{\beta''}
  \phi^{(k)} A'(v) \rangle \right|
\\
\leq
\| \partial_{y}^{\alpha'} \partial_{x}^{\beta'} \phi^{(k+2)}   A'(v)
\|_{0} \cdot \| \partial_{y}^{\alpha''} \partial_{x}^{\beta''}
\phi^{(k)} A'(v) \|_{0}
\\
\leq
M_{1} \left(
\frac{1}{\mu} \| X \partial_{y}^{\alpha_{-}} \partial_{x}^{\beta_{-}} \phi^{(k+2)}   A'(v)
\|_{0}^{2} + \mu \| X \partial_{y}^{\alpha} \partial_{x}^{\beta}
\phi^{(k)} A'(v) \|_{0}^{2} \right) ,
\end{multline*}
where $ \mu $ is a small fixed positive number, $ \alpha_{-}+
\beta_{-} = \alpha + \beta - 2 $.
Moreover in the above expression $ X $ denotes either the
vector field $ X_{1} $ or the field $ X_{3} $, since the $ y $-support
of the integrated function is bounded away from zero. The second term
above can be absorbed on the left of \eqref{eq:max1} provided $ \mu $
is chosen small but finite, so that we have to estimate the first, to
which we may apply the inductive hypothesis:
\begin{multline*}
\frac{M_{1}^{\frac{1}{2}}}{\mu^{\frac{1}{2}}} \| X \partial_{y}^{\alpha_{-}} \partial_{x}^{\beta_{-}} \phi^{(k+2)}   A'(v)
\|_{0}
\leq
\frac{M_{1}^{\frac{1}{2}}}{\mu^{\frac{1}{2}}} C_{\epsilon_{1} \epsilon_{2}} \epsilon_{1}^{\alpha_{-} +
\beta_{-}} \epsilon_{2}^{k+2} (\alpha_{-} + \beta_{-} + k +
2)!^{s_{0}}
\\
=
\frac{1}{\nu} C_{\epsilon_{1} \epsilon_{2}} \epsilon_{1}^{\alpha+\beta}
\epsilon_{2}^{k} (\alpha+\beta+k)!^{s_{0}} \cdot
\frac{M_{1}^{\frac{1}{2}} \nu \epsilon_{2}^{2}}{\mu^{\frac{1}{2}} \epsilon_{1}^{2}} ,
\end{multline*}
which gives the desired estimate provided $ \epsilon_{1} $, $
\epsilon_{2} $ are such that
\begin{equation}
\label{eq:abk4}
\frac{M_{1}^{\frac{1}{2}} \nu \epsilon_{2}^{2}}{\mu^{\frac{1}{2}} \epsilon_{1}^{2}} < 1.
\end{equation}
Consider now $ B_{32} $.  Since
\begin{multline*}
[X_{2}, [X_{2} ,
\partial_{y}^{\alpha} \partial_{x}^{\beta} \phi^{(k)} ] ] = [ X_{2},
\partial_{y}^{\alpha} [x^{q-1}, \partial_{x}^{\beta}] \phi^{(k)} D_{y}
]
\\
=\partial_{y}^{\alpha} [ x^{q-1} , [ x^{q-1}, \partial_{x}^{\beta}] ]
\phi^{(k)} D_{y}^{2}
\\
=
\partial_{y}^{\alpha} \sum_{\ell_{1}=1}^{\min\{\beta, q-1\}}
\binom{\beta}{\ell_{1}}  (\ad \partial_{x})^{\ell_{1}}(x^{q-1})
[ \partial_{x}^{\beta-\ell_{1}} , x^{q-1} ] \phi^{(k)} D_{y}^{2}
\\
=
- \sum_{\ell_{1}=1}^{\min\{\beta, q-1\}}
\sum_{\ell_{2}=1}^{\min\{\beta-\ell_{1}, q-1\}}
\binom{\beta}{\ell_{1}} \binom{\beta - \ell_{1}}{\ell_{2}}
\\
\cdot
(\ad \partial_{x})^{\ell_{1}}(x^{q-1})
(\ad \partial_{x})^{\ell_{2}}(x^{q-1})
\partial_{x}^{\beta-\ell_{1} -\ell_{2}}
\phi^{(k)} \partial_{y}^{\alpha+2} .
\end{multline*}
Hence, $ \beta $ being bound by $ N $,
\begin{multline*}
\left| \langle [X_{2} , [ X_{2} , \partial_{y}^{\alpha}
  \partial_{x}^{\beta} \phi^{(k)} ] ]  A'(v) ,  \partial_{y}^{\alpha} \partial_{x}^{\beta}
  \phi^{(k)} A'(v) \rangle \right|
\\
\leq
M_{2} \sum_{\ell=2}^{\min\{\beta, 2(q-1)\}} \left| \langle 
  \partial_{x}^{\beta - \ell} \partial_{y}^{\alpha+2}  \phi^{(k)}  A'(v) ,  \partial_{y}^{\alpha} \partial_{x}^{\beta}  \phi^{(k)} A'(v) \rangle \right|
\\
\leq
M_{3}
\sum_{\ell=2}^{\min\{\beta, 2(q-1)\}} \left ( \| X_{3}
  \partial_{x}^{\beta - \ell} \partial_{y}^{\alpha+1}  \phi^{(k)}
  A'(v) \|_{0}^{2}
\right.
\\
\left.
  + \| X_{3} \partial_{y}^{\alpha-1}
  \partial_{x}^{\beta}  \phi^{(k)} A'(v) \|_{0}^{2} \right) .
\end{multline*}
We may apply the inductive hypothesis to both summands above and treat
each of them as we did in \eqref{eq:commX2}.

Finally consider $ B_{33} $. We may suppose that $ \alpha > 0 $, otherwise
the double commutator is identically zero. We have
\begin{multline*}
[X_{3}, [X_{3} ,
\partial_{y}^{\alpha} \partial_{x}^{\beta} \phi^{(k)} ] ] = -
[ y^{a} \partial_{y} , [ y^{a} , \partial_{y}^{\alpha} ] \partial_{y} ]
\partial_{x}^{\beta} \phi^{(k)}
\\
=
\sum_{\ell=2}^{\min\{\alpha, 2a\}} c_{\alpha \ell} \alpha^{\ell}
y^{2a-\ell} \partial^{\alpha+2-\ell}_{y} \partial_{x}^{\beta}
\phi^{(k)} ,
\end{multline*}
where the $ c_{\alpha \ell} $ are absolute constants. Hence
\begin{multline*}
\left| \langle [X_{3} , [ X_{3} , \partial_{y}^{\alpha}
  \partial_{x}^{\beta} \phi^{(k)} ] ]  A'(v) ,  \partial_{y}^{\alpha} \partial_{x}^{\beta}
  \phi^{(k)} A'(v) \rangle \right|
\\
\leq
\sum_{\ell=2}^{\min\{\alpha, 2a\}} |c_{\alpha \ell}| \alpha^{\ell}
\left| \langle y^{2a-\ell} \partial^{\alpha+2-\ell}_{y}
  \partial_{x}^{\beta} \phi^{(k)} A'(v), \partial_{y}^{\alpha} \partial_{x}^{\beta}
  \phi^{(k)} A'(v) \rangle \right|
\\
\leq
M_{4} \sum_{\ell=2}^{\min\{\alpha, 2a\}} \alpha^{\ell - 1} \| X_{3} \partial^{\alpha+1-\ell}_{y}
\partial_{x}^{\beta} \phi^{(k)} A'(v) \|_{0}
\\
\cdot
\alpha
\| X_{3} \partial_{y}^{\alpha-1} \partial_{x}^{\beta}
  \phi^{(k)} A'(v) \|_{0} .
\end{multline*}
The inductive hypothesis can be applied to both factor in the sum
above and we proceed as for \eqref{eq:commX3}.

Choosing $ \nu $ sufficiently large achieves the proof of
\eqref{eq:XA'} and hence of \eqref{eq:XA}. 

Via the Sobolev immersion theorem as well as Lemma \ref{lemma:global},
we get a pointwise estimate of the
same type, from which setting $ x=0 $ we deduce the conclusion of the
lemma. 

\end{proof}
Let us now prove \eqref{eq:Fav1}.
\begin{lemma}
\label{lemma:fourier}
For any $ \epsilon > 0 $ there exist positive constants $ K_{\epsilon}
$, $ M $, such that
\begin{equation}
\label{eq:fourier}
\left | \mathscr{F}(\partial_{x}^{\beta} A(v))(0, \eta) \right| \leq K_{\epsilon} e^{-
 M \left(\frac{ |\eta|}{\epsilon}\right)^{\frac{1}{s_{0}}}} ,
\end{equation}
for $ \beta = 0, 1 $.
\end{lemma}
\begin{proof}
Set $ \beta = 0 $. The argument is exactly the same when $ \beta = 1
$. We have,
$$ 
| D^{\alpha}_{y} A(v)(0, y) | \leq C_{\epsilon} \langle y
\rangle^{-a} \epsilon^{\alpha} \alpha!^{s_{0}},
$$
by Lemma \ref{lemma:yd}. Then we obtain
\begin{equation}
  \label{eq:FAv}
|\mathscr{F}(A(v))(0, \eta) | \leq \frac{1}{|\eta|^{\alpha}} C_{\epsilon} \epsilon^{\alpha}
\alpha!^{s_{0}} \int_{\R} \langle y \rangle^{-a} dy \leq C_{1
  \epsilon} \frac{1}{|\eta|^{\alpha}} \epsilon^{\alpha}
\alpha!^{s_{0}} .
\end{equation}
Then
$$ 
| \mathscr{F}(A(v))(0, \eta) |^{\frac{1}{s_{0}}}
\frac{\left(\left(\frac{|\eta|}{2\epsilon}\right)^{\frac{1}{s_{0}}}
  \right)^{\alpha}}{\alpha!} \leq C_{1 \epsilon}^{\frac{1}{s_{0}}}
\left( \frac{1}{2^{\frac{1}{s_{0}}}} \right)^{\alpha},
$$
whence, summing on $ \alpha $, we deduce the assertion of the lemma. 
\end{proof}
Next we prove an estimate from below of $ \left |
  \mathscr{F}(\partial_{x}^{\beta} A(v))(0, \eta) \right| $, where $
\beta $ takes either the value 0 or 1, depending on the ground state,
$ \phi_{0} $, of the operator $ Q $.
\begin{lemma}
\label{lemma:below}
There exist positive constants $ \lambda $, $ C_{\lambda} $ and a real
constant $ \lambda' $ such that for $ \eta \geq (6 R_{0})^{s_{0}} $
and $ R_{0} $ large enough, we have
\begin{equation}
\label{eq:below}
\left | \mathscr{F}(\partial_{x}^{\beta} A(v))(0, \eta) \right| \geq
C_{\lambda} \eta^{\lambda'} e^{- \lambda |\eta|^{\frac{1}{s_{0}}}} .
\end{equation}
\end{lemma}
\begin{proof}
We have
$$ 
\partial_{x}^{\beta} A(v) (0, y) = \int_{0}^{+\infty} e^{i y
  \rho^{s_{0}}} \rho^{r + \beta \frac{s_{0}}{q}} 
  \partial_{x}^{\beta} v (0, \rho) d\rho ,
$$
where $ v $ is given by \eqref{eq:v}. By Lemma \ref{lemma:psij} we see
that the $ \rho $-support of $ v $ is contained in $ [3 R_{0}, +\infty
[$, hence changing variables according to $ \eta = \rho^{s_{0}} $ we
obtain
$$ 
\partial_{x}^{\beta} A(v) (0, y) = \frac{1}{s_{0}} \int_{-\infty}^{+\infty} e^{i y \eta}
\eta^{\frac{r}{s_{0}} + \frac{\beta}{q}} 
  \partial_{x}^{\beta} v (0, \eta^{\frac{1}{s_{0}}})
  \eta^{\frac{1}{s_{0}} -1}  d\eta .
$$
By \eqref{eq:convv} we have that
\begin{multline}
  \label{eq:Ft}
\mathscr{F}(\partial_{x}^{\beta} A(v))(0, \eta) = \frac{2\pi}{s_{0}}
\eta^{\frac{r+1}{s_{0}} +  \frac{\beta}{q} -1} 
\partial_{x}^{\beta} v (0, \eta^{\frac{1}{s_{0}}})
\\
=
\frac{2\pi}{s_{0}}
\eta^{\frac{r+1}{s_{0}} +  \frac{\beta}{q} -1}
\sum_{j \geq 0} \psi_{j}(\eta^{\frac{1}{s_{0}}} ) \partial_{x}^{\beta}
u_{j}(0, \eta^{\frac{1}{s_{0}}}) .
\end{multline}
Moreover by \eqref{eq:u0-}, \eqref{eq:mitlda0} we have
\begin{multline*}
\left | \sum_{j \geq 0} \psi_{j}(\eta^{\frac{1}{s_{0}}} ) \partial_{x}^{\beta}
  u_{j}(0, \eta^{\frac{1}{s_{0}}}) \right|
\geq
\psi_{0}(\eta^{\frac{1}{s_{0}}} ) \left| e^{\mu_{0 i^{*}} \eta^{\frac{1}{s_{0}}}}
  \partial_{x}^{\beta}\phi_{0}(0) \right|
\\
-
\sum_{j > 0} \psi_{j}(\eta^{\frac{1}{s_{0}}} ) \left |\partial_{x}^{\beta}
  u_{j}(0, \eta^{\frac{1}{s_{0}}}) \right|
\\
\geq
\psi_{0}(\eta^{\frac{1}{s_{0}}} )  e^{\tilde{\mu}_{0} \eta^{\frac{1}{s_{0}}}}
\left|  \partial_{x}^{\beta}\phi_{0}(0) \right|
-
\sum_{j > 0} \psi_{j}(\eta^{\frac{1}{s_{0}}} ) e^{\tilde{\mu}_{0}
  \eta^{\frac{1}{s_{0}}}}  \left |e^{- \tilde{\mu}_{0} \eta^{\frac{1}{s_{0}}}} \partial_{x}^{\beta}
  u_{j}(0, \eta^{\frac{1}{s_{0}}}) \right| .
\end{multline*}
In view of the estimates \eqref{eq:pest} and Lemma \ref{lemma:psij},
we examine the sum on the right hand side of the above inequality.
\begin{multline*}
\sum_{j > 0} \psi_{j}(\eta^{\frac{1}{s_{0}}} ) e^{\tilde{\mu}_{0}
  \eta^{\frac{1}{s_{0}}}}  \left |e^{- \tilde{\mu}_{0} \eta^{\frac{1}{s_{0}}}} \partial_{x}^{\beta}
  u_{j}(0, \eta^{\frac{1}{s_{0}}}) \right|
\\
\leq
\sum_{j > 0} \psi_{j}(\eta^{\frac{1}{s_{0}}} ) e^{\tilde{\mu}_{0}
  \eta^{\frac{1}{s_{0}}}} K_{u}^{1+j+\beta}
\eta^{-\frac{(j-\delta)\kappa}{s_{0}}} \left( \frac{j}{s_{0}} + \beta
  \frac{q-1}{q}\right)!
\\
\leq
C_{\psi}  e^{\tilde{\mu}_{0} \eta^{\frac{1}{s_{0}}}}
\sum_{j > 0}   \tilde{K}_{u}^{1+j}
\left(3 R_{0}(j+1)\right)^{-(j-\delta)\kappa} (j+1)^{\frac{j}{s_{0}}}
\\
\leq
C_{\psi}  e^{\tilde{\mu}_{0} \eta^{\frac{1}{s_{0}}}}
(3 R_{0})^{-(1-\delta)\kappa} \sum_{j > 0} \tilde{K}_{u}^{1+j}
(j+1)^{j(\frac{1}{s_{0}} -\kappa) + \delta \kappa}
\\
=
C_{1} (3 R_{0})^{-(1-\delta)\kappa} e^{\tilde{\mu}_{0}
  \eta^{\frac{1}{s_{0}}}} .
\end{multline*}
In the last line above we used the fact that $ \delta < 1 $, $ \kappa >
\frac{1}{s_{0}} $ (see the statement of Theorem \ref{th:pest-t}).

As a consequence, taking $ \eta^{\frac{1}{s_{0}}} \geq 6 R_{0} $, we
have the following lower bound
\begin{multline}
\label{eq:lowerbd}
\left | \sum_{j \geq 0} \psi_{j}(\eta^{\frac{1}{s_{0}}} ) \partial_{x}^{\beta}
  u_{j}(0, \eta^{\frac{1}{s_{0}}}) \right|
\\
\geq   e^{\tilde{\mu}_{0} \eta^{\frac{1}{s_{0}}}} \left[  \left|
    \partial_{x}^{\beta}\phi_{0}(0) \right| -
C_{1} (3 R_{0})^{-(1-\delta)\kappa} \right].
\end{multline}
Here $ \phi_{0} $ denotes the ground state of the operator $ Q $
defined in \eqref{eq:Q}. We point out that either $ \phi_{0}(0) $ or $
\partial_{x} \phi_{0}(0) $ are non zero. This is easily seen by
remarking that if both $ \phi_{0}(0) $ and $\partial_{x} \phi_{0}(0) $
are zero then all derivatives of $ \phi_{0} $ vanish at the origin,
which is false due to the analyticity of $ \phi_{0} $.

Choosing $ R_{0} $ large enough we deduce that there is a $ \beta \in
\{0, 1\} $ and there is a constant $ C_{2} > 0 $, such that, for
$ \eta^{\frac{1}{s_{0}}} \geq 6 R_{0} $,
\begin{equation}
\label{eq:lowerbd-1}
\left | \sum_{j \geq 0} \psi_{j}(\eta^{\frac{1}{s_{0}}} ) \partial_{x}^{\beta}
  u_{j}(0, \eta^{\frac{1}{s_{0}}}) \right|
\geq
C_{2} e^{\tilde{\mu}_{0} \eta^{\frac{1}{s_{0}}}} .
\end{equation}
Then, plugging this into \eqref{eq:Ft} we get
$$ 
\left| \mathscr{F}(\partial_{x}^{\beta} A(v))(0, \eta) \right| \geq
C_{3} \eta^{\frac{\re r+1}{s_{0}} + \frac{\beta}{q} -1} e^{\tilde{\mu}_{0} \eta^{\frac{1}{s_{0}}}} ,
$$
for $ \eta^{\frac{1}{s_{0}}} \geq 6 R_{0} $, and a suitable positive
constant $ C_{3} $.

This proves the lemma.
\end{proof}
Now Lemma \ref{lemma:fourier} as well as Lemma \ref{lemma:below} yield
a contradiction, thus proving Theorem \ref{th:1}.

\setcounter{section}{0}
\renewcommand\thesection{\Alph{section}}
\section{Appendix: Solution of an ODE}

\setcounter{equation}{0}
\setcounter{theorem}{0}
\setcounter{proposition}{0}
\setcounter{lemma}{0}
\setcounter{corollary}{0}
\setcounter{definition}{0}
\renewcommand{\thetheorem}{\thesection.\arabic{theorem}}
\renewcommand{\theproposition}{\thesection.\arabic{proposition}}
\renewcommand{\thelemma}{\thesection.\arabic{lemma}}
\renewcommand{\thedefinition}{\thesection.\arabic{definition}}
\renewcommand{\thecorollary}{\thesection.\arabic{corollary}}
\renewcommand{\theequation}{\thesection.\arabic{equation}}
\renewcommand{\theremark}{\thesection.\arabic{remark}}

Let $ \mu > 0 $ be a positive number. We want to find the solution of
the ordinary differential equation
\begin{equation}
\label{eq:ode}
\left( D_{\rho}^{2a} + \mu \right) u = f,
\end{equation}
where $ f \in C^{\infty}(\R) $, $ \supp f \subset \R^{+} $, which is
rapidly decreasing at $ +\infty $. 

Let us denote by $ \mu_{i} $, $ i = 1, \ldots, 2a $, the $ 2a $-roots
of $ (-1)^{a+1} \mu $. We observe that $ \re \mu_{i} \neq 0 $ for
every $ i $. In fact if $ a $ is an even integer, in order to have $
\re \mu_{k+1} = 0 $, for some $ k = 0, 1, \ldots, 2a-1 $, we should
have
$$ 
\frac{\pi}{2a} + k \frac{\pi}{a} = \ell \frac{\pi}{2},
$$
for some odd integer $ \ell $. This would imply $ 2k + 1 = a \ell $,
which is impossible. Analogously, assume $ a $ is an odd integer. Then
if $\re \mu_{k+1} = 0$ we must have
$$ 
k \frac{\pi}{a} = \ell \frac{\pi}{2}.
$$
This would imply $ 2 k = a \ell $, which is also impossible, since $ a
\ell$ is odd.

For any $ j = 1, \ldots, 2a $, define
\begin{multline}
\label{eq:If}
I_{j}(f)(\rho)
\\
= - \sgn\left( \re \mu_{j}\right) \int_{\R} e^{\mu_{j}
  (\rho - \sigma)} H\left(-\sgn\left(\re \mu_{j}\right) (\rho -
  \sigma)\right) f(\sigma) d\sigma ,
\end{multline}
where $ H $ denotes the Heaviside function.

Since
$$ 
\frac{1}{\sigma^{2a} + (-1)^{a} \mu} = \prod_{j=1}^{2a}
\frac{1}{\sigma - \mu_{j}},
$$
define the positive numbers $ A_{j} $ by the relation
\begin{equation}
\label{eq:Aj}
\sum_{j=1}^{2a} A_{j} \frac{1}{\sigma - \mu_{j}} = \prod_{j=1}^{2a}
\frac{1}{\sigma - \mu_{j}}.
\end{equation}
Multiplying both sides of \eqref{eq:Aj} by $ \sigma - \mu_{\ell} $, $
\ell \in \{1, \ldots, 2a\} $, we get
$$ 
A_{\ell} + \sum_{\substack{j = 1\\j\neq\ell}}^{2a} A_{j} \frac{\sigma
  - \mu_{\ell}}{\sigma - \mu_{j}} =
\prod_{\substack{j=1\\j\neq\ell}}^{2a} \frac{1}{\sigma - \mu_{j}} .
$$
Computing the above identity for $ \sigma = \mu_{\ell} $ we obtain an
expression for the $ A_{j} $:
\begin{equation}
  \label{eq:Al}
A_{\ell} = \prod_{\substack{j=1\\j\neq\ell}}^{2a} \frac{1}{\mu_{\ell}
  - \mu_{j}} .
\end{equation}
Define now
\begin{equation}
\label{eq:u}
u(\rho) = \sum_{j=1}^{2a} A_{j} I_{j}(f)(\rho) .
\end{equation}
We want to show that $ u $ is the desired solution.

First observe that
$$ 
\partial_{\rho} I_{j}(f) = \mu_{j} I_{j}(f) + f .
$$
Hence we deduce
\begin{align}
\label{eq:dku}
\partial_{\rho} u & = \sum_{j=1}^{2a} A_{j} \left(f +\mu_{j} I_{j}(f)
                    \right) \\
\partial_{\rho}^{2} u & = \sum_{j=1}^{2a} A_{j} \left( f' + \mu_{j} f
                        + \mu_{j}^{2} I_{j}(f) \right) \notag \\  
  \vdots         &\hphantom{=}  \hspace{2.cm} \vdots  \notag \\
\partial_{\rho}^{2a} u & = \sum_{j=1}^{2a} A_{j} \left( f^{(2a-1)} +
                         \mu_{j} f^{(2a-2)} + \cdots + \mu_{j}^{2a-1}
                         f + \mu_{j}^{2a} I_{j}(f) \right).
                         \notag
\end{align}
Now, because of \eqref{eq:Aj}, we have
\begin{equation}
  \label{eq:Aj:2}
\sum_{i=1}^{2a} A_{i} \prod_{\substack{j=1\\j\neq i}}^{2a} (\sigma -
\mu_{j}) = 1 .
\end{equation}
Moreover the polynomial multiplying a single $ A_{i} $ above is
written as
\begin{multline}
  \label{eq:Aj:3}
\prod_{\substack{j=1\\j\neq i}}^{2a} (\sigma - \mu_{j}) =
\sigma^{2a-1} - \sigma^{2a-2} \sum_{\substack{j=1\\j\neq i}}^{2a} \mu_{j}
+\sigma^{2a-3} \sum_{\substack{j_{1} < j_{2} =1 \\ j_{1}, j_{2} \neq
    i}}^{2a} \mu_{j_{1}} \mu_{j_{2}} \\
+ \cdots + (-1)^{k} \sum_{\substack{j_{1} < \cdots < j_{k} =1 \\
    j_{1}, \ldots, j_{k} \neq i }}^{2a} \mu_{j_{1}} \ldots \mu_{j_{k}}
\sigma^{2a-k-1} + \cdots + (-1)^{2a-1} \prod_{\substack{j=1\\j \neq
    i}}\mu_{j} .
\end{multline}
Defining as $ s_{k}^{(i)}= s_{k}^{(i)}(\mu_{1}, \ldots, \mu_{2a}) $
the symmetric function of degree $ k $ of the $ 2a-1 $ arguments $
\mu_{1}, \ldots, \mu_{i-1}, \mu_{i+1}, \ldots, \mu_{2a} $, with $
s_{0}^{(i)} = 1 $, the above identity can be rewritten as
$$ 
\prod_{\substack{j=1\\j\neq i}}^{2a} (\sigma - \mu_{j}) =
\sigma^{2a-1} + \sum_{k=1}^{2a-1} (-1)^{k} s_{k}^{(i)} \sigma^{2a-1-k} .
$$
Let us show inductively that for $ 1 \leq k \leq 2a $ we have
\begin{equation}
\label{eq:ski}
s_{k}^{(i)} = s_{k} - \mu_{i} s_{k-1} + \mu_{i}^{2} s_{k-2} + \cdots +
(-1)^{k} \mu_{i}^{k},
\end{equation}
where $ s_{k} $ denotes the symmetric function of $ k $ of the $ 2a $
arguments $ \mu_{1}, \ldots,$  $ \mu_{2a} $, and we make the convention
that $ s_{0} = 1 $. 
For $ k=1 $ it is obviously true. Consider
\begin{align*}
s_{k}^{(i)} &= \sum_{\substack{j_{1} < \cdots <j_{k} =1 \\ j_{1},
    \cdots, j_{k} \neq i}}^{2a} \mu_{j_{1}} \cdots \mu_{j_{k}}
\\
&= \sum_{\substack{j_{1} < \cdots <j_{k} =1 \\ j_{2},
    \cdots, j_{k} \neq i}}^{2a} \mu_{j_{1}} \cdots \mu_{j_{k}} -
\mu_{i} \sum_{j_{2} < \cdots <j_{k} = i+1}^{2a} \mu_{j_{2}} \cdots
\mu_{j_{k}}
\end{align*}
Iterating the argument we get
\begin{align*}
s_{k}^{(i)} &= \sum_{\substack{j_{1} < \cdots <j_{k} =1 \\ j_{3},
    \cdots, j_{k} \neq i}}^{2a} \mu_{j_{1}} \cdots \mu_{j_{k}} -
\mu_{i} \sum_{j_{2} < \cdots <j_{k} = i+1}^{2a} \mu_{j_{2}} \cdots
\mu_{j_{k}}
  \\
  &\phantom{=}\ 
- \mu_{i} \sum_{j_{1} < i < j_{3} < \cdots < j_{k}}
    \mu_{j_{1}} \mu_{j_{3}} \cdots \mu_{j_{k}}
  \\
&= \cdots =
\sum_{j_{1} < \cdots <j_{k}} \mu_{j_{1}} \cdots \mu_{j_{k}} - \mu_{i}
\left( \sum_{i < j_{2} < \cdots <j_{k}} \mu_{j_{2}} \cdots
  \mu_{j_{k}} \right.
  \\
  &\phantom{=}\ 
\left. +
  \sum_{j_{1} < i < j_{3} < \cdots <j_{k}} \mu_{j_{1}} \mu_{j_{3}} \cdots
  \mu_{j_{k}} + \cdots +
\sum_{j_{1} < \cdots <j_{k-1} < i} \mu_{j_{1}} \cdots \mu_{j_{k-1}}
\right)
\\
&=
s_{k} - \mu_{i} s_{k-1}^{(i)}.
\end{align*}
Applying the inductive hypothesis
\begin{align*}
s_{k}^{(i)} &= s_{k} -\mu_{i} \left( s_{k-1} - \mu_{i} s_{k-2} + \cdots
              + (-1)^{k} \mu_{i}^{k-1}\right) \\
  &= \sum_{\ell=0}^{k} (-1)^{\ell} \mu_{i}^{\ell} s_{k-\ell},
\end{align*}
which is the desired conclusion.

Going back to \eqref{eq:Aj:2}, \eqref{eq:Aj:3} we see that
\eqref{eq:Aj:3} can be written as
$$ 
\prod_{\substack{j=1\\j\neq i}}^{2a} (\sigma - \mu_{j}) =
\sum_{\ell=0}^{2a-1} (-1)^{\ell} \sigma^{2a-1-\ell} s_{\ell}^{(i)} ,
$$
so that, identifying the polynomials on both sides of \eqref{eq:Aj:2},
we obtain
$$ 
\sum_{i=1}^{2a} A_{i} s_{k}^{(i)} = 0, \qquad \text{for } k = 0,
\ldots, 2a-2.
$$
When $ k=0 $ we immediately get that $ \sum_{i=1}^{2a}A_{i} = 0
$. Assume $ k=1 $. Since $ s_{1}^{(i)} = s_{1} - \mu_{i} $ we have
that $ \sum_{i=1}^{2a} A_{i} s_{1}^{(i)} = s_{1} \sum_{i=1}^{2a} A_{i}
-  \sum_{i=1}^{2a} A_{i} \mu_{i} = 0$, giving that $ \sum_{i=1}^{2a}
A_{i} \mu_{i} = 0 $.

Iterating this kind of argument we obtain that
\begin{equation}
  \label{eq:Amuk}
\sum_{i=1}^{2a} A_{i} \mu_{i}^{k} = 0, \qquad \text{for } k = 0,
\ldots, 2a-2.
\end{equation}
Finally, again from \eqref{eq:Aj:2}, \eqref{eq:Aj:3}, the identity
$$ 
- \sum_{i=1}^{2a} A_{i} \mu_{1} \cdots \mu_{i-1} \mu_{i+1} \cdots
\mu_{2a} = 1
$$
implies that
\begin{multline}
  \label{eq:Amu2a-1}
1 = - \sum_{i=1}^{2a} A_{i} \frac{1}{\mu_{i}} \mu_{1} \cdots \mu_{2a}
= - \sum_{i=1}^{2a} A_{i} \frac{1}{\mu_{i}} (-1)^{a} \mu
\\
= \sum_{i=1}^{2a} A_{i} \frac{1}{\mu_{i}} (-1)^{a+1} \mu
= \sum_{i=1}^{2a} A_{i} \frac{1}{\mu_{i}} \mu_{i}^{2a}
= \sum_{i=1}^{2a} A_{i} \mu_{i}^{2a-1} .
\end{multline}
Thus $ u $, as defined in \eqref{eq:u}, is a solution of
\eqref{eq:ode}.

\section{Appendix: Estimate of the Ground Level Eigenfunction}

\setcounter{equation}{0}
\setcounter{theorem}{0}
\setcounter{proposition}{0}
\setcounter{lemma}{0}
\setcounter{corollary}{0}
\setcounter{definition}{0}

This section contains the proof of the estimates of the function $
\phi_{0} \in \ker Q_{\lambda_{0}} $, where $ Q_{\lambda_{0}} =
D_{x}^{2} + x^{2(q-1)} - \lambda_{0} $. The method for
obtaining such estimates has been introduced by M\'etivier in
\cite{metivier-duke-80} in a homogeneous, i.e. quadratic, case. 

Let us start by defining
\begin{equation}
\label{eq:X+}
\left \{
\begin{array}{ccc}
  X_{+} & = & \partial_{x}  \\[7pt]
  X_{-} & = & x 
\end{array}
\right .
\end{equation}
For $ k \in \N $ we denote by $ I $ the multiindex $ I = (i_{1},
\ldots, i_{k}) $, where $ i_{j} \in \{ \pm \} $ for $ j = 1, \ldots, k
$. We also write $ k = | I | $. Define
\begin{equation}
\label{eq:XI}
X_{I} = X_{i_{1}} \cdots X_{i_{k}}.
\end{equation}
Set
$$ 
I_{+} = (i_{1}^{+}, \ldots, i_{k}^{+}) ,
$$
where
\begin{equation}
\label{eq:i+}
i_{\nu}^{+} =
\begin{cases}
  + , & \text{ if } i_{\nu} = + \\[5pt]
  0 , & \text{ if } i_{\nu} = -
\end{cases}
,
\end{equation}
and analogously for $ I_{-} $ and $ i_{\nu}^{-} $. Define $ |I_{+}| $
as the number of the non-zero components of $ I_{+} $ and similarly
for $ |I_{-}| $.

We finally set
\begin{equation}
\label{eq:bracketsI}
\langle I \rangle = |I_{+}| + \frac{| I_{-}|}{q-1} .
\end{equation}
We are going to need the spaces $ H^{k}_{q}(\R) $, which for suitable
$ k $ are natural domains of the operator $ Q_{\lambda_{0}} $ in $
L^{2}(\R) $:
\begin{equation}
\label{eq:Hk}
H^{k}_{q}(\R) = \{ u \in L^{2}(\R) \ | \ X_{I}u \in L^{2}(\R), \text{
  for every } I, \langle I \rangle \leq k \},
\end{equation}
for $ k \in \N \cup \{0\} $.

We equip $ H^{k}_{q}(\R) $ with the norm
\begin{equation}
  \label{eq:normq}
  \| u \|_{k} = \max_{0 \leq \ell \leq k} | u |_{\ell} ,
\end{equation}
where
$$ 
|u|_{\ell} = \max_{\langle I \rangle = \ell} \| X_{I} u \|_{L^{2}(\R)} .
$$
For the sake of simplicity we write $ \| u \|_{0} $ for the $
L^{2}(\R) $ norm of $ u $.

Due to Gru\v sin, \cite{grushin70}, or rather Lemma \ref{lemma:apr}
below for the anisotropic case, we know that the following a
priori estimate is satisfied
\begin{equation}
\label{eq:apriori}
\| u \|_{2} \leq C_{0} \left( \| Q_{\lambda_{0}} u \|_{0} + \| u
  \|_{0} \right),
\end{equation}
for a suitable positive constant $ C_{0} $.

We want to prove the following
\begin{proposition}
\label{prop:1st-est}
Assume that $ u \in \ker Q_{\lambda_{0}} $.  
Then there exist positive constants $ C $, $ R $, depending only on the
operator $ Q_{\lambda_{0}} $, such that, for every multiindex $ I $,
we have the inequality
\begin{equation}
\label{eq:XIu}
\| X_{I} u \| \leq C R^{\langle I \rangle} \|u \|_{0} (\langle I
\rangle !)^{\frac{q-1}{q}},
\end{equation}
where, for $ x > 0 $, $ x! $ means $ \Gamma(x+1) $ and $ 0! = 1 $.
\end{proposition}
\begin{corollary}
\label{cor:xd}
We have for any $ u \in \ker Q_{\lambda_{0}} $
\begin{equation}
\label{eq:xd}
\| x^{\beta} \partial_{x}^{\alpha} u \|_{0} \leq C^{\alpha+\beta+1}
\left( \alpha \frac{q-1}{q} + \frac{\beta}{q}\right)!
\end{equation}
\end{corollary}
Before proving the proposition we state a couple of lemmas that are
used in its proof.
\begin{lemma}
  \label{lemma:comm}
  Let $ I $ be a multliindex. Then
\begin{equation}
\label{eq:comm}
\| [ Q_{\lambda_{0}} , X_{I} ] u \|_{0} \leq C |I| \| u \|_{2 +
  \langle I \rangle - \frac{q}{q-1}},
\end{equation}
for any $ u $ in the $ L^{2} $ domain of the operator on the left hand
side.  
\end{lemma}
\begin{proof}[Proof of Lemma \ref{lemma:comm}]
The assertion is proved by remarking that
$$ 
[ Q_{\lambda_{0}} , X_{I} ] = \sum_{I_{1}} b_{I_{1}} X_{I_{1}} ,
$$
where $ b_{I_{1}} \in \C $ and are bounded by a quantity depending
only on the problem and $ \langle I_{1} \rangle = \langle I \rangle +
2 - \frac{q}{q-1} $. Moreover the number of terms in the summation
above is bounded by $ \kappa | I | $, $ \kappa $ denoting a positive
constant depending only on the problem data.

Inequality \eqref{eq:comm} then immediately follows.
\end{proof}
\begin{lemma}
\label{lemma:dec}
Let $ I $ a multiindex such that $ \langle I \rangle > 2 $. Then we
may decompose $ X_{I} $ as
$$ 
X_{I} = X_{I''} X_{I'} + A ,
$$
where $ \langle I'' \rangle = 2 $, $ \langle I' \rangle = \langle I
\rangle - 2 $ and $ A $ is a finite sum of the form
$$ 
A = \sum_{J} c_{J} X_{J} ,
$$
with $ \langle J \rangle = \langle I \rangle - \frac{q}{q-1} $. Here
both the coefficients $ c_{J} $ and the number of summands are bounded
by an absolute constant.
\end{lemma}
The proof of the lemma is straightforward and we skip it.
\begin{proof}[Proof of Proposition \ref{prop:1st-est}]
Instead of proving \eqref{eq:XIu} directly we are going to show that
the following inequality holds:
\begin{equation}
\label{eq:XIu'}
\| X_{I} u \| \leq C_{1} R_{1}^{\langle I \rangle} \|u \|_{0} \langle I
\rangle^{\langle I \rangle \frac{q-1}{q}},
\end{equation}
for certain positive constants $ C_{1} $, $ R_{1} $, for every
multiindex $ I $. It is then obvious that \eqref{eq:XIu'} implies
\eqref{eq:XIu} slightly modifying the constants, because of the
Stirling formula for the Euler Gamma function.

We argue by induction on $ k $ where $ \langle I \rangle = \frac{k}{q-1} $.

First of all we observe that when $ \langle I \rangle \leq 2 $ we have
by \eqref{eq:apriori} that
$$ 
\| X_{I} u \| \leq C_{0}\left(\|Q_{\lambda_{0}} u \|_{0} + \| u \|_{0}
\right) = C_{0} \| u \|_{0}.
$$
Assume now that the assertion holds for any $ I $ with $ \langle I
\rangle = \frac{\ell}{q-1} $, $ \ell = 0, 1, \ldots, k $, $ k > 2(q-1)
$. We want to show that the assertion is true for $ I $ with $ \langle
I \rangle = \frac{k+1}{q-1} $. 

Let $ I $ be a multiindex with $ \langle I \rangle = \frac{k+1}{q-1}
$. Using Lemma \ref{lemma:dec} we write $ X_{I} = X_{I'} X_{J} + A $,
with $ \langle I' \rangle = 2 $, $ \langle J \rangle = \langle I
\rangle - 2 $ and $ A $ of the form specified in Lemma
\ref{lemma:dec}.

By \eqref{eq:apriori} we have
\begin{eqnarray*}
\| X_{I} u \|_{0} & \leq & \|X_{I'} X_{J} u \|_{0} + \| Au \|_{0} 
  \\[5pt]
   & \leq & C_{0} \left ( \| Q_{\lambda_{0}} X_{J} u \|_{0} + \| X_{J} u
          \|_{0} + \| A u \|_{0} \right).
\end{eqnarray*}
Since $ u \in \ker Q_{\lambda_{0}} $,
$$
Q_{\lambda_{0}} X_{J} u = X_{J} Q_{\lambda_{0}} u + [ Q_{\lambda_{0}}
, X_{J} ] u =  [ Q_{\lambda_{0}} , X_{J} ] u
$$
Hence
\begin{equation}
  \label{eq:apriori-comm}
\| X_{I} u \|_{0} \leq C_{0} \left ( \| [ Q_{\lambda_{0}} , X_{J} ] u
  \|_{0} + \| X_{J} u \|_{0} + \| A u \|_{0} \right).
\end{equation}
Since $ \langle J \rangle = \langle I \rangle -2 = \frac{k+1}{q-1} - 2
\leq \frac{k}{q-1}$ and $ A $ is a sum of terms involving $ X_{J'} $,
with $ \langle J' \rangle = \langle I \rangle - \frac{q}{q-1} =
\frac{k- (q-1)}{q-1} = \frac{k}{q-1} - 1 < \frac{k}{q-1} $, we see that both terms $ \|
X_{J}u \| $ and $ \| A u \| $ satisfy the inductive hypothesis:
\begin{multline*}
  \| X_{J} u \| \leq C_{1} R_{1}^{\langle J \rangle} \| u \|_{0} \langle
  J \rangle^{\langle J \rangle \frac{q-1}{q}}
\\
= C_{1}
R_{1}^{\frac{k+1}{q-1} -2} \| u \|_{0} \left(\frac{k+1}{q-1}
  -2\right)^{\left(\frac{k+1}{q-1} -2\right) \frac{q-1}{q}}
\leq
(C_{1} R_{1}^{-2}) R_{1}^{\langle I \rangle} \| u \|_{0} \langle I
\rangle^{\langle I \rangle \frac{q-1}{q}} .
\end{multline*}
Furthermore writing $ A = \sum_{J'} c_{J'} X_{J'} $, with $ \langle J'
\rangle = \frac{k}{q-1} -1 $ and where both the number of summands and
the constants $ c_{J'} $ are bounded by a universal constant, say $ M
$, we have 
\begin{multline*}
\| Au \|_{0} \leq \sum_{J'} | c_{J'} | \| X_{J'} u \|_{0}
\leq \sum_{J'} | c_{J'} | C_{1} R_{1}^{\langle J' \rangle} \| u \|_{0}
\langle J' \rangle^{\langle J' \rangle \frac{q-1}{q}}
\\
\leq
( C_{1} M^{2} R_{1}^{- \frac{q}{q-1}} ) R_{1}^{\langle I \rangle} \|
u \|_{0}  \langle I \rangle^{\langle I \rangle \frac{q-1}{q}} .
\end{multline*}
Note that there is always a gain of a negative power of $
R_{1} $ in the above terms. 

Consider now the norm with the commutator in the right hand side of
\eqref{eq:apriori-comm}. By Lemma \ref{lemma:comm} we have
\begin{eqnarray}
  \label{eq:comm-est}
\| [ Q_{\lambda_{0}} , X_{J} ] u \|_{0} &\leq & C |J| \| u \|_{2 +
\langle J \rangle - \frac{q}{q-1}}  \notag
\\
 & \leq & C |J| \max_{0 \leq \ell \leq
2+\langle J \rangle - \frac{q}{q-1}} \max_{\langle I' \rangle = \ell}
\|X_{I'} u \|_{0} .
\end{eqnarray}
Observe that
$$ 
2 + \langle J \rangle - \frac{q}{q-1} = \langle I \rangle -
\frac{q}{q-1} = \frac{k}{q-1} - 1.
$$
Thus the norms in the right hand side of \eqref{eq:comm-est} satisfy
the inductive hypothesis. We deduce that
\begin{eqnarray*}
\| [ Q_{\lambda_{0}} , X_{J} ] u \|_{0}  & \leq & C k C_{1}
R_{1}^{\frac{k}{q-1} -1} \| u \|_{0} \left( \frac{k}{q-1} -1
\right)^{\left(\frac{k}{q-1} -1\right) \frac{q-1}{q}}
\\
& \leq &
\left( \frac{C C_{1} (q-1)}{R_{1}^{\frac{q}{q-1}}} \right) \| u \|_{0} 
R_{1}^{\frac{k+1}{q-1}} \left(
\frac{k+1}{q-1}\right)^{\frac{k+1-q}{q} +1}
\\
 &  = &
\left( \frac{C C_{1} (q-1)}{R_{1}^{\frac{q}{q-1}}} \right) \| u \|_{0} 
R_{1}^{\frac{k+1}{q-1}} \left(
        \frac{k+1}{q-1}\right)^{\frac{k+1}{q}} .
\end{eqnarray*}
Choosing $ R_{1} $ in such a way that the constant in parentheses
is bounded by $ C_{1} $ achieves the proof of the proposition.
\end{proof}
Taking $ |I_{-}| = \beta $, $ |I_{+}| = \alpha $ and
using the Sobolev embedding theorem (in one dimension,) we
obtain a bound for the derivatives of the functions in $ \ker
Q_{\lambda_{0}} $:
\begin{corollary}
\label{cor:Dphi}
Let $ \phi \in \mathscr{S}(\R) $ be such that $ Q_{\lambda_{0}} \phi =
0$. Then for every $ \alpha, \beta \in \N \cup \{ 0 \} $ there exists a
positive constant $ C_{\phi} $ such that 
\begin{equation}
\label{eq:Dphi}
|x^{\beta} \partial_{x}^{\alpha} \phi(x) | \leq
C_{\phi}^{\alpha+\beta+ 1} \alpha!^{\frac{q-1}{q}} \beta!^{\frac{1}{q}},
\end{equation}
for every $ x \in \R $.
\end{corollary}
\section{Appendix: An Inequality Involving Powers of $ Q $ }

\setcounter{equation}{0}
\setcounter{theorem}{0}
\setcounter{proposition}{0}
\setcounter{lemma}{0}
\setcounter{corollary}{0}
\setcounter{definition}{0}

We prove here the following inequality
\begin{proposition}
\label{prop:ineqQ}
Let $ u \in \mathscr{S}(\R) $, 
and $ Q $ the operator $ D_{x}^{2} + x^{2(q-1)} $, $ \theta \in
\Q^{+} \cup \{ 0 \} $, and $ X_{I} $ an operator of the type defined
in \eqref{eq:XI}.

Set
$$ 
m(\theta, I) = \left[ (2 \theta + \langle I \rangle) \frac{q-1}{q} \right],
$$
where the square brackets denote the integer part. Then there exists a
positive constant $ C $ such that
\begin{equation}
\label{eq:ineqQ}
\| Q^{\theta} X_{I}  u \|_{0} \leq  \sum_{\ell = 0}^{m(\theta, I)}
C^{\ell + 1} \binom{m(\theta, I)}{\ell} p^{\ell} \|
Q^{\frac{1}{2}\left(p - \ell \frac{q}{q-1} \right)} u \|_{0},
\end{equation}
where
\begin{equation}
\label{eq:hk}
p = 2\theta + \langle I \rangle.
\end{equation}
Here we used the fact that $ Q $ is a positive globally elliptic
operator with discrete spectrum and its rational powers are defined
via the spectral mapping theorem, see Helffer \cite{helffer84}.
\end{proposition}
\begin{proof}
The proof is carried out via a number of lemmas.
\begin{lemma}
\label{lemma:A}
Let $ \mu \in \Q^{+} $, $ \mu \geq 1 $, $ I $ a multiindex, then
\begin{equation}
\label{eq:A}
\| Q^{\mu} X_{I} u \|_{0} \leq \| Q^{\mu - 1} X_{I} Q u \|_{0} +
\sum_{I_{1}} c_{I_{1}} \| Q^{\mu - 1} X_{I_{1}} u \|_{0},
\end{equation}
where $ I_{1} $ is a multiindex such that
$$ 
\langle I_{1} \rangle = \langle I \rangle + 2 - \frac{q}{q-1},
$$
the constants $ c_{I_{1}} $ are uniformly bounded by a constant and
the number of summands is bounded by $ M \langle I \rangle $, $ M > 0
$ independent of $ I $. 
\end{lemma}
\begin{proof}[Proof of Lemma \ref{lemma:A}]
We have $ Q^{\mu} X_{I} = Q^{\mu - 1}X_{I} Q + Q^{\mu-1} [ Q,
X_{I}] $. We then remark that
$$ 
[ Q, X_{I}] = \sum_{I_{1}} c_{I_{1}} X_{I_{1}},
$$
where $ \langle I_{1} \rangle = \langle I \rangle + 2 - \frac{q}{q-1}
$ and the bounds in the above statement hold.
\end{proof}
\begin{lemma}
\label{lemma:apr}
Let $ I $ be such that $ \langle I \rangle
= 2$. Then for any $ u \in \mathscr{S}(\R) $ we have the estimate
\begin{equation}
\label{eq:apr}
\| X_{I} u \|_{0} \leq \| Q u \|_{0} + C_{q} \|
x^{q-2} u \|_{0},
\end{equation}
where $ C_{q} $ denotes a positive constant.
\end{lemma}
\begin{proof}
It is a very simple computation. First remark that
$$ 
\| Q u \|_{0}^{2} = \| D_{x}^{2} u \|_{0}^{2} + \| x^{2(q-1)} u
\|_{0}^{2} + \langle \left( x^{2(q-1)} D_{x}^{2} + D_{x}^{2}
  x^{2(q-1)} \right) u, u \rangle.
$$
Now
$$ 
x^{2(q-1)} D_{x}^{2} + D_{x}^{2} x^{2(q-1)} = 2 D_{x} x^{2(q-1)} D_{x}
+ [ D_{x} , [ D_{x} , x^{2(q-1)}] ],
$$
so that we get the identity
$$
\| D_{x}^{2} u \|_{0}^{2} + \| x^{2(q-1)} u \|_{0}^{2} + 2 \|
x^{q-1}D_{x} u \|_{0}^{2}
= \| Q u \|_{0}^{2} + (2q-2)(2q-3) \| x^{q-2} u \|_{0}^{2},
$$
and \eqref{eq:apr} immediately follows remarking that
$$ 
\| X_{I} u \|_{0} \leq \| D_{x}^{2} u \|_{0} + \| x^{2(q-1)} u \|_{0} +  \|
x^{q-1}D_{x} u \|_{0} + C \| x^{q-2} u \|_{0} .
$$
\end{proof}
\begin{lemma}
\label{lemma:B}
Let $ \mu \in \Q^{+} $, $ \mu < 1 $, $ I $ a multiindex, then
\begin{equation}
\label{eq:AC}
\| Q^{\mu} X_{I} u \|_{0} \leq \| Q^{\mu} X_{\tilde{I}} Q u \|_{0} +
\sum_{I_{1}} c_{I_{1}} \| Q^{\mu} X_{I_{1}} u \|_{0} + C \| Q^{\mu + 1
- \frac{q}{2(q-1)}} X_{\tilde{I}} u \|_{0} ,
\end{equation}
where $ \tilde{I} $ is a multiindex such that $ \langle \tilde{I} \rangle
= \langle I \rangle - 2 $, $ \langle I_{1}\rangle = \langle I \rangle -
\frac{q}{q-1} $, 
and for the sum we have bounds analogous to
those in Lemma \ref{lemma:A}.
\end{lemma}
\begin{proof}
We may write
$$ 
X_{I} = X_{\hat{I}} X_{\tilde{I}} + \sum_{I_{1}} c_{I_{1}}' X_{I_{1}} ,
$$
where $ \langle \hat{I} \rangle = 2 $, $ \langle \tilde{I} \rangle =
\langle I \rangle - 2 $, $ \langle I_{1} \rangle = \langle I \rangle -
\frac{q}{q-1}$. Moreover both the constants $ c_{I_{1}}' $ and the
number of the summands are bounded by a universal constant.

Now $ Q^{\mu} X_{\hat{I}} X_{\tilde{I}} = X_{\hat{I}} Q^{\mu}
X_{\tilde{I}} + [ Q^{\mu} , X_{\hat{I}}] X_{\tilde{I}} $.
By lemma \ref{lemma:apr} we have
\begin{eqnarray*} 
\|  X_{\hat{I}} Q^{\mu} X_{\tilde{I}} u \|_{0} & \leq & \| Q^{\mu + 1}
X_{\tilde{I}} u \|_{0} + C_{q} \| x^{q-2} Q^{\mu}  X_{\tilde{I}} u
                       \|_{0}
  \\
& \leq &
\| Q^{\mu + 1} X_{\tilde{I}} u \|_{0} + C_{q} \| Q^{\mu} X_{I_{3}} u
         \|_{0}
\\
&     &  
         + C_{q} \| [ Q^{\mu}, x^{q-2}] X_{\tilde{I}} u \|_{0},
\end{eqnarray*}
where $ \langle I_{3} \rangle = \langle I \rangle - \frac{q}{q-1} $ so
that the second term in the last line above has the same weight as the
terms containing $ X_{I_{1}} $ above. Adapting to the anisotropic case
the calculus of globally elliptic pseudodifferential operators we
get---see Helffer \cite{helffer84}, Theorem 1.11.2 and Proposition
1.6.11---that 
$$ 
\| [ Q^{\mu}, x^{q-2}] X_{\tilde{I}} u \|_{0} \leq C_{1} \| Q^{\mu -
  \frac{1}{q-1}} X_{\tilde{I}} u \|_{0} \leq C_{1} \| Q^{\mu + 1
- \frac{q}{2(q-1)}} X_{\tilde{I}} u \|_{0} ,
$$
$$ 
\| [ Q^{\mu} , X_{\hat{I}}] X_{\tilde{I}} u \|_{0} \leq C_{2} \|  Q^{\mu + 1
- \frac{q}{2(q-1)}} X_{\tilde{I}} u \|_{0}.
$$
The conclusion follows applying Lemma \ref{lemma:A} to the term $\|
Q^{\mu + 1} X_{\tilde{I}} u \|_{0}  $.
\end{proof}
Let us now go back to the proof of inequality \eqref{eq:ineqQ}. We
proceed by induction with respect to $ p = 2 \theta + \langle I
\rangle $. Since $ \theta \in \Q $ and $ \langle I \rangle $ is a
rational number whose denominator is $ q-1 $, we may write $ p $ as a
fraction $ p = p_{1}/ d(\theta, I) $, so that, proceeding by
induction actually means inducing with respect to $ p_{1} $.

If $ p_{1} = 0 $ there is nothing to prove. Assume thus that
\eqref{eq:ineqQ} holds for any $ q_{1}/d(\theta, I) $, $ q_{1} < p_{1} $
and let us show that it is true when $ q_{1} = p_{1} $.

Assume $ \theta < 1 $. Applying Lemma \ref{lemma:B} we obtain
\begin{eqnarray*}
\| Q^{\theta} X_{I} u \|_{0} & \leq & \| Q^{\theta} X_{\tilde{I}} Q u \|_{0} +
\sum_{I_{1}} c_{I_{1}} \| Q^{\theta} X_{I_{1}} u \|_{0}
\\
&   &
+ C \| Q^{\theta + 1 - \frac{q}{2(q-1)}} X_{\tilde{I}} u \|_{0}
\\  
& = &
A_{1} + A_{2} + A_{3},
\end{eqnarray*}
where $ \langle I_{1} \rangle = \langle I \rangle - \frac{q}{q-1} $
and $ \langle \tilde{I} \rangle = \langle I \rangle - 2 $. 
Let us start by examining $ A_{1} $. Write
$$ 
\left( 2 \theta + \langle \tilde{I} \rangle \right) \frac{q-1}{q} =
m(\theta, \tilde{I}) + \tilde{\sigma}, \qquad 0 \leq \tilde{\sigma} < 1,
$$
with $ m(\theta, \tilde{I}) \in \N \cup \{ 0 \} $. Since
$$ 
\left( 2 \theta + \langle \tilde{I} \rangle \right) \frac{q-1}{q} =
\left( 2 \theta + \langle I \rangle \right) \frac{q-1}{q} -1 - \left(
  1 - \frac{2}{q}\right),
$$
we deduce that
$$ 
m(\theta, I) - 2 \leq m(\theta, \tilde{I}) \leq m(\theta, I) -1 .
$$
Hence applying the induction to $ A_{1} $ we obtain
\begin{equation}
\label{eq:A1}
\| Q^{\theta} X_{\tilde{I}} (Q u) \|_{0} \leq \sum_{\ell=0}^{m(\theta,
  \tilde{I})} C^{\ell+1} \binom{m(\theta, \tilde{I})}{\ell}
\tilde{p}^{\ell} \|(Q^{\frac{1}{2}})^{\tilde{p} + 2 - \ell
  \frac{q}{q-1}} u \|_{0} ,
\end{equation}
where $ \tilde{p} = 2 \theta + \langle \tilde{I} \rangle $. 

Since
$$ 
\tilde{p} + 2 = 2 \theta + \langle \tilde{I} \rangle + 2 = 2 \theta +
\langle I \rangle = p ,
$$
we may write
\begin{eqnarray}
\label{eq:A1fin}
\| Q^{\theta} X_{\tilde{I}} (Q u) \|_{0} & \leq &
\sum_{\ell=0}^{m(\theta, \tilde{I})} C^{\ell+1} \binom{m(\theta, \tilde{I})}{\ell}
p^{\ell} \|(Q^{\frac{1}{2}})^{p - \ell \frac{q}{q-1}} u \|_{0}                                                
\notag
\\
&  \leq &
\sum_{\ell=0}^{m(\theta,I) -1} C^{\ell+1} \binom{m(\theta,I) -1}{\ell}
p^{\ell} \|(Q^{\frac{1}{2}})^{p - \ell \frac{q}{q-1}} u \|_{0}
\end{eqnarray}
Next let us examine $ A_{2} $. Applying the inductive hypothesis we
may write
\begin{multline}
\label{eq:A2}
\sum_{I_{1}} c_{I_{1}} \| Q^{\theta} X_{I_{1}} u \|_{0}
\\
\leq 
\sum_{I_{1}} C_{1} \sum_{\ell=0}^{m(\theta, I_{1})}
C^{\ell+1} \binom{m(\theta, I_{1})}{\ell} p^{\ell} \|
(Q^{\frac{1}{2}})^{2\theta + \langle I_{1} \rangle - \ell
\frac{q}{q-1} } u \|_{0} 
\end{multline}
Since
$$ 
2 \theta + \langle I_{1} \rangle = 2 \theta + \langle I \rangle -
\frac{q}{q-1} = p - \frac{q}{q-1} ,
$$
we have that
$$
m(\theta, I_{1}) = \left[ \frac{q-1}{q} \left( 2 \theta + \langle I
    \rangle - \frac{q}{q-1}\right) \right]
= \left[ \frac{q-1}{q}
  p\right] -1 = m(\theta, I) -1.
$$
Hence we may conclude that
\begin{multline*}
\sum_{I_{1}} c_{I_{1}} \| Q^{\theta} X_{I_{1}} u \|_{0}
\\
\leq
M C_{1} \sum_{\ell=0}^{m(\theta, I) - 1}
C^{\ell+1} \binom{m(\theta, I) - 1}{\ell} p^{\ell+1} \|
(Q^{\frac{1}{2}})^{2\theta + \langle I \rangle - (\ell + 1)
  \frac{q}{q-1} } u \|_{0}
\\
=
M C_{1} \sum_{\ell=1}^{m(\theta, I)}
C^{\ell} \binom{m(\theta, I) - 1}{\ell-1} p^{\ell} \|
(Q^{\frac{1}{2}})^{2\theta + \langle I \rangle - \ell
  \frac{q}{q-1} } u \|_{0}
\\
=
M C_{1}C^{-1} \sum_{\ell=1}^{m(\theta, I)}
C^{\ell+1} \binom{m(\theta, I) - 1}{\ell-1} p^{\ell} \|
(Q^{\frac{1}{2}})^{2\theta + \langle I \rangle - \ell
  \frac{q}{q-1} } u \|_{0}
\\
\leq
\frac{1}{2}
\sum_{\ell=1}^{m(\theta, I)}
C^{\ell+1} \binom{m(\theta, I) - 1}{\ell-1} p^{\ell} \|
(Q^{\frac{1}{2}})^{2\theta + \langle I \rangle - \ell
  \frac{q}{q-1} } u \|_{0} ,
\end{multline*}
provided $ C $ is chosen so that  $ M C_{1}C^{-1} \leq \frac{1}{2} $. 

Finally consider $ A_{3} $.
Since
$$
2 \theta + 2 - \frac{q}{q-1} + \langle \tilde{I} \rangle = 2 \theta +
2 - \frac{q}{q-1} + \langle I \rangle - 2
= 2 \theta + \langle I
\rangle - \frac{q}{q-1} < p,
$$
we may apply the inductive hypothesis. Moreover as above
$$ 
\left[ \frac{q-1}{q} \left( p - \frac{q}{q-1} \right) \right] = \left[
  \frac{q-1}{q} p \right] - 1 = m(\theta,I) - 1 ,
$$
so that, renaming $ \tilde{C} $ the constant $ C $ in the definition
of $ A_{3} $, we have
\begin{multline*}
  \tilde{C} \| Q^{\theta + 1 - \frac{q}{2(q-1)}} X_{\tilde{I}} u \|_{0}
  \\
  \leq \tilde{C}
  \sum_{\ell=0}^{m(\theta, I)-1}
C^{\ell+1} \binom{m(\theta, I) - 1}{\ell} p^{\ell} \|
(Q^{\frac{1}{2}})^{p - (\ell + 1) \frac{q}{q-1} } u \|_{0}
\\
=
\tilde{C}
  \sum_{\ell=1}^{m(\theta, I)}
C^{\ell} \binom{m(\theta, I) - 1}{\ell - 1} p^{\ell - 1} \|
(Q^{\frac{1}{2}})^{p - \ell \frac{q}{q-1} } u \|_{0}
\\
\leq
\frac{1}{2}
\sum_{\ell=1}^{m(\theta, I)}
C^{\ell+1} \binom{m(\theta, I) - 1}{\ell - 1} p^{\ell} \|
(Q^{\frac{1}{2}})^{p - \ell \frac{q}{q-1} } u \|_{0} ,
\end{multline*}
provided $ C $ is chosen so that $ \tilde{C} C^{-1} \leq \frac{1}{2}
$. Hence we conclude that
\begin{multline*}
  \| Q^{\theta} X_{I} u \|_{0}
  \leq
  \sum_{\ell=0}^{m(\theta,I) -1} C^{\ell+1} \binom{m(\theta,I) -1}{\ell}
  p^{\ell} \|(Q^{\frac{1}{2}})^{p - \ell \frac{q}{q-1}} u \|_{0}
  \\
  +
  \sum_{\ell=1}^{m(\theta, I)}
C^{\ell+1} \binom{m(\theta, I) - 1}{\ell - 1} p^{\ell} \|
(Q^{\frac{1}{2}})^{p - \ell \frac{q}{q-1} } u \|_{0}
\\
=
\sum_{\ell=0}^{m(\theta, I)}
C^{\ell+1} \binom{m(\theta, I)}{\ell} p^{\ell} \|
(Q^{\frac{1}{2}})^{p - \ell \frac{q}{q-1} } u \|_{0} ,
\end{multline*}
thus concluding the proof of Proposition \ref{prop:ineqQ}.

If $ \theta > 1 $ the proof is completely analogous, using Lemma
\ref{lemma:A}. 
\end{proof}
\section{Appendix: Reducing Vector Fields to Powers of $ Q $ }

\setcounter{equation}{0}
\setcounter{theorem}{0}
\setcounter{proposition}{0}
\setcounter{lemma}{0}
\setcounter{corollary}{0}
\setcounter{definition}{0}

We first prove th following
\begin{proposition}
\label{prop:A}
Let $ \mu \in \Q^{+} $. Then there exists a positive constant, $ C $,
independent of $ \mu $, such that
\begin{equation}
\label{eq:A1D}
\| Q^{\mu} x \partial_{x} u \|_{0} \leq C \left( \| Q^{\mu +
    \frac{q}{2(q-1)}} u\|_{0} + \mu^{2\frac{q-1}{q} \mu + 1} \| u
  \|_{0} \right).
\end{equation}
\end{proposition}
\begin{proof}
Let $ [ \mu ] = k $, so that $ \mu = k + \theta $, $ 0 \leq \theta < 1
$. Then
$$
Q^{\mu} x \partial_{x} =  x \partial_{x} Q^{\mu} + [ Q^{\theta} Q^{k}
, x \partial_{x} ]
=
x \partial_{x} Q^{\mu} + Q^{\theta} [ Q^{k}, x \partial_{x} ] + [
Q^{\theta}, x \partial_{x} ] Q^{k} .
$$
Since $ x \partial_{x} $ has weight $ \frac{q}{q-1} $ we have
$$ 
\| x \partial_{x} v \|_{0} \leq C_{1} \| Q^{\frac{q}{2(q-1)}} v \|_{0},
$$
for a suitable constant $ C_{1} $. Here we used the pseudodifferential
calculus adapted to the anharmonic oscillator $ Q $ (see \cite{helffer84} and
\cite{bm-apde1}, Definition 2.1.)
Hence
$$ 
\| x \partial_{x} Q^{\mu} u \|_{0} \leq C_{1} \| Q^{\mu + \frac{q}{2(q-1)}} u \|_{0}.
$$
Analogously for the third term we have the estimate
$$ 
\| [ Q^{\theta}, x \partial_{x} ] Q^{k} u \|_{0} \leq C_{2} \| Q^{\mu}
u \|_{0} \leq C_{3} \| Q^{\mu + \frac{q}{2(q-1)}} u \|_{0} ,
$$
where $ C_{2} $, $ C_{3} $ are independent of $ k $.
Here we also used the fact that if $\sigma $, $ \tau $ are rational
numbers such that $ 0 < \sigma < \tau $, then $ \| Q^{\sigma} u \|_{0}
\leq C \| Q^{\tau} u \|_{0} $, for a suitable constant $ C > 0 $,
independent of $ \sigma $, $ \tau $.

Let us consider the term $ Q^{\theta} [ Q^{k}, x \partial_{x} ] $.
We have
$$ 
Q^{\theta} [ Q^{k}, x \partial_{x} ] = \sum_{i=1}^{k} \binom{k}{i}
Q^{\theta} (\ad Q)^{i}(x \partial_{x}) Q^{k-i} .
$$
The iterated commutator above is a sum of products of $ X_{-} $, $
X_{+} $ (see equation \eqref{eq:X+})
$$ 
(\ad Q)^{i}(x \partial_{x}) = \sum_{I} c_{I} X_{I} ,
$$
where the $ |c_{I}| \leq C_{4}^{i} $, the number of summands is
bounded by $ C_{5}^{i} $ and
$$ 
\langle I \rangle = i \frac{q-2}{q-1} + \frac{q}{q-1}.
$$
Applying Proposition \ref{prop:ineqQ} to $ Q^{\theta} X_{I} $ we
obtain the estimate
\begin{equation}
\label{eq:D2}
\| Q^{\theta} (\ad Q)^{i}(x \partial_{x}) v \|_{0}
\leq C_{*}^{i}
\sum_{\ell=0}^{m(i)} C^{\ell+1} \binom{m(i)}{\ell} p(i)^{\ell}
\|(Q^{\frac{1}{2}})^{p(i) - \ell \frac{q}{q-1}} v \|_{0} ,
\end{equation}
where
$$ 
p(i) = 2\theta + i \frac{q-2}{q-1} + \frac{q}{q-1} , \qquad m(i) =
\left[ 2 \theta \frac{q-1}{q} + i \frac{q-2}{q}\right] + 1.
$$
It will be useful to simplify the sums of the above type by, roughly,
taking the terms with the maximum and the minimum power of $ Q $, and,
correspondingly, with the minimum and maximum power of $ p(i) $. This
is done by a sort of convexity estimate of the following type
\begin{proposition}
\label{prop:cvx}
Let $ 1 \leq p, p' < + \infty $, $ p^{-1} + p'^{-1} = 1 $, and let $
\lambda > 0 $ be a real number. Then
$$ 
\| \lambda Q v \|_{0} \leq \frac{1}{p} \| Q^{p} v \|_{0} +
\frac{1}{p'} \lambda ^{p'}  \| v \|_{0} .
$$
\end{proposition}
The proof is straightforward, by using the spectral mapping theorem
and we omit it.

Going back to \eqref{eq:D2} we may write
\begin{multline}
  \label{eq:D3}
\| Q^{\theta} (\ad Q)^{i}(x \partial_{x}) v \|_{0}
\\
\leq C_{*}^{i} \sum_{\ell=0}^{m(i)} C^{\ell+1} \binom{m(i)}{\ell}
\left( \| (Q^{\frac{1}{2}})^{p(i)} v \|_{0} + p(i)^{p(i)
    \frac{q-1}{q}} \| v \|_{0}\right)
\\
\leq
C_{6}^{i+1} \left( \| (Q^{\frac{1}{2}})^{p(i)} v \|_{0} + p(i)^{p(i)
    \frac{q-1}{q}} \| v \|_{0}\right),
\end{multline}
due to the bound
$$ 
C_{*}^{i} \sum_{\ell=0}^{m(i)} C^{\ell+1} \binom{m(i)}{\ell} = C
C_{*}^{i} (C+1)^{m(i)} \leq C_{6}^{i+1} .
$$
As a consequence
\begin{multline}
  \label{eq:D4}
  \| Q^{\mu} x \partial_{x} u \|_{0}
  \leq
  C_{7} \| Q^{\mu + \frac{q}{2(q-1)}} u \|_{0}
\\
  + \sum_{i=1}^{k}
  C_{6}^{i+1} \binom{k}{i} \left( \| (Q^{\frac{1}{2}})^{p(i) + 2(k-i)}
    u \|_{0} + p(i)^{\frac{q-1}{q} p(i)} \| Q^{k-i} u \|_{0} \right)
  \\
  = C_{7} \| Q^{\mu + \frac{q}{2(q-1)}} u \|_{0}
  \\
  + \sum_{i=1}^{k} C_{6}^{i+1} \binom{k}{i}
  \left( \| Q^{\mu - (i-1)
      \frac{q}{2(q-1)}} u \|_{0}
    + p(i)^{\frac{q-2}{q} i + 1 + 2 \theta \frac{q-1}{q}} \| Q^{k-i} u \|_{0} \right)  
\end{multline}
Observe now that
$$ 
p(i)^{\frac{q-2}{q} i + 1 + 2 \theta \frac{q-1}{q}} \leq C_{0}'^{i}
i^{\frac{q-2}{q} i} i^{1+2\theta \frac{q-1}{q}}  \leq C_{0}^{i}
\mu^{\frac{q-2}{q} i},
$$
since $ i \geq 1 $, for a suitable constant $ C_{0} $. Moreover $
\binom{k}{i} \leq \mu^{i} i!^{-1} $, since $ k = [\mu] $. Applying
Proposition \ref{prop:cvx} to both terms under the sum sign above, we
obtain
$$ 
\mu^{i} \| Q^{\mu + \frac{q}{2(q-1)} - i \frac{q}{2(q-1)}} u \|_{0}
\leq \| Q^{\mu + \frac{q}{2(q-1)}} u \|_{0} + \mu^{2\frac{q-1}{q} \mu
  + 1} \| u \|_{0} ,
$$
and
$$ 
\mu^{\frac{q-2}{q} i + i} \|Q^{\mu + \frac{q}{2(q-1)} - i} u \|_{0}
\leq \| Q^{\mu + \frac{q}{2(q-1)}} u \|_{0} + \mu^{2\frac{q-1}{q} \mu
  + 1} \| u \|_{0} .
$$
Plugging the above estimates into \eqref{eq:D4}, we find
\begin{multline*}
 \| Q^{\mu} x \partial_{x} u \|_{0}
 \leq
C_{8} \sum_{i=0}^{[\mu]} \frac{C^{i+1}}{i!} \left ( \| Q^{\mu +
    \frac{q}{2(q-1)}} u \|_{0} + \mu^{2\frac{q-1}{q} \mu + 1} \| u
  \|_{0} \right)
\\
\leq
C_{9}
\left(
  \| Q^{\mu + \frac{q}{2(q-1)}} u \|_{0} + \mu^{2\frac{q-1}{q} \mu
  + 1} \| u \|_{0} \right) .
\end{multline*}
This completes the proof of the proposition.

\end{proof}
Proposition \ref{prop:A} takes care of the action of (rational) powers
of $ Q $ on the principal part of the transport operator $ P_{1} $. We
need a similar result for the other transport operators, $ P_{k} $, $
k = 2,\ldots, 2a $, in \eqref{eq:Pj}.
\begin{proposition}
\label{prop:AjD}
Let $ \mu \in \Q^{+} $. Then there exists a positive constant, $ C $,
independent of $ \mu $, such that
\begin{equation}
\label{eq:AjD}
\| Q^{\mu} (x \partial_{x})^{j} u \|_{0} \leq C_{1}^{j} \left( \| Q^{\mu +
   j \frac{q}{2(q-1)}} u\|_{0} + \mu^{2\frac{q-1}{q} \mu + j} \| u
  \|_{0} \right) ,
\end{equation}
$ j = 1,\ldots, 2a $.
\end{proposition}
\begin{proof}
We proceed by induction with respect to $ j $. When $ j=1 $ the
assertion is just Proposition \ref{prop:A}. Assume that the assertion
for $ j $ holds and let us prove the assertion for $ j+1 $.
We have
$$ 
\| Q^{\mu} (x \partial_{x})^{j+1} u \|_{0} = \| Q^{\mu} (x
\partial_{x})^{j} (x \partial_{x}) u \|_{0} .
$$
Applying the above estimate for $ j $ we obtain
$$
  \| Q^{\mu} (x \partial_{x})^{j+1} u \|_{0}
  \leq
  C_{1}^{j} \left( \| Q^{\mu +
   j \frac{q}{2(q-1)}} (x \partial_{x}) u\|_{0} + \mu^{2\frac{q-1}{q}
   \mu + j} \| (x \partial_{x}) u \|_{0} \right) 
$$
For the first term we use Proposition \ref{prop:A}, while for the
second we apply Propositions \ref{prop:ineqQ} and \ref{prop:cvx}.
Now
\begin{multline*}
\| Q^{\mu} (x \partial_{x})^{j+1} u \|_{0}
\\
\leq
C_{1}^{j} \Bigg (
 C_{*}
\left(
  \| Q^{\mu + (j+1)\frac{q}{2(q-1)}} u \|_{0} + \mu^{2\frac{q-1}{q} \mu
    + 1 + j} \| u \|_{0} \right)
\\
+ C_{**} \left( \| Q^{\mu +
   (j + 1) \frac{q}{2(q-1)}} u\|_{0} + \mu^{2\frac{q-1}{q} \mu + j + 1} \| u
 \|_{0} \right)
\Bigg),
\end{multline*}
so that choosing $ C_{1} \geq 2 \max\{ C_{*}, C_{**}\} $ achieves the proof.
\end{proof}

\end{document}